\documentclass[11pt]{article}

\usepackage[english]{babel}
\usepackage{amssymb}
\usepackage{amsfonts}
\usepackage{cite}

\usepackage{bbold}
\renewcommand{\baselinestretch}{1.3}
\makeatletter
\def\@begintheorem#1#2{\trivlist%
 \item[\hskip \labelsep{\sffamily\bfseries #2\ #1}]\itshape}
\newtheorem{teo}{Theorem}[section]
\newtheorem{defi}[teo]{Definition}
\newtheorem{cor}[teo]{Corollary}
\newtheorem{lem}[teo]{Lemma}
\newtheorem{pro}[teo]{Proposition}
\newtheorem{_rem}[teo]{Remark}
\newtheorem{_rems}[teo]{Remarks}
\newtheorem{_eje}[teo]{Example}

\newenvironment{rem}{\def\@begintheorem##1##2{\trivlist%
 \item[\hskip\labelsep{\sffamily\bfseries ##2\ ##1}]}\begin{_rem}}{\end{_rem}}
\newenvironment{rems}{\def\@begintheorem##1##2{\trivlist%
 \item[\hskip\labelsep{\sffamily\bfseries ##2\ ##1}]}
\begin{_rems}}{\end{_rems}}

\makeatother

\newenvironment{beweis}{{\em Proof:}}{\hfill $\rule{2mm}{2mm}$
\vspace{3mm}}

\DeclareMathAlphabet{\Ma}{U}{msa}{m}{n}
\DeclareMathAlphabet{\Mb}{U}{msb}{m}{n}
\DeclareMathAlphabet{\Meuf}{U}{euf}{m}{n}

\def\got#1{\Meuf{#1}}

\DeclareSymbolFont{ASMa}{U}{msa}{m}{n}
\DeclareSymbolFont{ASMb}{U}{msb}{m}{n}
\DeclareMathSymbol{\hrist}{\mathord}{ASMa}{"16}
\DeclareMathSymbol{\varkappa}{\mathalpha}{ASMb}{"7B}
\DeclareMathSymbol{\CrPr}{\mathord}{ASMb}{"6F}

\def\restriction{\upharpoonright}

  \def\al #1.{{\mathcal{#1}}}
  \def\ot #1.{{\got{#1}}}
 
\def\CCRX{\overline{\Delta( X,\,\sigma)}}
  \def\ccr #1,#2.{\overline{\Delta(#1,\,#2)}}

  \def\b #1.{{\bf #1}}
  \def\cross#1.{\mathrel{\mathop{\times}\limits_{#1}}}
                  
  \def\C{\Mb{C}}
  \def\N{\Mb{N}}
  
  \def\R{\Mb{R}}
  \def\Z{\Mb{Z}}
 \def\un{{\mathbb 1}}

\def\f #1,#2.{\mathsurround=0pt \hbox{${#1\over #2}$}\mathsurround=5pt}

\def\hlf{{\f 1,2.}}

  \def\wt{\widetilde}

  \def\cross #1.{\mathrel{\raise 3pt\hbox{$\mathop\times\limits_{#1}$}}}
  \def\ol #1.{\overline{#1}}
\def\b #1.{{\bf #1}}
\def\slim{\mathop{\hbox{\rm s-lim}}}
\def\ker{{\rm Ker}\,}
\def\aut{{\rm Aut}\,}
\def\dom{{\rm Dom}\,}

\def\ran{{\rm Ran}\,}
\def\rep{{\rm Rep}\,}

\def\Ad{{\rm Ad}\,}
\def\rlf{{R(\lambda,f)}}
\def\rsl{\mathord{\al R.(X,\sigma)}}
\def\s #1.{_{\smash{\lower2pt\hbox{\mathsurround=0pt $\scriptstyle #1$}}\mathsurround=5pt}}
\def\set #1,#2.{\left\{\,#1\;\bigm|\;#2\,\right\}}
\def\maprightu #1;{\smash{\mathop{\longrightarrow}\limits^{#1}}}
\def\maprightd #1;{\smash{\mathop{\longrightarrow}\limits_{#1}}}
\def\maprightt #1,#2.{\mathrel{\smash{\mathop{\longrightarrow}\limits_{#1}^{#2}}}}
\def\chop{\hfill\break}
\def\j{\phi}

\def\XP#1!{\renewcommand{\baselinestretch}{.7}\marginpar{$\leftarrow${\footnotesize #1}\hfil}
 \renewcommand{\baselinestretch}{1}}
\def\XB{\marginpar{
{\footnotesize Change~starts-}\lower 11pt\hbox{\mathsurround=0pt$
\!\!\displaystyle{
\Bigg\downarrow}$\mathsurround=3pt}}}
\def\XE{\marginpar{{\footnotesize Change~ends-}\raise 10pt\hbox{\mathsurround=0pt$
\!\!\displaystyle{
\Bigg\uparrow}$\mathsurround=3pt}}}
\def\f #1,#2.{\mathsurround=0pt \hbox{${#1\over #2}$}\mathsurround=5pt}

\def\ie{\textit{i.e.\ }}
\def\eg{\textit{e.g.\ }}
\def\etc{\textit{etc\ }}
\def\margin #1.{\marginpar{#1}}

\def\WD{{\cal O}}
\def\Ob{{\cal P}}

\title{\bf The Resolvent Algebra:  \\ 
A New Approach to Canonical Quantum Systems}

\author{
 {\sc Detlev Buchholz}            \\[1mm]
{\footnotesize Universit\"at G\"ottingen} \\
\and
 {\sc Hendrik Grundling}                                            \\[1mm]
{\footnotesize University of New South Wales} 
}
\date{}

\begin{document}
\maketitle
\begin{abstract} 
\noindent The standard C*--algebraic version of the 
algebra of canonical commutation
relations, the Weyl algebra, frequently causes difficulties in
applications since it neither admits the formulation
of physically interesting dynamical laws nor does it incorporate pertinent 
physical observables such as (bounded functions of) the Hamiltonian.
Here a novel C*--algebra of the canonical 
commutation relations is presented which does not suffer from 
such problems.
It is based on the resolvents of the canonical operators 
and their algebraic relations. The resulting 
C*--algebra, the resolvent algebra, is shown to have many  
desirable analytic properties and the regularity
structure of its representations
is surprisingly simple. Moreover, the resolvent algebra is 
a convenient framework for applications to interacting and
to constrained quantum
systems, as we demonstrate by
several examples.
\end{abstract}

\bigskip\bigskip
\par\noindent\hrule\par

\section{Introduction}

Canonical systems of operators have always been a central ingredient  
in the modelling of quantum systems. There is an extensive
literature analyzing their properties, starting 
with the seminal paper of Born, Jordan and Heisenberg on the 
physical foundations and reaching a first 
mathematical satisfactory formulation in the works of von Neumann
and of Weyl. 
These canonical systems of operators may all be presented in the following
general form: there is
a real linear map $\j$  from a given symplectic space 
${(X,\,\sigma)}$ to a linear space of essentially selfadjoint operators 
on some common dense invariant core $\al D.$ in a Hilbert space $\al H.,$
satisfying the relations
\[
\big[\j(f),\,\j(g)\big]=i\sigma(f,\, g) \, {\un}, \quad \j(f)^* =\j(f)\quad\hbox{on}\quad
\al D. \, .
\]
In the case that $X$ is finite dimensional, we can reinterpret this 
relation in terms of the familiar quantum mechanical 
position and momentum operators, and if 
$X$ consists of Schwartz functions on a space--time 
manifold one may consider $\j$ to be a bosonic quantum field.
 The observables of the system are then constructed from the operators
 ${\{\j(f)\,\mid\,f\in X\}}\,,$ usually as polynomial expressions.
 Since one wants to study representations of such
systems, some care needs to be taken about the appropriate mathematical
framework to do this in, especially since there are known pathologies
\eg for the case that $X$ is infinite dimensional. 
 Here we will use C*-algebras to encode the algebraic
information of the canonical systems, given the rich source of mathematical
tools available.

As is well--known, if ${(X,\,\sigma)}$ is non--degenerate
then the operators $\j(f)$ cannot all be bounded.  
Thus starting from the polynomial 
algebra $\al P.$ generated by ${\{\j(f)\,\big|\,f\in X\}}$ 
we have to obtain a C*--algebra 
encoding the same algebraic information, necessarily in bounded form. 
The obvious way to take this step 
is to form suitable bounded functions of the fields $\j(f)$.
In the approach introduced by Weyl, this is done by 
considering the C*--algebra generated by the set of unitaries
\[
\big\{\exp\big(i\j(f)\big)\,\big|\,f\in X\big\}
\]
and this C*--algebra is simple.  It can be defined abstractly
as the C*--algebra generated by a set of unitaries  
${\set\delta_f , {f\in X}.}$ subject to the relations 
$\delta_f^*=\delta_{-f}$ and $\delta_f\delta_g=e^{-i\sigma(f,g)/2}
\delta_{f+g}$. 
This is the familiar Weyl (or CCR) algebra, often denoted $\CCRX$
\cite{Man}. 
By its definition, it has a representation in which 
the unitaries $\delta_f$ can be identified with
the exponentials $e^{i \j(f)}$, and hence we can obtain the 
concrete algebra $\al P.$ back from these.
Such representations $\pi:\CCRX\to\al B.(\al H.)$, \ie those for which
the one--parameter groups $\lambda\to\pi(\delta\s\lambda f.)$ are strong operator
continuous for all $f\in X$ are called regular.
Since for physical situations the quantum fields are defined as the generators
of the one--parameter groups $\lambda\to\pi(\delta\s\lambda f.),$ the
representations of interest are required to be regular. 

The Weyl algebra suffers from several well--known flaws. First and foremost,
the dynamics (one--parameter automorphism groups) of the Weyl algebra which 
are most naturally defined from symplectic transformations of $X$,  correspond
to the dynamics produced by quadratic Hamiltonians, and this excludes most physically 
interesting situations. As a matter of fact, 
in the case $X=\R^2$ in the Schr\"odinger representation
for $\CCRX$, the time evolutions obtained from Hamiltonians of the form
$H=P^2+V(Q)$ for potentials $V\in L^\infty(\R)\cap L^1(\R)$ do not preserve
$\CCRX$ unless $V=0$ (cf.~\cite{FaVB}). Thus the Weyl algebra does not allow
the definition of much interesting dynamics on it, and this limits its usefulness.
Second, in regular representations, natural observables
such as bounded functions of the Hamiltonian 
are not in $\CCRX$. Third, the Weyl algebra has a large number of nonregular 
representations \cite{Ace, Gr3, GrHu3}. Through the use of w*--limits
of regular states these are 
interpreted as situations where the field $\j(f)$ can have ``infinite field strength''.
Whilst this is useful for some idealizations e.g. plane waves cf.~\cite{Ace}
or for quantum constraints cf.~\cite{GrHu3}, for most physical situations
one wants to exclude such representations.

This motivates the consideration of alternative
versions of the C*--algebra of canonical commutation relations. 
Instead of taking the C*--algebra generated 
by exponentials of the underlying fields,
as for the Weyl algebra,
we will consider the C*--algebra generated by
the resolvents of the fields. These are formally given by 
$R(\lambda, f) := (i \lambda \un - \j(f))^{-1}$. 
All algebraic 
properties of the fields can be expressed in terms of 
simple relations amongst these resolvents. The unital C*-algebra
generated by the resolvents, 
\[
\al R. := {\rm C}^*\big\{ \, R(\lambda, f)  \mid\, 
f\in X, \, \lambda \in \R\backslash 0 \big\}\,,
\]
called the resolvent algebra, will be the subject of our investigation,
and below we will exhibit its main algebraic and analytic properties.

As a preview, consider for fixed $f \in X$ the abelian subalgebra
${\al R.}_f := {\rm C}^*\big\{\, R(\lambda, f) \mid\,
\lambda \in  \R \big\} $. This algebra is
isomorphic to the algebra of continuous functions on 
$\R$, vanishing at infinity. Thus all states of ``infinite field
strength'' for $\j(f)$ annihilate this algebra, which greatly simplifies
the representation theory of $\al R.$. The obvious price to pay 
is that $\al R.$ is not simple; but this non--trivial ideal 
structure turns out to be 
useful in applications. 
Moreover, in contrast to the Weyl algebra, the 
underlying unbounded field
operators are affiliated with the C*--algebra $\al R.$ in the sense of 
Damak and Georgescu \cite{Georg}. 
This feature allows one to
manipulate rigorously polynomial expressions of the fields 
by the use of ``mollifiers'' (explained below). As a matter
of fact, the latter observation motivated us to introduce the 
resolvent algebra in our recent study of supersymmetry \cite{BuGr}, 
where we needed these mollifiers to construct superderivations.
Thus, in many respects, the resolvent algebra provides a 
technically convenient framework for the study of 
canonical quantum systems.

The structure of this paper is as follows. In Sect.~\ref{Moll} we introduce
the notion of mollifiers to encode polynomials in
the bosonic fields in bounded form,
and we show that the Weyl algebra does not have any mollifiers, whereas
the resolvent algebra does. In Sect.~\ref{ResBasics} we define the resolvent
algebra abstractly as a C*--algebra, and establish basic algebraic properties.
In Sect.~\ref{StatesRep} we study states and representations of the resolvent
algebra, in particular we introduce the notion of a ``regular" representation
and show that these are in a natural bijection with the usual regular 
representations 
of the Weyl algebra. We obtain interesting decompositions of the 
symplectic space
associated to representations, and we prove that every regular representation
is faithful.    
In Sect.~\ref{FurStruc} we consider further algebraic properties
and find that the resolvent algebra is nonseparable, find a tensor product
structure for some of its subalgebras, and for the case of a finite
dimensional symplectic
space, find that it contains a copy of the algebra of compact operators.
 In Sect.~\ref{DynHam} we consider how to encode dynamics and Hamiltonians
 in the resolvent algebra. Many resolvents of Hamiltonians are already
 in the resolvent algebra, and using the copies of the compacts
in the resolvent algebra, we can encode many more dynamical systems
than for the Weyl algebra. We then produce two applications to illustrate
the usefulness of the resolvent algebra. In Sect.~\ref{InteractSect} we
develop a model of an infinite family of atoms  which
are confined around the points of a lattice by a pinning potential 
and interact with their nearest neighbours.
We construct a ground state for it by algebraic means.
 In Sect.~\ref{Qconstr} we take a brief look
at Dirac constraints theory for linear bosonic constraints in the
context of the resolvent algebra, and find that it is considerably
simpler than in the Weyl algebra. All proofs for our results are collected
in  Sect.~\ref{Proofs}.

\section{Mollifiers and Resolvents}
\label{Moll}

There are several concepts of when a selfadjoint operator $A$ on 
a Hilbert space $\al H.$
is {\it affiliated} with a concretely represented
C*-algebra $\al A.\subset\al B.(\al H.)$. 
One is that
the resolvent ${(i\lambda\un-A)^{-1}}\in\al A.$ for some 
$\lambda\in\R\backslash0$ (hence for all $\lambda\in\R\backslash0).$ 
This notion is used by Damak and Georgescu~\cite{Georg} 
(and is weaker than the one  used by  Woronowicz~\cite{Wor1}) 
and it implies the usual one, \ie that $A$ commutes with all unitaries
commuting with $\al A.$ (but not conversely).
Observe that then 
\[
A(i\lambda\un-A)^{-1}=\ol(i\lambda\un-A)^{-1}A.=i\lambda(i\lambda\un-A)^{-1}-\un
\in\al A.\;.
\]
Thus the resolvent ${(i\lambda\un-A)^{-1}}=M$ acts as a ``mollifier''
for $A,$ \ie $\ol MA.$ and $AM$ are bounded and in $\al A.,$ and $M$
is invertible such that $M^{-1}MA=A=AMM^{-1}.$
This suggests that as $AM$ and $\ol MA.$ in $\al A.$ carries the
information of $A$ in bounded form, 
we can ``forget'' the original representation, and study the affiliated
$A$ abstractly through these elements.
In the literature, the bounded operators $A_\lambda:=i\lambda A\,{(i\lambda\un-A)^{-1}}$
are called ``Yosida approximations'' of $A$ \cite[p 9]{Paz}.

We want to apply this idea to a 
bosonic field as above,
\ie for a fixed Hilbert space $\al H.$ we assume that there is
a common dense invariant core $\al D.\subset\al H.$  for the 
selfadjoint operators $\phi(f)$, $f\in X$ 
on which the $\phi(f)$ satisfy the canonical commutation relations.
One may be tempted to find mollifiers for the operators $\phi (f)$ in the 
(concretely represented) Weyl algebra
\[
\CCRX=C^*\set\exp(i\phi (f)),{f\in X}.\subset\al B.(\al H.)\,,
\]
but unfortunately this is not possible because \cite{BuGr}:
\begin{pro}
\label{CCRX}
The Weyl algebra $\CCRX$ contains no 
nonzero element $M$ such that $\phi (f)M$ is bounded 
for some $f\in X\backslash0$.
Thus $\CCRX$ contains no mollifier for any nonzero
$\phi (f),$  and $\phi (f)$ is not affiliated with $\CCRX$.
\end{pro}

Our solution is to abandon the Weyl algebra as the appropriate
C*-algebra to model the bosonic fields $\phi (f)$, and instead to
choose the unital C*-algebra generated by the resolvents, 
\[
C^*\set{R(\lambda,f)},\lambda\in\R\backslash 0,\;
f\in{X}.  \, ,
\]
where $R(\lambda,f):=(i\lambda\un-\phi (f))^{-1}$. Then by construction
all $\phi (f)$ are
 affiliated to this C*-algebra and it contains
mollifiers $R(\lambda,f)$ for all of them.

\section{Resolvent Algebra - Basics}
\label{ResBasics}

The above discussion took place in a concrete setting, 
\ie represented on a Hilbert space, and we would like to abstract this.
Just as the Weyl algebra can be abstractly defined by the Weyl relations,
we want to abstractly define the C*-algebra of resolvents
by its generators and relations.
\begin{defi}
\label{ResAlg}
Given a symplectic space $(X,\,\sigma),$ we define 
$\al R._0$ to be the universal unital 
*-algebra generated by
the set $\set{R(\lambda,f)},\lambda\in\R\backslash 0,\;f\in X.$
and the relations
\begin{eqnarray}
\label{Riden}
R(\lambda,0)&=&-{i\over\lambda}\,\un \\[1mm]
\label{Rinvol}
R(\lambda,f)^*&=&R(-\lambda,f) \\[1mm]
\label{Rhomog}
\nu \,R(\nu \lambda,\, \nu f)   &=&  \rlf \\[1mm]
\label{Resolv}
\rlf - R(\mu,f) &=& i(\mu-\lambda)\rlf R(\mu,f)  \\[1mm]
\label{Rccr}
\big[\rlf,\,R(\mu,g)\big] &=&
i\sigma(f,g)\,\rlf\,R(\mu,g)^2\rlf \\[1mm]
\label{Rsum}
\rlf R(\mu,g)&=& R(\lambda+\mu,\,f+g)[\rlf+R(\mu,g)
+i\sigma(f,g)\rlf^2R(\mu,g)]
\end{eqnarray}
where $\lambda,\, \mu, \, \nu \in\R\backslash 0$ and $f,\,g\in X$,
and for (\ref{Rsum}) we require $\lambda+\mu\not=0$.
That is, start with the free unital *-algebra generated by
 $\set{R(\lambda,f)},\lambda\in\R\backslash 0,\;f\in X.$
 and factor out by the ideal generated by the relations
(\ref{Riden}) to (\ref{Rsum}) to obtain the *-algebra $\al R._0$. 
\end{defi}
\begin{rems}
(a) \ The *-algebra $\al R._0$ is nontrivial, because it has nontrivial
representations. {}For instance, in a Fock representation 
$\pi$ of the
CCRs over $(X,\,\sigma)$ 
one has the selfadjoint CCR-fields $\phi_\pi (f)$, $f \in X$ 
from which one can define $\pi(\rlf) :=(i\lambda\un-\phi_\pi (f))^{-1}$
to obtain a representation of $\al R._0.$ \\[1mm]
(b) \ Obviously (\ref{Rinvol}) encodes the selfadjointness of 
$\phi_\pi (f),$ (\ref{Rhomog}) encodes $\phi_\pi (\nu f)
=\nu \phi_\pi (f),$
(\ref{Resolv}) encodes that $\rlf$ is a resolvent,
(\ref{Rccr}) encodes the canonical commutation relations
and (\ref{Rsum}) encodes additivity
$\phi_\pi (f+g)=\phi_\pi (f)+\phi_\pi (g)$. 
Note that $R(0,f)$ is undefined. \\[1mm]
(c) \ Let $\mu=-\lambda$ in Equation~(\ref{Resolv}) to get 
the useful equation
\begin{equation}
\label{RRstar}
 \rlf-\rlf^*=-2i\lambda\rlf\rlf^*\,.
\end{equation}
\end{rems}
To define our resolvent C*-algebra, we need to decide on which 
C*-seminorm to define on $\al R._0.$
The obvious choice is the enveloping
C*-seminorm, however since  $\al R._0$ is merely a *-algebra, we need to 
establish some uniform boundedness for its Hilbert space representations
before we can define this.
\begin{pro}
\label{Rbound}
Given a symplectic space $(X,\,\sigma),$ define 
$\al R._0$ as above.
\begin{itemize}
\item[(i)]
Let $\pi_0:\al R._0\to\al B.(\al H.)$ be a *-representation of $\al R._0,$ where 
$\al H.$ is a Hilbert space. Then ${\big\|\pi_0\big(\rlf\big)\big\|}\leq{|\lambda|^{-1}}.$
Thus, for each 
$A\in\al R._0$ there is a $c_A \geq 0$ such that 
${\|\pi(A)\|}\leq c_A$ for all (bounded) Hilbert space 
representations $\pi$ of $\al R._0.$
\item[(ii)] Let $\omega$ be a positive functional of $\al R._0,$ \ie
$\omega:\al R._0\to\C$ is linear and $\omega(A^*A)\geq 0$ for all $A\in\al R._0.$
Then the GNS--construction yields a cyclic *--representation of $\al R._0,$
denoted $\pi_\omega,$ consisting of {\it bounded} Hilbert space operators.
\end{itemize}
\end{pro}
Thus by (i)  we can form direct sum representations over infinite sets of representations
and still maintain boundedness. 
However the class of (nondegenerate) representations is not a set\footnote{
If the nondegenerate representations were a set, we could take the direct sum of
the representations which do not have themselves as a direct summand, to obtain Russell's paradox.},
 so we need to
find a suitable set of representations.
Let $\ot S.$ denote the set of positive functionals $\omega$ of $\al R._0$ for which 
$\omega(\un)=1,$ then by Proposition~\ref{Rbound}(ii) their GNS--representations are bounded,
and in fact by (i) they are uniformly bounded.
Since $\ot S.$ is a set, we can now sensibly define:
\begin{defi}
The {universal representation} 
$\pi_u:\al R._0\to\al B.(\al H._u)$ is given by
\[
\pi_u(A):=\bigoplus\big\{\pi_\omega(A)\,\big|\,\omega\in\ot S.\big\}\quad\hbox{and}\quad
\|A\|_u:=\|\pi_u(A)\|=\sup_{\omega\in\ot S.}\|\pi_\omega(A)\|
\]
denotes the {enveloping C*--seminorm} of $\al R._0$; 
note that $\|A\|_u =\sup\limits_{\omega\in\ot S.}\omega(A^*A)^{1/2}$ since
$\ot S.$ contains all vector states of all its GNS--representations. 
We define our {resolvent algebra} $\al R.(X,\,\sigma)$ as the  abstract
C*-algebra generated by $\pi_u(\al R._0),$ \ie we factor $\al R._0$ by
$\ker\pi_u$ and complete w.r.t.\ the enveloping C*--seminorm $\|\cdot\|_u.$
\end{defi}
\begin{rem}
Previously, in \cite{BuGr} we defined the resolvent algebra with a different
C*--seminorm, because we needed more analytic structure. 
Henceforth we will deal with the resolvent
algebra  $\al R.(X,\,\sigma)$ as defined above 
and denote its norm  by $\| \cdot \|$. 
Below we will prove
isomorphism with the previous version.
\end{rem}

We state some elementary properties of $\al R.(X,\,\sigma)$. 

\begin{teo}
\label{Relemen}
Let $(X,\,\sigma)$ be a given nondegenerate symplectic space, and define
$\al R.(X,\,\sigma)$ as above. Then for all $\lambda,\,\mu\in\R\backslash 0$
and $f,\, g\in X$ we have:
\begin{itemize}
\item[(i)] $[\rlf,\,R(\mu,f)]=0$. 
 Substitute $\mu=-\lambda$ to see that
$\rlf$ is normal.
\item[(ii)] $\rlf R(\mu,g)^2\rlf= R(\mu,g)\rlf^2 R(\mu,g)$.
\item[(iii)] $\big\| \rlf \big\| =|\lambda|^{-1}$.
\item[(iv)] $\rlf$ is analytic in $\lambda$. Explicitly, the series 
expansion (von Neumann series) 
\[
\rlf=\sum_{n=0}^\infty (\lambda_0-\lambda)^n\, i^n R(\lambda_0,f)^{n+1},
\qquad\lambda,\;\lambda_0\not=0\qquad\qquad
\]
converges absolutely in norm whenever $|\lambda_0-\lambda|<|\lambda_0|\;.$
\item[(v)] Let $T\in{\rm Sp}(X,\sigma)$ be a symplectic transformation.
Then $\alpha\big(\rlf\big):=R(\lambda,Tf)$ extends to an
automorphism $\alpha\in\aut\al R.(X,\,\sigma)$.
\end{itemize}
\end{teo}
Note that the von Neumann series for $\rlf$ converges for any
$\lambda\in\C$ with $|\lambda-\lambda_0|<|\lambda_0|,$ \ie on a disk
which stays off the imaginary axis. Using different $\lambda_0\in\R\backslash 0$ 
we can thus define ${R(z,f)}$ for any complex $z$ not on the imaginary axis,
\ie analytically continue $\rlf$ from $\R\backslash 0$ to $\C\backslash{i\R}.$
Thus $\al R.(X,\,\sigma)$
contains  also `resolvents' $R(z,f)$ for complex
$z\in\C\backslash{i\R}.$ In fact, we can define $\al R.(X,\,\sigma)$
as the enveloping C*--algebra of
the universal unital 
*-algebra generated by
the set $\set{R(z,f)},z\in\C\backslash{i\R},\;f\in X.$ and the analytic continuations 
of the relations~(\ref{Riden})--(\ref{Rsum}) \ie
\begin{eqnarray}
\label{Ridenz}
R(z,0)&=&-{i\over z}\,\un \\[1mm]
\label{Rinvolz}
R(z,f)^*&=&R(-\overline{z},f) \\[1mm]
\label{Rhomogz}
\nu R(\nu z,\nu f)
 &=& R(z,\,f)\;,\quad\nu \in\R\backslash 0  \\[1mm]
\label{Resolvz}
R(z,f)- R(w,f) &=& i(z-w)R(z,f)\,R(w,f)  \\[1mm]
\label{Rccrz}
\big[R(z,f),\,R(w,g)\big] &=&
i\sigma(f,g)\,R(z,f)\,R(w,g)^2R(z,f) \\[1mm]
\label{Rsumz}
R(z,f) R(w,g)&=& R(z+w,\,f+g)[R(z,f)+R(w,g)
+i\sigma(f,g)R(z,f)^2R(w,g)]
\end{eqnarray}
where $z,\, w\in\C\backslash{i\R}$ and $f,\,g\in X$,
and for (\ref{Rsumz}) we require $z+w\not\in{i\R}$.
Note that Eq.~(\ref{RRstar}) becomes
\begin{equation}
\label{RRstarz}
 R(z,f)-R(z,f)^*=i(z+\overline{z})R(z,f)R(z,f)^*\,,
\end{equation}
hence we get $\|R(z,f)\|=|{\rm Re}\,z|^{-1}.$
Since by Eq.~(\ref{Resolvz}) we can again (via a von~Neumann series)
prove analyticity off the imaginary axis, it follows that this C*--algebra
coincides with the previously defined $\al R.(X,\,\sigma)$.
Making use of this fact, we can exhibit another family of automorphisms
of the resolvent algebra which will be useful for proofs below.
\begin{pro}
\label{ImAut}
Let $(X,\,\sigma)$ be a given nondegenerate symplectic space. 
Then for each linear map $h:X\to\R$ there is an automorphism
$\beta_h\in\aut\al R.(X,\,\sigma)$ defined by
\[
\beta_h\big(R(z,f)\big) =
R(z+ih(f),\,f)
\]
for $f\in X$ and $z\in\C\backslash{i\R}.$
\end{pro}

As already stated, the resolvent
algebra is not simple, in contrast to the Weyl algebra. More specifically:
\begin{teo}
\label{Ideals0}
Let $(X,\,\sigma)$ be a given nondegenerate symplectic space. 
Then for all $\lambda\in\R\backslash 0$
and $f\in X\backslash 0$ we have that the
 closed two--sided ideal generated by $\rlf$ in 
$\al R.(X,\,\sigma)$ is
\[
\big[\rlf\rsl\big]=\big[\rsl\rlf\big]=\big[\rsl\rlf\rsl\big]
\]
where ${[\, \cdot \,]}$ indicates the closed linear span 
of its argument. This ideal is proper.
Moreover the intersection of the ideals 
$\big[R(\lambda_i,f_i)\rsl\big]$,
 $i=1,\ldots,\,n$ for distinct $f_i\in X\backslash 0$ is the ideal
${\big[R(\lambda_1,f_1)\cdots R(\lambda_n,f_n)\rsl\big]}$.
\end{teo}
{}From these ideals we can build other ideals, \eg for a set $S\subseteq X$
we can define the ideals ${\mathop{\bigcap}\limits_{f\in S}\big[\rsl R(\lambda,f)\big]}$
as well as ${\Big[\mathop{\bigcup}\limits_{f\in S}\big[\rsl R(\lambda,f)\big]\Big]}$. 
Ideals of a different structure will occur in the following sections.
Thus $\rsl$ has a very rich ideal structure.

\vspace*{1mm}
The fact that $\rsl$ is not simple, does not disqualify it from being used
as the observable algebra of a physical system, because we will show below that
its images in all physical (\ie regular) representations are isomorphic. 

\section{States, representations and regularity}
\label{StatesRep}

Any operator family $R_{\lambda}$, $\lambda \in \R \backslash 0$ 
on some Hilbert space satisfying the resolvent equation
(\ref{Resolv}) is called by Hille a pseudo-resolvent
and for such a family we know (cf.\ Theorem~1
in~\cite[p 216]{Yos}) that:
\begin{itemize}
\item{} All $R_\lambda$ have a common range and a common null space.
\item{} A pseudo-resolvent is the resolvent for an operator
$B$ iff $\ker R_{\lambda} =\{0\}$ for some (hence for
all) $\lambda \in \R \backslash 0$, and in this case
$\dom B=\ran R_\lambda$ for all $\lambda \in \R \backslash 0$.
\end{itemize}
This leads us to an examination of $\ker\pi\big(\rlf\big)$ for representations
$\pi$.
\begin{teo}
\label{Ideals1}
Let $(X,\,\sigma)$ be a given nondegenerate symplectic space, and define
$\al R.(X,\,\sigma)$ as above. Then for $\lambda\in\R\backslash 0$
and $f\in X\backslash 0$ we have:
\begin{itemize}
\item[(i)] If for a representation $\pi$ of $\rsl$ we have
$\ker\pi\big(\rlf\big)\not=\{0\} $,  
then $\ker\pi\big(\rlf\big)$ reduces $\pi(\rsl)$. Hence
there is a unique orthogonal decomposition
$\pi=\pi_1\oplus\pi_2$ such that $\pi_1(\rlf)=0$ and
$\pi_2(\rlf)$ is invertible.
\item[(ii)] Let $\pi$ be any nondegenerate representation of $\rsl$,
then  $$P_f:=\slim_{\lambda\to\infty}i\lambda\, \pi\big(\rlf\big) $$
exists, defines a central projection of $\pi(\rsl)''\!$, and it is the
range projection of $\pi\big(\rlf\big)$ as well as the projection
of the ideal ${\pi\left(\big[\rsl R(\lambda,f)\big]\right)}$.
\item[(iii)] If $\pi$ is a factorial representation of $\rsl$,
then $P_f=0$ or $\un$ and such $\pi$ are classified by the sets
${\big\{f\in X\backslash 0\,\mid\, P_f=\un\big\}}$.
\item[(iv)] There is a state $\omega\in\ot S.\big(\rsl\big)$
such that $\rlf\in\ker\omega$. Moreover, given a state $\omega$
with  $\rlf\in\ker\omega$, then $\rlf\in\ker\pi_\omega$. 
\end{itemize}
\end{teo}
Given a $\pi\in\rep\big(\al R.(X,\,\sigma),\al H._\pi\big)$ with
$\ker\pi\big(R(1,f)\big)=\{0\},$ we can define a field
operator by
\[
\j_\pi(f):=i\un-\pi\big(R(1,f)\big)^{-1}\;
\]
with domain $\dom \j_\pi(f)=\ran\pi\big(R(1,f)\big)$, and it has the following
properties:
\begin{teo} 
\label{RegThm}
Let $\al R.(X,\,\sigma)$ be as above, and let $\pi\in\rep\big(\al R.(X,\,\sigma),\al H._\pi\big)$
satisfy $\ker\pi\big(R(1,f)\big)=\{0\}= \ker\pi\big(R(1,h)\big)$ for given
$f,\; h\in X.$ Then
\begin{itemize}
\item[(i)]  $\j_\pi(f)$ is selfadjoint, and $\pi(\rlf)\dom \j_\pi(h)
\subseteq\dom \j_\pi(h)$.
\item[(ii)]
$\lim\limits_{\lambda\to\infty}i\lambda\pi(\rlf)\Psi=\Psi$ for all
$\Psi\in\al H._\pi$.
\item[(iii)]
$\lim\limits_{\mu \to 0}i\pi(R(1, \mu f))\Psi=\Psi$ for all
$\Psi\in\al H._\pi$.
\item[(iv)]
The space $\al D.:={\pi\big(R(1,f)R(1,h)\big)\al H._\pi}$
is a joint dense domain 
for $\j_\pi(f)$ and $\j_\pi(h)$ 
and we have:
$[\j_\pi(f),\,\j_\pi(h)]=i\sigma(f,h) \un $ on $\al D..$
\item[(v)] $\ker\pi\big(R(1,\nu f+h)\big)=\{0\}$
for $\nu \in\R$. 
Then $\j_\pi(\nu f+h)$ is defined,  
$\al D.$ 
is a core for  $\j_\pi(\nu f+h)$ and
$\j_\pi(\nu f+h)=\nu \j_\pi(f)+\j_\pi(h)$
 on $\al D..$ Moreover
${\pi\big(R(1,\nu f+h)\big)}\in
{\big\{\pi\big(R(1,f)\big),\,\pi\big(R(1,h)\big)\big\}''.}$
\item[(vi)]
$\j_\pi(f)\pi(\rlf)=\pi(\rlf)\j_\pi(f)=i\lambda\pi(\rlf)-\un$ on
$\dom \j_\pi(f)$.
\item[(vii)]
$\big[\j_\pi(f),\pi(R(\lambda,h))\big]=i\sigma(f,h)\pi(R(\lambda,h)^2)$
on $\dom \j_\pi(f)$.
\item[(viii)]
Denote $W(f):=\exp(i\j_\pi(f))$, then
\begin{eqnarray*}
W(f)W(h) &=& e^{- i\sigma(f,h)/2}\,W(f+h)  \\[1mm]
W(f)\pi\big(R(\lambda,h)\big)W(f)^{-1} &=&
\pi\big(R(\lambda+i\sigma(h,f),\,h)\big)=\pi\big(\beta_{f_\sigma}(R(\lambda,h))\big)
\end{eqnarray*}
where $f_\sigma:X\to\R$ is given by $f_\sigma(h):=\sigma(h,f)$, $h\in X$.
Moreover  $W(sf)\al D.\subseteq\al D.\supseteq W(th)\al D.$ for $s,\,t\in\R$, hence
 $\al D.$ is a common core
for $\j_\pi(f)$ and $\j_\pi(h)$.
\end{itemize}
\end{teo}
Thus we define:
\begin{defi}
A  representation $\pi\in\rep\big(\al R.(X,\,\sigma),\al H._\pi\big)$
is {regular on} a set $S\subset X$ if
\[
\ker\pi\big(R(1,f)\big)=\{0\}\qquad\mbox{for all} \ f\in
S\;.
\]
A state $\omega$ of $\rsl$ is {regular on} a set $S\subset X$ if its GNS--representation
$\pi_\omega$ is regular on $S\subset X$.
A {regular representation (resp.\ state)} is a
representation (resp.\ state) which is regular on $X$. 
Given a Hilbert space $\al H.,$ we denote 
the set of (nondegenerate) regular representations $\pi:\rsl\to\al B.(\al H.)$
by ${{\rm Reg}\big(\rsl,\al H.\big)}.$
The set of regular states of $\rsl$ is denoted by
$\ot S._r\big(\rsl\big).$
\end{defi}
Obviously many regular representations are known, \eg the
Fock representation.
Recall that the class of all regular representations of $\rsl$ is not a set,
hence the necessity to fix $\al H..$
Thus for $\pi\in{{\rm Reg}\big(\rsl,\al H.\big)},$ all the field operators
$\j_\pi(f),$ $f\in
X$ are defined, and we have the resolvents 
$\pi(\rlf)=(i\lambda\un-\j_\pi(f))^{-1}$.

{}From Theorem~\ref{RegThm}, we can now establish a bijection between the
regular representations of $\rsl$ and the regular representations
of the Weyl algebra $\CCRX:$
\begin{cor}
\label{RegBij} Let $\rsl$ be as above.
Given a regular representation $\pi\in{\rm Reg}\big(\rsl,\al H.\big),$ define a
regular representation $\wt\pi\in{\rm Reg}\big(\CCRX,\al H.\big)$
by $\wt\pi(\delta_f):=\exp(i\j_\pi(f))=W(f)$
(using Theorem~\ref{RegThm}(viii)). 
This correspondence establishes a bijection 
between ${\rm Reg}\big(\rsl,\al H.\big)$ and ${\rm Reg}\big(\CCRX,\al H.\big)$
which respects irreducibility and direct sums. Its inverse is given by
the Laplace transform,
\begin{equation}
\label{Laplace1}
\pi(R(\lambda,f)):= - i\int_0^{\, \sigma \infty}  e^{-\lambda t}\pi
( \delta\s -tf.)\,dt\;, \qquad \sigma := \mbox{sign} \, \lambda \,.
\end{equation}
By an application of this to the GNS--representations of regular states, 
we also obtain an affine bijection between $\ot S._r\big(\rsl\big)$
and the regular states $\ot S._r\big(\CCRX\big)$ of $\CCRX,$
and it restricts to a bijection between the pure regular states of
$\rsl$ and the pure regular states of $\CCRX.$
\end{cor}
Note that whilst we have a bijection between the regular states of
$\rsl$ and those of $\CCRX$, there is no such map between the
nonregular states of the two algebras. In fact, fix a nonzero $f\in X$ and
consider the two commutative subalgebras 
$C^*\{\rlf\, , \un \, \mid\,\lambda\in\R\backslash 0\} \subset\rsl$ and
$C^*\{\delta_{tf}\,\mid\,t\in \R\}\subset\CCRX$, then these are isomorphic
respectively to the continuous functions on
the one point compactification of $\R$, and the continuous functions on
the Bohr compactification of $\R$.
Note that the point measures on
the compactifications without $\R$
 produce nonregular states (after extending to the full C*--algebras by
 Hahn--Banach) and there are
many more of these for the Bohr compactification than for the one 
point compactification of $\R$,
(cf. Theorem~5 in~\cite[p 949]{DS2}).
So the Weyl algebra has many more nonregular states than
the resolvent algebra.

Some further properties of regular representations and states are:
\begin{pro} 
\label{RegAlg}
Let $\al R.(X,\,\sigma)$ be as above. 
\begin{itemize}
\item[(i)] If a representation $\pi$ of $\rsl$ is faithful and
factorial, it must be regular.
\item[(ii)] If a representation $\pi:\rsl\to\al B.(\al H.)$ is regular
then ${\|\pi(\rlf)\|}=\|\rlf\|=|\lambda|^{-1}$ for all $\lambda\in\R\backslash 0,$
$f\in X$.
 \item[(iii)] A state $\omega$ of $\rsl$ is regular iff
 $\omega(A)=\lim\limits_{\lambda\to\infty}i\lambda\,\omega\big(\rlf A\big)$
 for all $A\in\rsl$ and $f\in X$.
\end{itemize}
\end{pro} 
Thus regular states restrict to regular states on subalgebras.
Clearly a cyclic component of a regular representation is again regular,
so it makes sense to define
\begin{defi}
The {universal regular representation} of $\rsl$ is
\[
\pi_{r}:=\bigoplus\big\{\pi_\omega\,\mid\,\omega\in\ot S._r\big(\rsl\big)
\big\}\quad\hbox{and the {regular seminorm} is}\quad\|A\|_r:=\|\pi_r(A)\|
\]
and we denote $\al R._r(X,\sigma):=\pi_r\big(\rsl\big).$
\end{defi}
Now $\pi_r$ is a subrepresentation of $\pi_u$ and all
regular representations of $\rsl$ will factor through $\pi_r.$
We want to prove that $\pi_r$ is faithful, and hence that 
$\al R._r(X,\sigma)\cong\rsl$. {}For this, we need to develop some structure theory
for  general representations. First, some notation: for a
subspace $S\subset X$ its symplectic complement
will be denoted by $S^\perp:={\set f\in X,\sigma(f,S)=0.}$.
By $X=S_1\oplus S_2\oplus\cdots\oplus S_n$ we will mean that all
$S_i$ are nondegenerate and $S_i\subset S_j^\perp$ if $i\not=j,$
and each $f\in X$ has a unique decomposition
$f=f_1+f_2+\cdots+f_n$ such that $f_i\in S_i$ for all $i$.
\begin{pro}
\label{Xdecomp}
Let $\pi:\rsl\to\al B.(\al H.)$ be a nondegenerate representation. Then
\begin{itemize}
\item[(i)]
the set $X_R:=\set f\in X,\ker\pi\big(R(1,f)\big)=\{0\}.$ is a linear space.
Hence if $f\in X_S:=X\backslash X_R,$ then $f+g\in X_S$ for all $g\in X_R$.
\item[(ii)]
The set $X_T:=\set f\in X,\ker\pi\big(R(1,f)\big)=\{0\}\;\hbox{and}\;\pi\big(R(1,f)\big)^{-1}\in
{\al B.(\al H.)}.\subset X_R$ is a linear space. Moreover if $f\in X_T$ then
$\pi\big(R(1,g)\big)=0$ for all $g\in X$ with $\sigma(f,g)\not=0$.
Thus $\sigma(X_T,X_R)=0$.
\item[(iii)]
If $\pi$ is factorial, then $\pi\big(R(1,f)\big)=0$ for all $f\in X_S,$ and
$\pi\big(R(1,f)\big)\in\C\un\backslash 0$ for all $f\in X_T$. 
Moreover $X_T=X_R\cap X_R^\perp$.
\item[(iv)] Let $X$ be finite dimensional 
and let ${\big\{q_1,\ldots,\,q_n\big\}}$ be a basis for  $X_T$.
If $\pi$ is factorial, we can augment this basis of  $X_T$ by
${\big\{p_1,\ldots,\,p_n\big\}}\subset X_S$ into a symplectic
basis of $Q:={\rm Span}{\big\{q_1,p_1;\ldots;\,q_n,p_n\big\}},$
\ie $\sigma(p_i,q_j)=\delta_{ij}$, $0=\sigma(q_i,q_j)=\sigma(p_i,p_j)$.
Then we have the decomposition
\begin{equation}
\label{Qdecomp}
X=Q\oplus(Q^\perp\cap X_R)\oplus(Q^\perp\cap X_R^\perp)
\end{equation}
into nondegenerate spaces  such that
$Q^\perp\cap X_R\subset\{0\}\cup(X_R\backslash X_T)$ and
$Q^\perp\cap X_R^\perp\subset\{0\}\cup X_S$.
\end{itemize}
\end{pro}
Clearly $X_R$ is the part of $X$ on which $\pi$ is regular,
$X_T$ is the part on which it is ``trivially regular'', 
$X_S$ is the part on which it is singular, and these have a particularly nice form when
$\pi$ is factorial.
So by this Proposition we obtain an interesting regularity structure theory for representations
of $\rsl$ which we will exploit below to prove our uniqueness theorem.
Note that (ii) implies that if $\pi$ is regular for a pair $f,\,g\in X$
with $\sigma(f,g)\not=0$, then the field operators
$\j_\pi(f)$ and $\j_\pi(g)$ are both unbounded.
\begin{pro}
\label{NonregApprox}
Let $X$ be finite dimensional, and let $\omega$ be a nonregular 
pure state of $\rsl$.  Then
\begin{itemize}
\item[(i)]
there is a sequence $\{\omega_n\}\subset
\ot S._r\big(\rsl\big)$ of regular states such that
$\omega=\lim\limits_{n\to\infty}\omega_n$ in the w*--topology.
\item[(ii)] $\ker\pi_r\subseteq\ker\pi_\omega$ where
$\pi_r$ is the universal regular representation of $\rsl$.
\item[(iii)]
$\rsl\cong \al R._r(X,\sigma)$.
\end{itemize}
\end{pro}
Using this, we can now prove our desired uniqueness theorem, \ie that
$\rsl\cong \al R._r(X,\sigma)$ for any $X$.
\begin{teo}
\label{UniqueR}
Let $X$ be a nondegenerate symplectic space of arbitrary dimension. Then
\begin{itemize}
\item[(i)]
The norms of $\al R.(X,\,\sigma)$ and $\al R.(S,\,\sigma)$ coincide on 
 ${\hbox{*-alg}\{\rlf\,\mid\,f\in S,\,\lambda\in\R\backslash 0\}}$ for each
finite dimensional nondegenerate subspace $S\subset X$.
Thus we obtain a containment ${\al R.(S,\,\sigma)}\subset\al R.(X,\,\sigma).$
\item[(ii)]$\al R.(X,\,\sigma)$ is the inductive limit of the net of all
${\al R.(S,\,\sigma)}$ where $S\subset X$ ranges over all finite dimensional
nondegenerate subspaces of $X.$
\item[(iii)] We have that
$\rsl\cong \al R._r(X,\sigma)$.
\end{itemize}
\end{teo}
It follows therefore from Fell's theorem (cf. Theorem~1.2 in \cite{Fell}) 
that {\it any} state of $\rsl$ is in the w*-closure of the convex hull of
the vector states of $\pi_r$, hence of
the regular states. We can now prove the following result, 
which is relevant for physics.
\begin{teo}
\label{RegFaith}
Let $(X,\,\sigma)$ be any nondegenerate symplectic space, and 
$\rsl$ as above. Then every regular representation of $\rsl$ is faithful.
\end{teo}
The importance of this result lies in the fact that the regular representations
are taken to be the physically relevant ones, and the images of $\rsl$
in all regular representations are isomorphic. Thus, since we can obtain
the quantum fields from $\rsl$ in these representations, we are justified 
in taking $\rsl$ to be the observable algebra for bosonic fields.
Usually one argues that for a C*--algebra $\al A.$ to be an  
observable algebra of a physical system,
it must be simple (cf.~\cite[p 852]{HaKa}). 
The argument is that by Fell equivalence of
the physical representations, the image of $\al A.$ in all
physical representations must be isomorphic. However, if one restricts the class
of physical representations (as we do here to the regular representations
of $\rsl),$
then the latter isomorphism does not imply that $\al A.$ must be simple.

This theorem also has structural consequences, \eg it implies that
$\rsl$ has faithful irreducible representations, hence that its centre must
be trivial. Moreover, it 
also proves isomorphism with the previous version of the 
resolvent algebra which we developed in~\cite{BuGr}, because it was the image of $\rsl$ 
in the ``universal strongly regular representation'' which we now define.

{}For many applications one needs regular representations where there is
a dense invariant joint domain for all the fields  $\j_\pi(f),$ and this 
leads us to a subclass of the regular representations as follows.
We will say that a state $\omega$ on the Weyl algebra $\CCRX$
is {strongly regular} if the functions
\[
\R^n\ni(\lambda_1,\ldots,\lambda_n)\mapsto\omega\big(\delta_{\lambda_1f_1}
\cdots\delta_{\lambda_nf_n}\big)
\]
are smooth for all $f_1,\ldots,\,f_n\in X$ and all 
$n\in\N$. Of special importance is that the GNS-representation 
of a strongly regular state has a common dense invariant domain
for all the generators $\j_{\pi_\omega}(f)$ of the one parameter groups
$\lambda\to\pi_\omega(\delta_{\lambda f})$
(this domain is obtained by applying the polynomial algebra of the Weyl operators
$\set\pi_\omega(\delta_f),f\in X.$ to the cyclic GNS-vector).
By the bijection of Corollary~\ref{RegBij},
we then obtain the set of strongly regular states on $\rsl,$ 
and we denote this by ${\ot S._{sr}\big(\rsl\big)}.$
\begin{defi}
The {universal strongly regular representation} of $\rsl$ is
\[
\pi_{sr}:=\bigoplus\big\{\pi_\omega\,\mid\,\omega\in\ot S._{sr}\big(\rsl\big)
\big\} \, .
\]
\end{defi}
$\pi_{sr}\big(\rsl\big)$ is the version of the resolvent algebra which we used in \cite{BuGr},
which is obviously isomorphic to $\rsl$ by Theorem~\ref{RegFaith}.

An important set of strongly regular states on $\CCRX$ are the
quasifree states. They are given by 
\[
\omega(\delta_f) = \exp\big(-\hlf\langle f | f \rangle_\omega\big) , 
\quad f \in X,
\]
where $\langle \, \cdot \, | \, \cdot \, \rangle_\omega $ 
is a (possibly semi--definite) scalar product on the complex 
linear space $X + i X$ satisfying 
$$ \langle f | g \rangle_\omega- 
\langle g | f \rangle_\omega = i \sigma(f,g), 
\quad f,g \in X. $$
By a routine computation one can represent the expectation
values of products of Weyl operators in a quasifree state 
in the form
$$ \omega(\delta_{f_1} \cdots \delta_{f_n}) = 
\exp\Big(- \sum_{k<l}  \langle f_k | f_l \rangle_\omega 
-\hlf \sum_l   \langle f_l | f_l \rangle_\omega \Big). $$
Making use of the Laplace transform~(\ref{Laplace1}) for the GNS--representation
of the resolvents, we 
 have for $\lambda_1, \dots \lambda_n > 0$
\begin{eqnarray}
& & \omega(R(\lambda_1,f_1) \cdots R(\lambda_n,f_n))  \nonumber\\[1mm]
& & = (-i)^n 
\int_0^\infty \! dt_1 \dots  \int_0^\infty \! dt_n \,
e^{-{\sum}_k t_k \lambda_k} \omega(\delta\s t_1 f_1. \cdots
\delta\s t_n f_n.) \nonumber\\[1mm]
\label{QFresolvent}
& & =  (-i)^n  
\int_0^\infty \! dt_1 \dots  \int_0^\infty \! dt_n \,
\exp\Big(-\sum_k t_k \lambda_k
- \sum_{k<l} \, t_k t_l \langle f_k | f_l \rangle_\omega 
-\hlf \sum_l \,  t_l^2 \langle f_l | f_l \rangle_\omega\Big). 
\end{eqnarray}

\begin{rem} One can replace anywhere in this
equation $f_k$ by $-f_k$, thus it does not impose any restriction of
generality to assume that $\lambda_1, \dots \lambda_n > 0$. The
 relation (\ref{QFresolvent}) should be regarded as the definition 
of quasifree states on the resolvent algebra.
\end{rem}

Concerning unitary implementability of automorphisms, we have the following
easy fact:
\begin{pro}
\label{Uaut}
Let $\alpha\in\aut\rsl$ correspond to a symplectic transformation
 $T\in{\rm Sp}(X,\sigma)$ 
by $\alpha\big(\rlf\big):=R(\lambda,Tf)$ for $f\in X$.
Then $\alpha$ is implemented by a unitary in both $\pi_r$ and $\pi_{sr}$.
\end{pro}

\section{Further structure.}
\label{FurStruc}

Here we want to explore the algebraic structure of $\al R.(X,\,\sigma)$.
\begin{teo}
\label{TensorAlg}
Let $(X,\,\sigma)$ be a given nondegenerate symplectic space, and let 
$X=S\oplus S^\perp$ for $S\subset X$ a nondegenerate subspace. Then
$$\al R.(X,\,\sigma)
 \supset{\rm C}^*\big(\al R.(S,\,\sigma)\cup\al R.(S^\perp,\,\sigma)\big)
\cong\al R.(S,\,\sigma)\otimes\al R.(S^\perp,\,\sigma)$$
where the tensor product uses the minimal (spatial) tensor norm.
The containment is proper in general.
\end{teo}
Thus we cannot generate $\rsl$ from a basis alone, \ie if $\{q_1,p_1;\,q_2,p_2;\cdots\}$
is a symplectic basis of $X$,  then 
${\rm C}^*\big\{R(\lambda_i,q_i),\;R(\mu_i,p_i)\,\mid\,\lambda_i,\,\mu_i\in\R\backslash 0,\;
i=1,\,2,\ldots\big\}$ is in general a proper subalgebra of $\rsl,$ though in any
regular representation $\pi$ it is strong operator dense in
$\pi\big(\rsl\big)$ by Theorem~\ref{RegThm}(v).

Note that since ${C^*(\{\rlf\,\mid\,\lambda\in\R\backslash 0\})}\cong C_0(\R)$
(easily seen in any regular representation), and we have that 
$C_0(\R^{n+m})=C_0(\R^n)\otimes C_0(\R^m),$ it follows from Theorem~\ref{TensorAlg}
that any $C_0\hbox{--function}$ of a finite commuting set of variables 
is in $\rsl$. More concretely, we have the following result which will
be used later.
\begin{pro}
\label{Czero}
Let $\{q_1,\ldots,\,q_k\}\subset X$ satisfy $\sigma(q_i,q_j)=0$
for all $i,\,j$. Then for each
$F \in C_0(\R^k)$ there is a (unique) $R_F \in\rsl$ such that
in any regular representation $\pi$
we have $\pi(R_F)={F \big(\j\s\pi.(q_1),\ldots,\j\s\pi.(q_k)\big)}$.
\end{pro}
Thus the resolvent algebra contains in abstract form 
all $C_0\hbox{--functions}$ of commuting fields. 
Note that such a result neither holds for the 
Weyl algebra nor for the corresponding 
twisted group algebra (in the case of finite
dimensional $X$).

\begin{teo}
\label{Nonsep}
Let $(X,\,\sigma)$ be a given nondegenerate symplectic space 
and let $f,\,h\in X\backslash 0$ such that  
$f\not\in\R h$. Then
\begin{itemize}
\item[(i)] $R(1,f)\not\in\big[\rsl R(1,h)\big]$, \ie the ideals separate the rays 
of $X$,
\item[(ii)] $\big\|R(1,f)-R(1,h)\big\|\geq 1$, and if $\sigma(f,h)=0$ we have equality.
\item[(iii)] $\rsl$ is nonseparable.
\end{itemize}
\end{teo}
Next, let us assume that our symplectic space $X$
is finite dimensional, hence by the von Neumann uniqueness theorem there
is (up to unitary equivalence) a unique irreducible regular representation $\pi_0$
of $\rsl$.
\begin{teo}
\label{CompactId}
Let $(X,\,\sigma)$ be a given finite dimensional nondegenerate symplectic space
equipped with the symplectic basis
 $\big\{q_1,p_1;\ldots;\,q_n,p_n\big\}$,
and let $\pi_0:\rsl\to\al B.(\al H._0)$ be an irreducible 
regular representation of $\rsl$. Then 
\begin{itemize}
\item[(i)] $\pi_0\Big(\big(R(\lambda_1,p_1)R(\mu_1,q_1)\big)\cdots
\big(R(\lambda_n,p_n)R(\mu_n,q_n)\big)\Big)$ is a Hilbert--Schmidt operator
for all $\lambda_i$, $\mu_i\in\R\backslash 0$.
\item[(ii)] There is a unique closed two-sided ideal $\al K.$ of $\rsl$ which is isomorphic to
the algebra of compact operators $\al K.({\al H.}_0)$, and such that 
$\pi_0(\al K.)=\al K.(\al H._0)\subset\al B.(\al H._0)$.
\item[(iii)] A representation of $\rsl$ is regular iff its restriction to
$\al K.$ is nondegenerate. Thus the regular representations (resp. regular states) 
are exactly the unique extensions of the representations (resp. states)
of $\al K.$ to $\rsl$.
\item[(iv)] If $n=1$ then the factor algebra $\rsl\big/\al K.$ is commutative,
but not if $n>1$.
\item[(v)] $\al K.$ is an essential ideal, \ie if
 $A\al K.=0$ or $\al K.A=0$, then $A=0$.
\item[(vi)] $\al K.$ is a minimal (nonzero) closed two--sided ideal, and is contained
in every closed nonzero two--sided ideal of $\rsl$.
Thus all closed nonzero two--sided ideals of $\rsl$ are essential.
\end{itemize}
\end{teo}
\begin{rems}
\label{CompactRem}
(a) \ 
The statement in (iii) that the regular representation theory of $\rsl$
is the representation theory of the compacts, is of course just a paraphrasing
of the  von Neumann uniqueness theorem. 
In fact, it is well--known that the C*--closure of the 
twisted convolution algebra 
$L^1(X)$ w.r.t. $\exp(i\sigma)$ (which is the
twisted group algebra of $X$ thought of as an abelian group) is isomorphic to
the compacts, and that the Weyl algebra acts on it by multipliers (cf. 
first part of \cite{Gr3} for a discussion of this). The regular representations
on the Weyl algebra are likewise obtained 
for finite dimensional $X$ from the unique extensions of the
representations of the compacts. The main attractive feature of the resolvent
algebra is that the algebra of compacts is actually an ideal inside $\rsl$, in contrast
with the Weyl algebra where it is outside. This is very useful, and we will
utilize this fact below. \\[1mm]
(b) \  {}From the structure above, we note that the nonregular representations
are those which have a direct summand which is the restriction of a representation
of the Calkin algebra to $\rsl$.  \\[1mm]
(c) \ {}From (vi) and (iii) we obtain a quick proof of Theorem~\ref{RegFaith}
as follows. Let $\pi$ be a regular representation of $\rsl$ where $X$ is arbitrary.
Since $\rsl$ is an inductive limit of $\al R.(S,\sigma)$ where $S\subset X$ ranges over 
the finite dimensional nondegenerate subspaces, it suffices to show that
$\pi$ is faithful on each $\al R.(S,\sigma)$. But if $\pi$ is not 
faithful on a $\al R.(S,\sigma)$ with ${\rm dim}(S)<\infty,$ then by (vi)
$\pi$ must vanish on $\al K.$ hence cannot be regular by (iii) which
is a contradiction.  \\[1mm]
(d) \ If $n>1$ then $\rsl$ contains infinitely many copies of the compacts,
because for each nondegenerate proper subspace $S\subset X$, the copy of 
$\al R.(S,\sigma)$ in $\rsl$ will contain its own compact ideal
$\al K._S\subset\al R.(S,\sigma)$. In general these will not be ideals of $\rsl$,
nor will they map onto the compacts in the irreducible regular representation.
This can be seen from the fact that $X=S\oplus S^\perp,$ hence $\al R.(S,\sigma)$
is embedded as $\al R.(S,\sigma)\otimes\un$ in the subalgebra
$\al R.(S,\sigma)\otimes\al R.(S^\perp,\sigma)$ of $\rsl$. It is now clear that
$\al K._S\otimes\un$ is not an ideal, and in the Schr\"odinger representation
w.r.t. the union of a symplectic basis of $S$ and of $S^\perp$ that
$\al K._S\otimes\un$ maps to a tensor product of compacts with the identity,
which is not compact. These embedded copies of $\al K.$ are nevertheless useful
as indicators of partial regularity, \ie a representation
 $\pi$ is regular on a nondegenerate
subspace $S\subset X$ iff it is nondegenerate on $\al K._S$. This can be particularly
useful for the infinite dimensional case, where we do not have a $\al K.\subset\rsl,$
but we still have that
 $\pi$ is regular iff its restriction to $\al K._S$ is nondegenerate for
each finite dimensional nondegenerate subspace $S\subset X$.
\end{rems} 

If $X$ is infinite dimensional, the question naturally arises as to whether
there is a C*--algebra $\al L.$ which can play the role of $\al K..$
That is, we want at least that there is a faithful 
embedding of $\rsl$ in 
the multiplier algebra $M(\al L.)$ and such that the unique extensions of
representations from $\al L.$ to $\rsl\subset M(\al L.)$ produces an identification
of the representation theory of $\al L.$ with the regular representations
of $\rsl$. In the case that $S$ is countably dimensional such an algebra
has recently been constructed for the Weyl algebra, cf.~\cite{GrNe},
hence is a strong candidate.

\section{Dynamics and Hamiltonians}
\label{DynHam}

{}For the dynamics (\ie time evolution) of a quantum mechanical
system, one usually assumes a 
one--parameter automorphism group of the observable algebra, 
and considers distinguished representations 
in which it is implemented by strong operator continuous unitary groups.
Much analysis is then done of the generators (Hamiltonians).
On both the Weyl algebra $\CCRX$ and the resolvent algebra $\rsl$
one can define dynamical groups in terms of symplectic
transformations. However, 
since such dynamics correspond to quadratic Hamiltonians, 
it does not cover most physically interesting situations. 
In fact, for the Weyl algebra it is a well--known 
frustration that many natural time--evolutions cannot be defined on it. 
{}For example in \cite{FaVB} it is proven for the case $X=\R^2$
that time evolutions
obtained from Hamiltonians of the form
$H= P^2+V(Q)$ for potentials $V\in L^\infty(\R)\cap L^1(\R)$ 
do not preserve $\CCRX$ in the Schr\"odinger representation 
unless $V=0$.  Moreover, there is no adequate method known
by which one can somehow encode the dynamics (or even more desirably,
the Hamiltonian) in the Weyl algebra. 

We want to show here that for
the resolvent algebra $\rsl$ the situation is completely different.
We will not give a comprehensive analysis of
the issue of dynamics but will rather 
illustrate the technical benefits
of the resolvent algebra for the analysis of dynamics questions. To simplify
matters, we will assume in the rest of this section that $X=\R^2$ 
with its standard symplectic form $\sigma$ and work in the Schr\"odinger
representation of $\rsl$. 
With some effort, our results
can be extended to the case of arbitrary finite dimensional 
and non--degenerate symplectic spaces $X$.
In our example in the next section we will touch on dynamics (with interaction)
when $X$ is infinite dimensional.

\subsection{Dynamics on $\rsl$} \label{dynamics}

Henceforth, for this section, let $X=\R^2$ with
its standard symplectic form.
Since the Schr\"odinger representation $\pi_0$ is faithful on $\rsl$, 
the dynamics on $\rsl$ can be defined in $\pi_0$
as follows. Let $H$ be a selfadjoint operator (Hamiltonian) generating
the unitary group $U(t) = e^{itH}$, $t \in \R$. In favourable 
cases where
\[ U(t) \pi_0(\rsl) U(t)^{-1} \subset  \pi_0(\rsl) \, , \quad t \in \R,
\]
one can define a corresponding automorphic action $\alpha_t$ 
on $\rsl$, $t \in \R$ putting
\[ \alpha_t(R) := \pi_0^{-1} 
\big( U(t) \pi_0(R) U(t)^{-1} \big) \, , \quad R \in
\rsl \, .
\]
We will say in these cases that the Hamiltonian $H$ induces a dynamics 
on the resolvent algebra $\rsl$ (where the context of $\pi_0$ is assumed).
{}From the case when $H$ is quadratic, we see via Theorem~\ref{Nonsep}(ii)
that in general the actions $t\mapsto
\alpha_t$ need not be pointwise norm--continuous.

The Hamiltonians which induce a dynamics on $\rsl$
exist in abundance. Denote by $Q, P$ the canonical 
position and momentum operators in the Schr\"odinger representation,
then we have the following result in contrast with the no--go
theorem for the Weyl algebra above.

\begin{pro} \label{Dynamics}
Let $V \in C_0(\R)$ be any real function. 
Then the corresponding selfadjoint Hamiltonian $H = P^2 + V(Q)$ 
induces a dynamics on $\rsl$.
\end{pro}
Thus the resolvent algebra is an appropriate framework for the formulation
of dynamical laws.

\subsection{Hamiltonians affiliated with $\rsl$} \label{affil}

It is another interesting feature of the resolvent algebra $\rsl$ 
that it contains many observables of physical interest. We 
illustrate this fact within the preceding concrete setting. 

Let $H$ be a selfadjoint operator. When its resolvent is 
contained in $\pi_0(\rsl)$ we may proceed to its pre--image,
\[ R_\lambda := \pi_0^{-1} \big((i\lambda \un - H)^{-1}\big) \, ,
\quad  \lambda \in \R\backslash \{0 \} \, ,
\] 
defining a pseudo--resolvent in $\rsl$.
By a slight abuse of terminology we say then that $H$
is affiliated with  $\rsl$ and regard it as an observable 
which is predetermined by the resolvent algebra.
In the regular representation $\pi_0$ the pseudo--resolvent
$R_\lambda$ produces the selfadjoint Hamiltonian $H$ which 
is the energy operator of the underlying states. 
But $R_\lambda$ can also be used for the determination of 
the energy content of non--regular representations
including extreme cases where the energy of all
states is infinite and $R_\lambda$ is represented by $0$.
The next result exhibits 
a multitude of Hamiltonians which are affiliated 
with the resolvent algebra.

\begin{pro} \label{Affil}
Let $V \in C_0(\R)$ be any real function. 
Then the corresponding selfadjoint Hamiltonian $H = P^2 + V(Q)$ 
is affiliated with $\rsl$.
\end{pro}
Since $X$ is finite dimensional $\rsl$ 
also contains the compact operators and hence in particular
all one--dimensional projections. As these projections play a fundamental 
role in quantum systems, we conclude that
the resolvent algebra contains all necessary ingredients for 
the treatment of quantum mechanical systems. It therefore
qualifies as a genuine observable algebra.

\subsection{Hamiltonians not affiliated with $\rsl$} \label{nataffil}

Since $\pi_0(\rsl)$ is a proper subalgebra of the algebra of all
bounded operators on the underlying representation space it is clear
from the outset that there are many selfadjoint operators which are
not affiliated with $\rsl$. It is thus of  
interest to explore which operators do not qualify as observables
in the framework of the resolvent algebra. We do not have a complete
answer to this question, not even for the case of quantum mechanical
Hamiltonians. But the next result shows that
many Hamiltonians which seem physically unacceptable 
can be excluded this way. 

\begin{pro} \label{Notaffil}
The selfadjoint Hamiltonian $H = P^2 - Q^2$ 
is not affiliated with $\rsl$.
\end{pro}
Since this Hamiltonian is not
semibounded it does not describe a stable system and hence is of limited
physical interest. Note that, in contrast, the Hamiltonian  $H = P^2 + Q^2$
of the harmonic oscillator is an observable which 
is affiliated with $\rsl$ since its
resolvent is a compact operator. 

\section{Infinite dimensional dynamical systems}
\label{InteractSect}

The resolvent algebra exhibits its full power when dealing
with systems for which ${\rm dim}(X)$ is infinite.
{}For, in contrast to the finite dimensional case, there
exists then an abundance of disjoint regular representations. It is
a notorious problem in quantum physics to select
the representations which describe the states of interest for a 
given infinite dimensional system. This selection is usually done by
the specification of a 
dynamics, and searching for representations describing specific
situations such as ground states or thermal equilibrium states.
It is an asset
of the resolvent algebra that it admits the 
definition of interesting dynamics also in the case of infinite 
dimensional systems. Moreover, it simplifies the construction of states
and representations of
interest by the use of C*--algebra techniques.

Next, we develop a concrete model to illustrate these facts.
Our model
describes particles which are confined around the points of a
lattice by harmonic pinning forces and which interact with their nearest 
neighbours. Thus it resembles the situation in
quantum spin systems  \cite{BR2}, but in contrast to the latter 
well--understood class of theories, the present model has an 
infinite number of degrees of freedom at each lattice site. 

Turning to mathematics, let  ${(X,\,\sigma)}$ be a countably 
dimensional symplectic space  with a fixed
symplectic basis ${\big\{p_l ,q_l\;\mid\; l \in\Z\big\}}$,
where the index $l \in\Z$ labels the lattice points.
Denote $X_l :={{\rm Span}\{p_l ,q_l \}}$ and
$X_\Lambda := {{\rm Span}\{p_l ,q_l \; \mid l \in \Lambda} \}$,
where $\Lambda \subset \Z$ is any finite subset of lattice points.
Then we take the observable algebra of the model to be $\rsl$, and it
is the C*--inductive
limit of the algebras ${\cal R}(X_\Lambda, \sigma)$, $\Lambda \subset
\Z$.

To define the dynamics on $\rsl$ we 
fix a regular (hence faithful) representation $\pi_0$ of $\rsl$,
and here we will take $\pi_0$ to be the Fock representation
w.r.t. the given basis.
Recall that this Fock representation is the (irreducible) representation
which is determined by the fact that there exits a unit vector
$\Omega_0 \in {\cal H}_0$ in the domain of all polynomials 
of the fields satisfying
\[ \big( \phi_{\pi_0}(p_l) + i \phi_{\pi_0}(q_l) \big) \, \Omega_0 = 0
\, , \quad l \in \Z \, .
\] 
As is well--known, this vector defines a product state on the algebra
$\ccr X,\sigma.$ in the sense that if $S$ and $T$ are nondegenerate subspaces of
$X$ with $S\cap T=\{0\}$ then
\[ \big( \Omega_0, \pi_0(W_1) \pi_0(W_2) \Omega_0 \big) =
\big( \Omega_0, \pi_0(W_1) \Omega_0 \big) \, 
\big( \Omega_0, \pi_0(W_2) \Omega_0 \big) \, , \quad 
W_1 \in \ccr S,\sigma. ,\; 
W_2 \in \ccr T,\sigma.\,.  
\] 
By strong operator limits and continuity we obtain the analogous statement
for $\rsl$ \ie $\Omega_0$ also defines a product state for $\rsl$
in a similar sense.

As a consequence, for any nondegenerate subspace $S\subset X$ the restriction 
$\pi_0 \upharpoonright {\cal R}(S, \sigma) $ acts
irreducibly on the corresponding subspace
$ {\cal H}_0 (S) = 
\overline{\pi_0({\cal R}(S, \sigma)) \, \Omega_0}$ 
and it is equivalent to the Fock
representation of ${\cal R}(S, \sigma)$. 
Moreover $\pi_0({\cal R}(S, \sigma))^{\prime \prime} $
is a factor.
In particular, these statements hold for $S=X_\Lambda$,
$\Lambda\subset\Z$.
Even though our results do not depend on the choice of      
representation $\pi_0$, its specific features will greatly simplify 
the necessary computations. 

\subsection{Interacting dynamical systems} \label{interaction}

Turning to the problem of defining 
the dynamics, denote:
\[ Q_l := \phi_{\pi_0} (p_l), 
\ P_l := \phi_{\pi_0} (q_l) \, , \quad l \in \Z \, .
\] 
Let $V \in C_0(\R)$ be any real function which models the interaction
potential.
We consider for each $\Lambda \subset \Z$ 
the selfadjoint ``local Hamiltonian'' 
\[ 
H_{\Lambda} :=
\sum_{l \in \Lambda} (P_l^2 + Q_l^2) +
\sum_{l, \, l+1 \in \Lambda} V(Q_l - Q_{l+1}) \, ,
\]
describing the interaction amongst the particles   
in $\Lambda$. 
The dynamics on the full system is then obtained as follows.

\begin{pro} \label{InteractDynamics}
Let $V \in C_0(\R)$ be a real function 
and let $\{ H_\Lambda \}_{\Lambda \subset \Z}$
be the corresponding set of Hamiltonians as above, and define
the unitary groups $U_\Lambda(t) :=e^{\, it H_\Lambda}$, $t \in \R$. 
Then
\begin{itemize}
\item[(i)] $U_{\Lambda} (t) \, \pi_0({\cal R}(X_\Lambda, \sigma)) \, U_{\Lambda} (t)^{-1}
\subset \pi_0 ({\cal R}(X_\Lambda, \sigma))$ for all $t \in \R$ and $\Lambda \subset \Z$.
\item[(ii)] For any $R \in  {\cal R}(X_{\Lambda_0}, \sigma)$,
$\Lambda_0 \subset \Z$ the net 
$\{U_{\Lambda} (t) \pi_0(R) U_{\Lambda} (t)^{-1} \}_{\Lambda \subset  \Z}$ 
converges in norm to an element of $\pi_0 (\rsl)$ 
as $\Lambda \nearrow \Z$.
\item[(iii)] There is a unique automorphic action (dynamics)
$\alpha:\R\to\aut\big(\rsl\big)$ such that for any $\Lambda_0 \subset \Z$ 
and $R \in  {\cal R}(X_{\Lambda_0} , \sigma) $ we have
\[ 
\alpha_t (R) := \mbox{n--}\!\! \lim_{\Lambda \nearrow \Z} 
\pi_0^{-1} \big(  U_{\Lambda} (t) \pi_0(R) U_{\Lambda} (t)^{-1}
\big)
\] 
where $\mbox{n--}\! \lim$ denotes the norm limit.
\end{itemize}
\end{pro}
This result can be extended to lattice models in
higher dimensions. Moreover, the harmonic pinning
potentials $Q^2$ may be dropped from the local 
Hamiltonians without changing the result. {}For 
the modelling of interacting systems of indistinguishable Bosons 
one would have to  
replace the nearest neighbour interaction by a full two--body 
potential. Whereas the stability of the C*--algebras 
$ \pi_0 ({\cal R}(X_\Lambda, \sigma))$ under the 
action of the corresponding local dynamics 
can still be established in that case by the present methods, 
the existence of the thermodynamic limit $\Lambda \nearrow
\Z$ is not settled. We hope to return to
this problem elsewhere.

We would like to emphasize in conclusion 
that the action $\alpha:\R\to\aut\big(\rsl\big)$ is \textit{not} 
pointwise norm--continuous, \ie  we are not in the usual   
setting of C*--dynamical systems. However there is an 
important substitute: Let
${\cal K}_\Lambda$ be the compact ideal 
in ${\cal R}(X_\Lambda, \sigma)$, $\Lambda \subset \Z$ 
and consider the C*--algebra  
${\cal K} = C^*\{{\cal K}_\Lambda \, | \,  \Lambda \subset \Z
\} + {\C \un} 
\subset \rsl$. This algebra is a proper subalgebra of $\rsl$,
but it is weakly dense in  $\rsl$ in any regular representation.
(This fact together with the next result 
will be useful for the analysis of states
in the next subsection.) Then: 
\begin{pro} \label{ContinuousAction}
Let $V \in C_0(\R)$ be a real function 
and let $\{ H_\Lambda \}_{\Lambda \subset \Z}$
be the corresponding family of Hamiltonians. 
Then the action $\alpha:\R\to\aut\big(\rsl\big)$ defined above  
is pointwise norm--continuous
on the elements of ${\cal K} \subset \rsl$.
(Note that ${\cal K}$ is in general 
not stable under the action of $\alpha_\R$.) 
\end{pro}

It is an immediate consequence of this proposition that the
C*--algebra ${\cal L}$ generated by the algebras 
$\alpha_t({\cal K})$, $t \in \R$ is stable under the
action of  of $\alpha_\R$, and $\alpha:\R\to\aut\big(\rsl\big)$
acts pointwise norm continuously on ${\cal L}$. 
Thus ${\cal L}$ depends on the chosen dynamics
and will be different, for example, for lattice systems
with a next to nearest neighbour interaction. It seems that
there is no universal subalgebra of the resolvent
algebra which is stable and pointwise norm continuous
under the action of all possible dynamics. 

\subsection{Ground states of interacting dynamical systems}  \label{groundstates}

Having established the existence of a large family of 
interacting dynamics $\alpha_\R$ on $\rsl$ we want to
analyze now whether there exist corresponding states
of physical interest. We will focus on 
ground states; but our arguments also apply to 
 thermal equilibrium states. 
 It is clear from the outset 
 by Haag's theorem \cite{Emch} that we will have to
 leave the Fock representation $\pi_0$ in order to describe the
 desired states. 
As we will see, it is very simple 
to identify the correct representations and establish 
their desired properties in the present framework. 

Fix the potential $V \in C_0(\R)$ and let
$\{ H_\Lambda \}_{\Lambda \subset \Z}$ be the corresponding
set of local Hamiltonians. The operators 
$H_\Lambda$ differ from the Hamiltonians $H^{(0)}_\Lambda$
of the harmonic oscillator for $\Lambda$ 
by the bounded interaction potential which is an 
element of ${\cal R}(X_\Lambda, \sigma)$. 
Thus, since the restriction of the resolvent of 
 $H^{(0)}_\Lambda$ to ${\cal H}_0 (X_\Lambda)$
is a compact operator, so is the restriction 
of the resolvent of $H_\Lambda$. 
Hence, by the 
same reasoning as in Sec.\ \ref{affil}, it follows
that $H_\Lambda$ is an observable which is affiliated
with the compact ideal ${\cal K}_\Lambda \subset {\cal R}(X_\Lambda, \sigma)$. 
In particular, each $ H_\Lambda$ has 
discrete spectrum and, moreover, is bounded from below. 
It suffices to consider here 
the Hamiltonians $H_n := H_{\Lambda_n}$ corresponding to the 
sets $\Lambda_n := \{ l \in \Z \ | \ |l| \leq n \}$,  $n \in \N$.  
Let $E_n$ be the smallest eigenvalue
of $H_n$ and let  $\widetilde{H}_n := H_n - E_n \un$  
be the corresponding ``renormalized'' Hamiltonian 
whose smallest eigenvalue is $0$. 

Fix for each $n \in \N$ a 
normalized eigenvector $\Omega_n$ corresponding to the eigenvalue
$0$ of  $\widetilde{H}_n$, and define
a sequence of states $\{ \omega_n \}_{n \in \N}$ on
$\rsl$ by
\[ 
\omega_n(R) := \big(\Omega_n, \pi_0(R) \Omega_n\big) \, , \quad
R \in \rsl \, .
\]
This sequence need not converge, but by w*--compactness of the unit
ball of the dual of $\rsl$ it has limit points. Let 
$\omega_\infty$ be any such w*--limit point. Since 
$\omega_\infty$ is a limit of ground states of the local Hamiltonians
it is a candidate for a ground state of the full theory;
but for infinite systems such heuristic expectations are known to
be taken with a pinch of salt and to require careful analysis. 
Our first result shows that 
$\omega_\infty$ is a physically acceptable state. It is
essential for its proof that the local Hamiltonians 
are affiliated with the resolvent algebra. 

\begin{lem} \label{RegularState}
Let $\omega_\infty$ be any w*--limit point
of the sequence  $\{ \omega_n \}_{n \in \N}$ defined above.
Then $\omega_\infty$ is a regular state on $\rsl$.
\end{lem}
With this information we can 
prove that  $\omega_\infty$ is a ground state for the
dynamics  $\alpha:\R\to\aut\big(\rsl\big)$. We first show that
$\omega_\infty$ is invariant under the adjoint action of $\alpha$. 
Fix $R \in \rsl$ and $t \in \R$, then by
definition of  $\omega_\infty$ there is a subsequence
$\{\omega_{n_k} \}_{k \in \N}$ such that 
\[
\omega_\infty(\alpha_t(R) - R) = \lim_{k \rightarrow \infty}
\omega_{n_k}(\alpha_t(R)  - R) =
\lim_{k \rightarrow \infty} \big( \Omega_{n_k}, \pi_0(\alpha_t(R)  - R)
\Omega_{n_k} \big) \, .
\] 
By the results of the preceding
section we have the norm limits
\[ 
\pi_0(\alpha_t(R)) = \mbox{n--}\!\!\lim_{k
  \rightarrow \infty} ({\Ad}\, e^{\, it H_{{n_k}}})(\pi_0(R))
=  \mbox{n--}\!\!\lim_{k
  \rightarrow \infty} ({\Ad} \, e^{\, it
  \widetilde{H}_{{n_k}}})(\pi_0(R)) \, ,
\] 
where the second equality follows from the fact that the additive
renormalization term in $ \widetilde{H}_{{n_k}}$ drops out
in the adjoint action. But 
$ \big( \Omega_{n_k},  ({\Ad} \, e^{\, it
  \widetilde{H}_{{n_k}}})(\pi_0(R))  
\Omega_{n_k} \big) =  \big( \Omega_{n_k}, \pi_0(R)
\Omega_{n_k} \big) $ since $\Omega_{n_k}$ is a 
ground state vector for $\widetilde{H}_{{n_k}}$.
Hence $\omega_\infty(\alpha_t(R) - R) = 0$ 
for $R \in \rsl$, $t \in \R$. 

In the GNS--representation $(\pi_\infty, {\cal H}_\infty,
\Omega_\infty)$ of the $\alpha$--invariant state
$\omega_\infty$, we can now define as usual 
a unitary representation $U_\infty$ 
of $\R$ implementing the action $\alpha$ by
\[ 
U_\infty(t) \pi_\infty(R) \Omega_\infty :=
 \pi_\infty(\alpha_t(R)) \Omega_\infty \, , \quad R \in \rsl, \, t \in \R \, .
\]
It follows from Proposition \ref{ContinuousAction} that 
$U_\infty(\R)$ acts continuously on the subspace 
$\pi_\infty({\cal K}) \Omega_\infty \subset  {\cal H}_\infty$. 
But $\omega_\infty$ is a regular state, so this subspace is
dense in ${\cal H}_\infty$, proving that the representation 
$U_\infty$ is continuous in the strong operator topology. 

It remains to determine the spectrum of $U_\infty(\R)$.
Let $h \in {\cal S}(\R)$ be a test function
whose Fourier transform 
has support on the negative real axis and let $K \in {\cal
  K}$. By Proposition \ref{ContinuousAction} the integral 
$K(h) := \int \! dt \, h(t) \alpha_t(K)$ is defined in
the norm topology and hence an element of $\rsl$.
Picking any other $R \in \rsl$
there is a subsequence
$\{\omega_{n_k} \}_{k \in \N}$ such that 
\[
\omega_\infty(R K(h)) = \lim_{k \rightarrow \infty}
\omega_{n_k}(R K(h)) =
\lim_{k \rightarrow \infty} \big( \Omega_{n_k}, \pi_0(R K(h))
\Omega_{n_k} \big) \, .
\] 
Furthermore, making use once more of the results in the preceding 
section and the dominated convergence theorem we have 
\[
\pi_0(K(h)) = \int \! dt \, h(t) \pi_0(\alpha_t(K)) =
\mbox{n--}\!\!\lim_{k \rightarrow \infty} 
 \int \! dt \, h(t) ({\Ad} \, e^{\, it
  \widetilde{H}_{{n_k}}})(\pi_0(K)) \, .
\] 
Combining these facts and $ e^{\, - it
  \widetilde{H}_{{n_k}}} \, \Omega_{n_k} = \Omega_{n_k}$ we get 
\[
\omega_\infty(R K(h)) =
\lim_{k \rightarrow \infty}  \int \! dt \, h(t) \,
\big( \Omega_{n_k}, \pi_0(R)  e^{\, it
  \widetilde{H}_{{n_k}}} \pi_0(K) 
\Omega_{n_k} \big) = 0 \, ,
\]
where for the second equality we made use of the 
support properties of the Fourier transform of 
$h$ and the fact that each 
$ \widetilde{H}_{{n_k}}$ is a positive operator. 
Since $\omega_\infty(R K(h)) = \int \! dt \, h(t)
\omega_\infty(R \alpha_t(K))$ by Proposition \ref{ContinuousAction}
it follows that 
$ \int \! dt \, h(t) \, U_\infty(t) \! \upharpoonright  \! \pi_0({\cal
K}) \Omega_\infty = 0$. But $\pi_0({\cal
K}) \Omega_\infty$ is dense in ${\cal H}_\infty$, so 
this equality holds on the whole Hilbert space. 
We have thus shown that the generator of $U_\infty(\R)$
has spectrum on the positive real axis, proving that
$\Omega_\infty$ is a ground state vector. The preceding results 
are summarized in the next proposition.

\begin{pro}
Let $\omega_\infty$ be any w*--limit point
of the sequence  $\{ \omega_n \}_{n \in \N}$ constructed above.
Then $\omega_\infty$ is a regular ground state for the 
dynamics $\alpha:\R\to\aut\big(\rsl\big)$.
\end{pro}

Here we end our discussion of infinite
dynamical systems in the framework of the resolvent algebra.
There are many more intriguing questions which can be
addressed in the present model, such as the uniqueness and  purity
of $\omega_\infty$, its behaviour under lattice translations
(which are a symmetry of the present model) \etc\hspace*{-4pt}.
But we think that
the results presented so far provide sufficient evidence 
that the
resolvent algebra is a natural and powerful setting for the analysis of
infinite dynamical quantum systems.

\section{Constraint theory.}
\label{Qconstr}

In this section we illustrate the usefulness of
the resolvent algebra in the study of the
non--regular representations which occur naturally in some
applications of physical interest.
Nonregular representations of the CCR--algebra 
have been used in a number of
situations, cf.\ \cite{Ace}. In particular they occur in
the context of constraint theory, cf.\ \cite{GrHu3,GrLl},
and they have been used to circumvent indefinite metric
representations in Gupta--Bleuler electromagnetism, cf. \cite{GrLl}.
In this section, we wish to develop constraint theory for linear
bosonic constraints to see what form it takes in the 
resolvent algebra, and to investigate
the nonregular representations which occur.
Since the constraints eliminate nonphysical information, the nonregularity
of the representations is not a problem. One only needs regularity on
the final physical algebra.

{}For linear bosonic constraints, we start with a nondegenerate symplectic space
${(X,\sigma)}$ and specify a nonzero {\it constraint subspace} $C\subset X$.
Our task is to implement the heuristic constraint conditions 
\[
\j(f)\,\Psi=0\quad f\in C
\]
to select the subspace spanned by the physical vectors $\Psi$.
There are many examples where these occur, \eg in quantum electromagnetism,
cf. \cite{GrHu1, GrHu2, GrLl}.
Now in a representation $\pi$ of $\rsl$ for which $\ker\pi\big(\rlf\big)=\{0\}$
we have by Theorem~\ref{RegThm}(vi) that
$\pi(\rlf)\j_\pi(f)=i\lambda\pi(\rlf)-\un$ on $\dom \j_\pi(f)$.
Hence the appropriate form in which to impose the heuristic constraint condition
in the resolvent algebra is to select the set of physical (``Dirac'') states by
\begin{equation}
\label{Rconstraint}
\ot S._D:=\left\{\omega\in\ot S.(\rsl)\,\mid\,
\pi_\omega\big(i\lambda\rlf-\un\big)\,\Omega_\omega=0\quad 
 f\in C, \, \lambda \in \R \backslash 0  \right\} \, ,
\end{equation}
where $\pi_\omega$ and $\Omega_\omega$ denote the GNS--representation and 
GNS--cyclic vector of $\omega.$ Thus $\omega\in\ot S._D$ iff
$\al C.\subset
{\cal N}_\omega:={\big\{A\in\rsl\,\mid\,\omega(A^*A)=0\big\}}$,
where $\al C.:={\big\{i\lambda\rlf-\un\,\mid\, f\in
  C, \, \lambda \in \R \backslash 0 \big\}}$. Note that 
${\cal C}^* = {\cal C}$. 
\begin{pro}
\label{RDirac}
Given the data above, we have:
\begin{itemize}
\item[(i)] $\ot S._D=\left\{\omega\in\ot S.(\rsl)\,\mid\,
\omega\big(R(1,f)\big)=-i, \ f\in C\right\}$.
\item[(ii)] If $\omega\in\ot S._D$, then it is not regular. 
In particular, if $\sigma(g,C)\not=0$ for some $g\in X,$ 
then $\pi_\omega(R(\lambda,g))=0$ for all $\lambda\in\R\backslash 0$.
\item[(iii)]  $\ot S._D\not=\emptyset$ iff $\sigma(C,C)=0$.
\end{itemize}
\end{pro}
Henceforth we will assume that $\sigma(C,\,C)=0$ and hence $\ot S._D\not=
\emptyset$.
We examine the algebraic structures produced by these constraints,
cf.\ \cite{GrHu1,GrLl}. {}For the left ideal generated by the constraints
 we have $\al N.:=[\rsl\al C.]
=\bigcap\; \{\al N._{\omega}\mid\omega\in{\ot S.}_D \}$,
cf.\ Theorem 3.13.5 in \cite{Ped}, 
where ${[\, \cdot \,]}$ denotes the closed span of its argument.
Then there is the following structure, cf.\ \cite{GrLl}.
\begin{pro}
\label{Diracfacts}
\chop
Let $\al D.:=\al N.\bigcap\al N.^*$ and 
 ${\WD} := \{ F\in \rsl  \mid [F,\, \al D.] \subset {\al D.} \}.$ Then 
\begin{itemize}
\item[{\it(i)}] $\al D.=\al N.\bigcap \al N.^*$ is the unique maximal 
  {\rm C}$^*$--algebra in $\, \bigcap\; \{ {\rm Ker}\,\omega\mid \omega\in 
  {\got S}_{{D}} \}$. Moreover $\al D.$ is a hereditary 
  {\rm C}$^*$--subalgebra of $\rsl  $.
\item[{\it(ii)}] ${\WD} = 
\{ F\in{\rsl  }\mid F{\al D.}\subset{\cal D}\supset {\cal D}F \}$, 
  \ie it is the relative mul\-ti\-plier algebra of ${\al D.}$ in ${\rsl  }$.
\item[{\it(iii)}] $\WD=\{F\in\rsl  \mid\; [F,\,\al C.]\subset\al D.\}$.
\item[{\it(iv)}] $\al D.=[\WD \, \al C.]=[\al C. \, \WD]$.
\end{itemize}
\end{pro}
Then ${\WD}$ is the C*--algebraic analogue of Dirac's observables 
(the weak commutant of the constraints), and the algebra of physical
observables is $\Ob :=\WD/\al D.$, using the fact that 
$\al D.$ is a closed two--sided ideal of $\WD$. Note that by (iii)
the relative commutant $\al C.'$ of 
$\al C.$ in $\rsl$ 
(traditionally 
regarded as algebra of observables) is contained in $\WD$. 

Having introduced the general concepts, 
let us  determine these algebras more explicitly. If a $g\in X$
satisfies $\sigma(g,C)=0$, then $R(\mu,g)\in\al C.'\subset\WD$.
On the other hand, if $\sigma(g,C)\not=0$, then by 
Proposition~\ref{RDirac}(ii) we get that $\pi_\omega\big(R(\mu,g)\big)=0$
for all $\omega\in \ot S._D$. 
But then $AR(\mu,g)B\in {\cal N}_\omega\bigcap {\cal N}_\omega^*$ for all $\omega\in \ot S._D$
and $A,\,B\in\rsl$,
\ie $AR(\mu,g)B\in\al D.\subset\WD$. Thus all the generating elements of 
$\rsl$ are in $\WD$. As $\WD$ is a C*--algebra it follows that
$\WD = \rsl$, and hence that $\al D.$ is a proper ideal of $\rsl$. 
Then by Proposition~\ref{Diracfacts}(iv) we can write 
 $\al D.=[\rsl\al C.]=[\al C.\rsl]$.
Moreover, as any monomial in the resolvents containing a resolvent 
not in $\al C.'$
is in $\al D.,$ we conclude that any $A\in\rsl$ can be approximated in norm
by elements in $\al C.'$ modulo elements of  $\al D.$.
Thus $\Ob =\al C.'\big/(\al C.'\bigcap\al D.)$.
Thus, if $\pi_\omega$ is the GNS--representation of a Dirac state
$\omega\in\ot S._D$, then $\pi_\omega\big(\rsl\big)$ is a homomorphic image
of the traditional observables $\al C.'$.
To summarize, we have shown:
\begin{pro}
\chop
 $\WD=\rsl$ with the proper ideal  $\al D.=[\rsl\al C.]=[\al C.\rsl]$,
and $\Ob =\al C.'\big/(\al C.'\bigcap\al D.)$.
\end{pro}
So Dirac constraining of linear bosonic constraints is considerably simpler
in the resolvent algebra $\rsl$ than in the CCR--algebra $\CCRX$~cf.\cite{GrHu1}.

\section{Discussion}

Starting from the basic relations of the canonical
observables in quantum physics in resolvent form, we were led 
in a natural manner to an intriguing mathematical structure, the 
resolvent algebra. This C*--algebra has many desirable features
for the treatment of finite and infinite dimensional quantum systems:
It allows for the formulation of physically relevant 
dynamical laws, it contains a multitude of physically 
significant observables and it provides a powerful framework 
for the analysis of physical states. Moreover, it 
is a convenient setting for the discussion of
singular representations of the canonical observables which
appear naturally in the discussion of quantum constraints
and in the modelling of Fermionic symmetries  \cite{BuGr}  
in the context of supersymmetry or of 
BRST--constraint theory. 
 
The modest price for these conveniences is the 
 fact that the resolvent algebra has a 
non--trivial ideal structure (cf.\ the remark at the end of the proof of
Proposition \ref{dynamics}). These ideals correspond to the kernels of
representations of the resolvent algebra in which some of the
underlying fields have ``infinite values'' (and thus become physically
meaningless). The fact that the resolvents of these fields simply
disappear in the respective representations makes the representation 
theory of the resolvent algebra particularly simple 
(cf.\ Proposition \ref{Xdecomp}). On the other hand, the physically 
significant regular representations 
of the resolvent algebra are faithful,
and are thus in perfect agreement with the 
principle of physical equivalence \cite{HaKa}. 

The resolvent algebra competes with several other approaches
to the treatment of canonical quantum systems, such as 
the Weyl algebra, its twisted convolution form and other 
possible variants, cf.\ \cite{Ka}. Surprisingly, its   
existence seems to have escaped observation so far;  
for both from a conceptual and technical point
of view the resolvent algebra seems superior to these other approaches.
We could only illustrate here some of its many advantages, but our 
results suggest that further study and applications of the 
resolvent algebra to concrete quantum systems 
are worthwhile.

\section{Proofs}
\label{Proofs}

\subsection*{Proof of Proposition~\ref{CCRX}}
Assume that $M\in\CCRX$ is nonzero such that
$\phi (f)M$ is bounded for some nonzero
$f\in X$.
Let $W_t :=\exp(it\phi (f))$, $t \in \R$ and denote the spectral resolution
of $\phi (f)$ by $\phi (f)={\int\lambda\,dP(\lambda)},$ then
\begin{eqnarray*}
\left\|{(W_t -\un)M}\right\| &=&
\Big\|\int(e^{it\lambda}-1)dP(\lambda)M\Big\|  \\[1mm]
&=& |t| \, \Big\|\int{(e^{it\lambda}-1)\over t\lambda}\,dP(\lambda)
\int\lambda'\,dP(\lambda')M\Big\|   \\[1mm] 
&\leq& |t|\|\phi (f)M\|\longrightarrow 0
\end{eqnarray*}
as $t\to 0,$ where we used the bound
${|{e^{ix}-1\over x}|} \leq 1$. 
Let $\al J.\subset\CCRX$ consist of all elements $M$
such that $\left\|{(W_t - \un)M}\right\|\to 0$ as $t\to 0.$
This is clearly a  norm-closed linear space,
and by the inequality $\left\|{(W_t - \un)MA}\right\| 
\leq\left\|{(W_t -\un)M}\right\|\,\|A\|$ it is also 
a right ideal. To see that it is a two sided ideal note that
\[
\left\|{(W_t -\un)e^{i\j(g)}M}\right\|=
\left\|{(W_t e^{it\sigma(f,g)}-\un)M}\right\|
\]
still converges to $0$ as $t\to 0,$
and use the fact that $\CCRX$ is the norm closure of the span 
of $\set e^{i\j(g)},g\in{X}.$.
As $\CCRX$ is simple and $\| W_t  - 1 \| = 2, \, t \neq 0$, 
it follows that $\al J.\ni M$ is zero.

\subsection*{Proof of Proposition~\ref{Rbound}}
(i) Using $\rlf^*=R(-\lambda,f)$ and the C*--property of the norm, we get 
\[
2|\lambda|\|\pi(\rlf)\|^2=\|\pi(2\lambda\rlf\rlf^*)\|=\|\pi(\rlf-\rlf^*)\|
\leq 2\|\pi(\rlf)\|\,.
\]
Thus we get that either $\pi(\rlf)=0$ which implies $\|\pi(\rlf)\|\leq 1/|\lambda|$,
or $\pi(\rlf)\not=0$ in which case a cancellation gives 
$\|\pi(\rlf)\|\leq 1/|\lambda|$.
Since $\al R._0$ is generated polynomially by the elements $\rlf,$ it follows that for each 
$A\in\al R._0$ there is a $c_A \geq 0$ such that 
${\|\pi(A)\|}\leq c_A$ for all Hilbert space representations $\pi$ of $\al R._0.$\chop
(ii) Let ${\cal N}_\omega:=\{A\in\al R._0\,\mid\,\omega(A^*A)=0\}$,
then the image of the factor map $\xi:\al R._0\to\al R._0\big/{\cal N}_\omega$ 
is a pre--Hilbert space, equipped with the inner product
${(\xi(A),\xi(B))}:={\omega(A^*B)}$. The Hilbert closure of 
$\al R._0\big/{\cal N}_\omega$ is the GNS--space $\al H._\omega.$
The GNS--representation is defined on $\al R._0\big/{\cal N}_\omega$ by
$\pi_\omega(A)\,\xi(B)=\xi(AB)$ which is well--defined because
${\cal N}_\omega$ is a left ideal. It is clear that this is a *-representation
on the dense invariant domain $\al R._0\big/{\cal N}_\omega\subset\al H._\omega$
and that $\xi(\un)=:\Omega_\omega$ is a cyclic vector for it.\chop
Now by Eq.~(\ref{RRstar}) we have for
$\Psi\in\al R._0\big/{\cal N}_\omega$ that
\begin{eqnarray*}
\big\|\pi_\omega(\rlf)\Psi\big\|^2&=&\big(\pi_\omega(\rlf)\Psi,\,
\pi_\omega(\rlf)\Psi\big)=\big(\Psi,\,\pi_\omega(\rlf^*\rlf)\Psi\big) \\[1mm]
&=&\Big|(2i\lambda)^{-1}\big(\Psi,\,[\pi_\omega(\rlf)-\pi_\omega(\rlf)^*]\Psi\big)
\Big| \\[1mm]
&\leq&|\lambda|^{-1}\big\|\pi_\omega(\rlf)\Psi\big\|\cdot\|\Psi\| 
\end{eqnarray*}
by the Cauchy--Schwartz inequality. Thus 
$\big\|\pi_\omega(\rlf)\big\|\leq |\lambda|^{-1}$
and so $\pi_\omega$ is bounded.\chop

\subsection*{Proof of Theorem~\ref{Relemen}}
(i) By (\ref{Resolv}) we have that $i(\mu-\lambda)\rlf R(\mu,f)=
\rlf-R(\mu,f)=-\big(R(\mu,f)-\rlf\big)=i(\mu-\lambda)R(\mu,f)\rlf$,
\ie ${[\rlf,\,R(\mu,f)]}=0$.\chop
(ii) This follows directly from Eq.~(\ref{Rccr})
by interchanging $\lambda$ and $f$ with $\mu$ and $g$ resp.\chop
(iii) The Fock representation $\pi$ defines a bounded representation of $\al R._0,$ 
where the resolvents of the fields
give the Fock representation induced on
$\al R.(X,\sigma),$ \ie $\pi(\rlf)=(i\lambda-\phi_\pi   (f))^{-1}$.
Since this is bounded, it defines a unique representation of 
$\rsl$.  Now
\[
\|\rlf\|\geq\|\pi(\rlf)\|=\|(i\lambda-\phi_\pi   (f))^{-1}\|
=\sup_{t\in\sigma(\phi_\pi   (f))}\Big|{1\over i\lambda-t}\Big|=
{1\over|\lambda|}
\]
using the fact that the spectrum $\sigma(\phi_\pi(f)) = \R$,
cf.\ for example Chapter 3.1 in \cite{Emch}, in particular
the proof of Theorem 5. 
Thus by Proposition~\ref{Rbound}(i) one arrives at 
$\|\rlf\|=1/|\lambda|$.\chop
(iv) Rearrange Eq.~(\ref{Resolv}) to get: 
\[
\rlf\big(\un-i(\lambda_0-\lambda)R(\lambda_0,f)\big)=R(\lambda_0,f)\,.
\]
Now by (iii), if $\big|\lambda_0-\lambda\big|<\big|\lambda_0\big|$
then ${\big\|i(\lambda_0-\lambda)R(\lambda_0,f)\big\|}<1,$ and hence
${\big(\un-i(\lambda_0-\lambda)R(\lambda_0,f)\big)^{-1}}$ exists, and is
given by a norm convergent power series in ${i(\lambda_0-\lambda)R(\lambda_0,f)}$.
That is, we have that
\[
\rlf=R(\lambda_0,f)\big(\un-i(\lambda_0-\lambda)R(\lambda_0,f)\big)^{-1}
=\sum_{n=0}^\infty (\lambda_0-\lambda)^n\, i^n R(\lambda_0,f)^{n+1} 
\]
when $\big|\lambda_0-\lambda\big|<\big|\lambda_0\big|,$
as claimed. \chop
(v)
The map $\alpha_T\big(\rlf\big)=R(\lambda,Tf)$ permutes the generating
elements of the free unital *-algebra generated by ${\{\rlf\,\mid\,\lambda\in\R\backslash0,\;
f\in X\}}$, hence defines an automorphism of it. 
Since $T$ is symplectic, this automorphism preserves the relations~(\ref{Riden})
to (\ref{Rsum}), hence the ideal they generate in the free algebra, hence it factors through to an
automorphism of $\al R._0.$ Moreover 
since $\alpha$ is a *--automorphism, it
 maps $\ot S.$ to itself bijectively, hence it preserves the enveloping
C*--seminorm, and so defines an automorphism of $\rsl.$

\subsection*{Proof of Proposition~\ref{ImAut}}

The argument is similar to the last part of the preceding proof:
The map $\beta_h\big(R(z,f)\big) :=
R(z+ih(f),\,f)$ permutes the generating
elements of the free unital *-algebra  generated by
the set $\set{R(z,f)},z\in\C\backslash{i\R},\;f\in X.$ 
hence defines an automorphism of it. 
This automorphism preserves the relations~(\ref{Ridenz})
to (\ref{Rsumz}), hence the ideal they generate in the free algebra, hence 
it factors through to an
automorphism of the *--algebra obtained from factoring out this ideal.
Since $\beta_h$ is a *--automorphism, it
maps $\ot S.$ to itself bijectively, hence preserves the enveloping
C*--seminorm, and hence it defines an automorphism of $\rsl.$

\subsection*{Proof of Theorem \ref{Ideals0}}

Since the monomials $\prod\limits_{j=1}^kR(\lambda_j,f_j)$
span a dense subspace of $\rsl$, we get that
the right ideal $\big[\rlf\rsl\big]$ 
is the closed span of the monomials
${\rlf\prod\limits_{j=1}^kR(\lambda_j,f_j)}$. But by Eq.~(\ref{Rsum})
\[
\rlf \, R(\mu,g)=R(\mu,g)\,\rlf+i\sigma(f,g)\,\rlf\,R(\mu,g)^2\rlf\in
\rsl\rlf
\]
so in each monomial we can progressively move $\rlf$ past all the 
$R(\lambda_j,f_j)$ until it stands on the right, and hence
 $\big[\rlf\rsl\big]\subseteq\big[\rsl\rlf\big]$.
 Likewise, by using Eq.~(\ref{Rsum}) to move $\rlf$ to the left,
 we get the reverse inclusion and hence equality
 $\big[\rlf\rsl\big]=\big[\rsl\rlf\big]$.
 By a similar argument we also get $\big[\rlf\rsl\big]=\big[\rsl\rlf\rsl\big]$
 which is obviously the closed two--sided ideal generated by $\rlf$.

 Next, we need to prove that $\big[\rsl\rlf\big]$ is proper.
 If it is not proper, then $\un\in\big[\rsl\rlf\big]$ and hence there is a sequence
 $\{A_n\}\subset\rsl$ such that ${A_n\rlf}\to\un$ in norm.
 Let $\pi$ be the Fock representation on $\al H._\pi,$  then $\j_\pi(f):=
i\un-\pi(R(1,f))^{-1}$ is selfadjoint and
 $\pi(\rlf)={\int_\R(i\lambda-t)^{-1}\,dP(t)}$ where $dP$ indicates the spectral measure 
 of $\j_\pi(f)$. Choose $\Psi_k\in P[k,\,k+1]\al H._\pi$ with $\|\Psi_k\|=1$,
 then since the spectrum of $\j_\pi(f)$ is $\R$,
\begin{equation}
\label{psikEst}
\big\|\pi(\rlf)\,\Psi_k\big\|\leq\sup_{t\in[k,k+1]}\left|(i\lambda-t)^{-1}\right|
=(\lambda^2+k^2)^{-1/2}\,.
\end{equation}
Now $\pi\big(A_n\rlf\big)\to\un$ uniformly as $n\to\infty$, so
for $\varepsilon<1$ there is an $N$ such that for all $k>0$,
$\|\pi\big(A_n\rlf\big)\Psi_k\big\| \geq 1 - \varepsilon$.
By Eq.~(\ref{psikEst})
$\|\pi\big(A_n\rlf\big)\Psi_k\big\|\leq\big\|\pi(A_n)\big\|\big/\sqrt{
\lambda^2+k^2}$ and consequently 
$\|\pi(A_n)\big\|\geq (1-\varepsilon)\sqrt{\lambda^2+k^2}$. 
As $\varepsilon<1$ and $k$ arbitrary we conclude that $\pi(A_n)$ is
unbounded. But this is a contradiction, hence  $\big[\rsl\rlf\big]$ is proper.

{}For the last statement, consider the intersection
$\big[R(\lambda_1,f_1)\rsl\big]\cap\big[R(\lambda_2,f_2)\rsl\big]\ni A$.
Since $A\in\big[R(\lambda_1,f_1)\rsl\big]$ there is a sequence $\{B_n\}
\subset \rsl$ such that $R(\lambda_1,f_1)\,B_n\to A$ in norm.
Let $\{E_\iota \}$ be an approximate identity of $\big[R(\lambda_2,f_2)\rsl\big]$
then we can construct a sequence
$R(\lambda_1,f_1)\,B_n\,E\s\iota_n.\to A$ in norm, using the fact
that $A\in\big[R(\lambda_2,f_2)\rsl\big]$. 
Since $\{E_\iota\}\subset\big[R(\lambda_2,f_2)\rsl\big]$ we can find 
$F_n\in\rsl$ such that $F_nR(\lambda_2,f_2)$ are arbitrarily close
to $E\s\iota_n.$ and hence we can find such $F_n$ such that
the sequence 
$R(\lambda_1,f_1)\,B_n\,F_n\,R(\lambda_2,f_2)\to A$ in norm,
and hence $A\in\big[R(\lambda_1,f_1)\,R(\lambda_2,f_2)\rsl\big]$.
Since it is trivial that 
\[
\big[R(\lambda_1,f_1)\,R(\lambda_2,f_2)\rsl\big]
\subseteq \big[R(\lambda_1,f_1)\rsl\big]\cap\big[R(\lambda_2,f_2)\rsl\big]
\]
the equality now follows. {}For more intersections we repeat an
inductive version of the argument.

\subsection*{Proof of Theorem \ref{Ideals1}}

(i)
Let $K:=\ker\pi\big(\rlf\big)$ and note that from Eq.~(\ref{Rccr})
\[
\pi\left(\rlf\,R(\mu,g)\right)=
\pi\left(R(\mu,g)+i\sigma(f,g)\,\rlf\,R(\mu,g)^2\right)\,
\pi(\rlf)
\]
hence $\pi\left(\rlf\,R(\mu,g)\right)\,K=0$ and so 
$\pi(R(\mu,g))\,K\subseteq K$ for all $\mu$ and $g$.
Since $R(\mu,g)^*=R(-\mu,g),$ the subspace $K$ reduces
$\pi(\rsl)$. Put $\pi_1:=\pi\restriction K$ and
$\pi_2:=\pi\restriction K^\perp$ to obtain the desired
decomposition. Uniqueness is clear.
\chop
(ii) According to (i)  $\pi=\pi_1\oplus\pi_2$ 
where $\pi_1(\rlf)=0$ and
$\ker\pi_2(\rlf)=\{0\}$. Thus by Theorem~\ref{RegThm}(ii), to be proved
subsequently, 
we obtain $\slim\limits_{\lambda\to\infty}i\lambda\pi(\rlf)
=\un\s{K}^\perp. =:P_f \in {\pi\big(\rsl\big)}^{\prime \prime}$, 
and it is obvious that $P_f$
commutes with ${\pi\big(\rsl\big)}^{\prime \prime}$.
Since $\pi_2(\rlf))$ is invertible, it has dense range, so the closure
of its range is $K^\perp$ and the projection on this is $P_f$.
Moreover 
\[
\pi\left(\big[\rlf\,\rsl\big]\right)\al H.
=\pi_2\left(\big[\rlf\,\rsl\big]\right) K^\perp 
\]
and this is dense in $K^\perp$ because
$\rlf\in\big[\rlf\,\rsl\big]$.
\chop
(iii) If $\pi$ is factorial, its center is trivial, so its
central projections can only be $0$ or $\un$. 
\chop
(iv)
By Theorem \ref{Ideals0},
$\big[R(\lambda,f)\rsl\big]$ is a proper ideal of $\rsl$, so
any state of the factor algebra ${\rsl\big/\big[R(\lambda,f)\rsl\big]}$
lifts to a state $\omega$ of $\rsl$ with  
$R(\lambda,f)\in\big[R(\lambda,f)\rsl\big]$
in its kernel.
On the other hand, if we are given a state $\omega$ of $\rsl$ with
$R(\lambda,f)\in\ker\omega$ then also $R(\lambda,f)^*\in\ker\omega$
and hence by $\rlf-\rlf^*=-2i\lambda\rlf\rlf^*$
we get that $\rlf\rlf^*\in\ker\omega$, \ie $\rlf\in {\cal N}_\omega^*$
which is a right ideal. Thus
$\big[R(\lambda,f)\rsl\big]\subset {\cal N}_\omega^*\subset\ker\omega$,
and as $\big[R(\lambda,f)\rsl\big]$ is a two--sided ideal
it must be in $\ker\pi_\omega$.

\subsection*{Proof of Theorem~\ref{RegThm}}

(i) Observe that by Theorem 1 in \cite[p 216]{Yos}, we deduce from 
${\ker\pi(R(1,f))}=\{0\}$ that $\pi(\rlf)$ is the resolvent of $\j_\pi(f),$
\ie we have now for all $\lambda\not=0$ that
$\j_\pi(f)=i\lambda\un-\pi(\rlf))^{-1}$. Then
\begin{eqnarray*}
\j_\pi(\mu f) &=& i\un -\pi(R(1,\mu f))^{-1}=i\un- \mu \pi(R(\f 1, \mu.,f))^{-1} \\[1mm]
&=& \mu \left(i {\f 1, \mu.} \un -\pi\big(R(\f 1,\mu.,f)\big)^{-1}\right)
=\mu\,\j_\pi(f)\,.
\end{eqnarray*}
Thus
\begin{eqnarray*}
\j_\pi(f)^* &=& {\left(i\un-\pi(R(1,f))^{-1}\right)^*}
\supseteq -i {\un} - \left(\pi(R(1,f))^{-1}\right)^*\\[1mm]
&=&-i\un -\pi(R(1,f)^*)^{-1}
=-i\un-\pi(R(-1,f))^{-1} \\[1mm]
&=&-i\un +\pi(R(1,-f))^{-1} = -\j_\pi(-f)=\j_\pi(f)
\end{eqnarray*}
and hence $\j_\pi(f)$ is symmetric. To see that it is selfadjoint
note that:
\[
\ran\left(\j_\pi(f)\pm i\un\right)=\ran\left(-\pi\big(R(\pm 1,f)\big)^{-1}\right)
=\dom\left(\pi\big(R(\pm 1,f)\big)\right)=\al H._\pi
\]
hence the deficiency spaces $\left(\ran\left(\j_\pi(f)\pm i\un\right)\right)^\perp
=\{0\}$ and so $\j_\pi(f)$ is selfadjoint.\chop
{}For the domain claim, recall that 
$\dom \j_\pi(f)=\ran\pi\big(R(1,f)\big)$.
So 
\begin{eqnarray*}
\pi(\rlf)\dom \j_\pi(h)&=&\pi(\rlf)\pi\big(R(1,h)\big)\al H._\pi  \\[1mm]
&=&\pi\left(R(1,h)\rlf+i\sigma(f,h)R(1,h)\rlf^2R(1,h)\right)\al H._\pi \\[1mm]
&\subseteq&\pi\left(R(1,h)\right)\al H._\pi=\dom \j_\pi(h).
\end{eqnarray*}
(ii) Let $\j_\pi(f)=\int\mu \, dP(\mu)$ be the spectral resolution of $\j_\pi(f).$
Then $\pi(R(\lambda,f))=\int{1\over i\lambda -\mu}dP(\mu)$ hence
\[
i\lambda\pi(R(\lambda,f))\Psi=
\int{i\lambda \over i\lambda -\mu}\,dP(\mu)\Psi, \quad \Psi
\in\al H._\pi\,.
\]
Since $\left|{i\lambda \over i\lambda -\mu}\right|<1$ 
(for $\lambda \in\R\backslash 0)$
 the integrand is dominated by $1$ which is an $L^1\hbox{--function}$
with respect to $dP(\mu),$ and as we have pointwise that
$\lim\limits_{\lambda \to\infty}{i\lambda \over i\lambda -\mu} =1,$ we can apply
the dominated convergence theorem to get that
\[
\lim_{\lambda \to\infty}i\lambda \pi(R(\lambda ,f))\Psi=
\int dP(\mu)\Psi=\Psi\,.
\]
(iii) $i\pi(R(1, \mu f))\Psi
=\int{i\over i- \mu \lambda}\,dP(\lambda)\Psi\to\Psi$ as $\mu \to 0$ by the same
argument as in (ii)\,.\chop
(iv) Let $\al D.:={\pi\big(R(1,f)R(1,h)\big)\al H._\pi},$ then by definition
$\al D.\subseteq\ran\pi(R(1,f))=\dom \j_\pi(f)$.
Moreover ${\pi\big(R(1,f)R(1,h)\big)\al H._\pi} =
{\pi\big(R(1,h)[R(1,f)+i\sigma(f,h)R(1,f)^2R(1,h)]\big)\al H._\pi}\subseteq\ran\pi(R(1,h))=
\dom \j_\pi(h),$
\ie $\al D.\subseteq\dom \j_\pi(f)\cap\dom \j_\pi(h)$.
That $\al D.$ is dense, follows from (iii) of this theorem, using
\[
\lim_{\mu \to 0}\lim_{\nu \to 0}\pi\big(R(1,\mu f)R(1, \nu h)\big)\Psi=-\Psi
\]
for all $\Psi\in\al H._\pi$, as well as $\mu R(1, \mu f)=R(1/\mu ,\,f)$
and the fact mentioned before that
all $\pi(\rlf)$ have the same range for $f$ fixed.

Let $\Psi\in\al D.,$ \ie $\Psi=\pi\big(R(1,f)R(1,h)\big)\Phi   $ for some
$\Phi   \in\al H._\pi$. Then
\begin{eqnarray*}
& &\pi\big(R(1,h)R(1,f)\big)\big[\j_\pi(f),\j_\pi(h)\big]\Psi  \\[1mm]
& &=\pi\big(R(1,h)R(1,f)\big)\big[\pi(R(1,f))^{-1},\pi(R(1,h))^{-1}\big]
\pi\big(R(1,f)R(1,h)\big)\Phi     \\[1mm]
& & = \pi\big(R(1,f)R(1,h)-R(1,h)R(1,f)\big)\Phi   
= i\sigma(f,h)\pi\big(R(1,h)R(1,f)^2R(1,h)\big)\Phi    \\[1mm]
& & =  i\sigma(f,h)\pi\big(R(1,h)R(1,f)\big)\Psi\;.
\end{eqnarray*}
Since $\ker\pi\big(R(1,h)R(1,f)\big)=\{0\}$ it follows that 
$\big[\j_\pi(f),\j_\pi(h)\big]=i\sigma(f,h) \un$ on $\al D.$.\chop
(v) {}From Eq.~(\ref{Rhomog}) we have that
\[
\pi(R(\nu, f))=(i\nu \un -\j_\pi(f))^{-1}=\f 1,\nu. \, \pi\big(R(1,\f 1,\nu.f)\big)
=\f 1,\nu.\, \left(i\smash{\un}-\j_\pi(\f 1,\nu.f)\right)^{-1}
\]
and hence that $\j_\pi(f)=\nu \,\j_\pi(\f 1,\nu.f),$ \ie
$\j_\pi(\nu f)=\nu \j_\pi(f)$ for all $\nu \in\R\backslash 0$
and hence the claim is established for $h=0,$ and we only need to
prove it for $\nu =1$. 
Consider Eq.~(\ref{Rsum}):
\[
\pi\left(\rlf R(\mu,h)\right)= \pi\left(R(\lambda+\mu,\,f+h)[\rlf+R(\mu,h)
+i\sigma(f,h)\rlf^2R(\mu,h)]\right)
\]
and note that as $K:=\ker\pi\big(R(1,f+h)\big)$ reduces $\pi(\rsl)$,
it is also in the kernel of the left hand side ${\pi\left(\rlf R(\mu,h)\right)}$ 
of the equation. However the latter is invertible, hence $K=\{0\}$
(thus the term in square brackets on the rhs is also invertible).
It is also clear from the equation that 
$\dom\j_\pi(f+h)=\ran\pi\left(R(\lambda+\mu,\,f+h)\right)
\supset\al D.$ which is the range of the left hand side;
moreover, the invertibility of the term in the square 
brackets implies that $\al D.$ is a core for the 
selfadjoint operator $\j_\pi(f+h)$. \chop
Next we multiply the equation above  
on the left by $i(\lambda+\mu)\un-\j_\pi(f+h)$ and apply this to
${(i\mu\un-\j_\pi(h))}{(i\lambda\un-\j_\pi(f))}\Psi,$ $\Psi\in\al D.$ to get
\[
\big(i(\lambda+\mu)\un-\j_\pi(f+h)\big)\Psi
=\left((i\mu\un-\j_\pi(h))+(i\lambda\un-\j_\pi(f))\right)\Psi
\]
making use of $\left[(i\mu\un-\j_\pi(h)) ,  
(i\lambda\un-\j_\pi(f))\right]\Psi=i\sigma(h,f)\Psi$. The additivity
of $\phi_\pi$ on $\al D.$ then follows. \chop
{}For the proof that  ${\pi\big(R(1,\nu f+h)\big)}$ is contained in
the von Neumann algebra generated by  
${\big\{\pi\big(R(1,f)\big),\,\pi\big(R(1,h)\big)\big\}}$, we begin by 
showing that the commutants of $\pi\big(R(\lambda, \nu f)\big)$ coincide
for all $\lambda, \nu \in \R \backslash 0$. Since 
$\pi\big(R(\lambda, \nu f)\big) = {1 \over \nu} 
\pi\big(R(\lambda / \nu, f) \big)$
it suffices to establish this for the case $\nu = 1$. Let $A \in {\al
  B.}({\al H.}_\pi)$ such that 
$[A, \pi\big(R(\lambda,f)\big) ] = 0$. Then it follows from
(4) that 
\[ 
[A, \pi\big(R(\mu ,f)\big) ] \, 
\big( 1 + i(\mu - \lambda) \pi\big(R(\lambda,f)\big) \big)=
\big[A,\pi( R(\lambda,f))\big]=0 \, , \qquad
\mu \in \R \backslash 0 \, .
\] 
Now the spectrum of $\pi\big(R(\lambda,f)\big)$ is contained in 
$\{ (i\lambda - x)^{-1} \mid x \in \R \}$, so 
$\big( \un + i(\mu - \lambda) \pi\big(R(\lambda,f)\big) \big)$ 
has dense range, hence $[A, \pi\big(R(\mu,f)\big) ] = 0$.
The same result clearly holds also for the resolvents
$\pi\big(R(\lambda, \nu h)\big)$, $\lambda, \nu \in \R \backslash 0$. 
To complete the proof, it therefore suffices to show that the commutant of 
${\pi\big(R(\lambda + \mu ,f+h)\big)}$ contains the commutant of 
$\big\{\pi\big(R(\lambda ,f)\big),\,\pi\big(R(\mu ,h)\big)\big\}$ for 
$\lambda, \mu, \lambda + \mu \in \R \backslash 0$. 
Let $A \in\al B.(\al H._\pi)$ commute with both $\pi\big(R(\lambda,f)\big)$
and $\pi\big(R(\mu,h)\big).$ Apply the commutator with $A$ to both sides of
the represented Eq.~(6) to get:
\[
0=\big[A,\,\pi\left(R(\lambda+\mu,\,f+h)\right)\big]\,\pi\Big(\rlf+R(\mu,h)
+i\sigma(f,h)\rlf^2R(\mu,h)\Big) \, .
\]
As mentioned above, 
$\pi\big(\rlf+R(\mu,h) +i\sigma(f,h)\rlf^2R(\mu,h)\big)$ is invertible. 
So, by the preceding results, 
$0=\big[A,\,\pi\left(R(1, \nu f+h)\right)\big]$, proving  
${\pi\big(R(1,\nu f+h)\big)}\in
 {\big\{\pi\big(R(1,f)\big),\,\pi\big(R(1,h)\big)\big\}''.}$
 
\noindent 
(vi) {}From the spectral resolution for $\j_\pi(f)$ we have trivially that on $\dom \j_\pi(f)$
\[
\j_\pi(f)\pi(R(\mu,f))=\pi(R(\mu,f))\j_\pi(f)=
\int{\lambda\over i\mu -\lambda}dP(\lambda)
=i\mu\pi(R(\mu,f))-\un\,.
\]
(vii) Let $\Psi\in\dom \j_\pi(f)=\ran\pi(\rlf),$ \ie $\Psi=\pi(\rlf)\Phi   $ for some
$\Phi   \in\al H._\pi.$ Then
\begin{eqnarray*}
\pi(\rlf)\big[\j_\pi(f),\,\pi(R(\lambda,h))\big]\Psi & \! = \! &
\pi(\rlf)\big[\j_\pi(f),\,\pi(R(\lambda,h))\big]\pi(\rlf)\Phi     \\[1mm]
& \! = \! & \pi\left(\big[\rlf, R(\lambda,h) \big]\right)\Phi   
=i\sigma(f,h)\pi\big(\rlf R(\lambda,h)^2\rlf\big)\Phi     \\[1mm]
& \! = \! & i\sigma(f,h)\pi\big(\rlf R(\lambda,h)^2\big)\Psi\;.
\end{eqnarray*}
Since $\ker\pi(\rlf)=\{0\}$, it follows that
\[
\big[\j_\pi(f),\,\pi(R(\lambda, h))\big]= i\sigma(f,h)\pi\big( R(\lambda,h)^2\big)
\]
on $\dom \j_\pi(f)$.\chop
(viii) We first prove the second equality.
Let $\Psi,\, \Phi    \in \widetilde{\al D.}
:={\rm Span}\set{\chi\s[-a,a].\big(\j_\pi(f)\big) \al H._\pi},a>0.$
where $\chi\s[-a,a].$ indicates the characteristic function of 
${[-a,a]},$ and note that
$\widetilde{\al D.}$
is a dense subspace.  
Since
${\left\|\j_\pi(f)^n\restriction\chi\s[-a,a].\big(\j_\pi(f)\big)\al 
H._\pi\right\|}
\leq a^n$ for $n\in\N,$ we can use the exponential series, \ie
\[
W(f)\Psi:=\exp\big(i\j_\pi(f)\big)\Psi
=\sum_{n=0}^\infty{\big(i\j_\pi(f)\big)^n\over n!}\Psi,  
\quad \Psi\in \widetilde{\al D.} \;.
\]
By the usual rearrangement of series we then have
\[
\left(\Phi   ,\,W(f)\pi\big(R(\lambda,h)\big)W(f)^{-1}
\Psi\right)=
\sum_{n=0}^\infty{1\over n!}\left(\Phi   ,\,\big({\rm ad}\,i\j_\pi(f)\big)^n\big(\pi(R(\lambda,h))
\big)\Psi\right)
\]
for all $\Phi   ,\,\Psi\in\widetilde{\al D.}$. Using part (vii) 
repeatedly we have 
\[
\big({\rm ad}\,i\j_\pi(f)\big)^n\big(\pi(R(\lambda,h))
= n!\,\sigma(h,f)^n \, \pi\big(R(\lambda,h)^{n+1}\big) \\[1mm]
\]
and consequently 
\begin{eqnarray*}
 \left(\Phi, (\Ad \,W(tf)) (\pi\big(R(\lambda,h)\big) \Psi\right)
& = & \sum_{n=0}^\infty t^n\sigma(h,f)^n\big(\Phi   ,\,\pi(R(\lambda,h)^{n+1})
\Psi\big)  \\[1mm]
&  = & \left(\Phi   ,\,\pi(R(\lambda+it\sigma(h,f),\,h))
\Psi\right) 
\end{eqnarray*}
whenever $\big|t\sigma(h,f)\big|<|\lambda|$ and
where we made use of the von Neumann series (Theorem~\ref{Relemen}(iii)) in the last step.
Since the operators involved are bounded and $\widetilde{\al D.}$ 
is dense, it follows that
$W(tf)\pi\big(R(\lambda,h)\big)W(tf)^{-1}=\pi(R(\lambda+it\sigma(h,f),\,h))$
for $\big|t\sigma(h,f)\big|<|\lambda|.$ By analyticity in $\lambda$ this can be
extended to complex $\lambda$ such that $\lambda\not\in i\R.$ Using the group property of
$t\mapsto W(tf)$ we then obtain for $\lambda\in\R\backslash 0$ that
\begin{equation}
\label{AdWR}
W(f)\pi\big(R(\lambda,h)\big)W(f)^{-1}=\pi(R(\lambda+i\sigma(h,f),\,h))\;.
\end{equation}
To prove the first equation, let us write $W(h)$ in terms of resolvents.
Note that $\lim\limits_{n\to\infty}(1-it/n)^{-n}=e^{it},$ $t\in\R$ and
since  
$\sup\limits_{t\in\R}\left|(1-it/n)^{-n}\right|=1,$ it follows from spectral theory 
(cf.\ Theorem~VIII.5(d) in \cite[p 262]{RS1}) that
\[
W(h)=e^{i\j_\pi(h)}=\lim_{n\to\infty}\left(1-i\j_\pi(h)/n\right)^{-n}
=\lim_{n\to\infty}\pi\left(iR(1,-h/n)\right)^n
\]
in strong operator topology. Apply Eq.~(\ref{AdWR}) to this to get
\begin{eqnarray*}
W(f)W(h)W(f)^{-1} &=& 
\slim_{n\to\infty}\pi\left(iR(1+i\sigma(-\f h,n., f),-\f h,n.)
\right)^n \\[1mm]
&=& \slim_{n\to\infty} \Big( \un- (i/n) \big(\sigma(h,f) \un + 
\j_\pi(h) \big)  \Big)^{-n}  \\[1mm]
&=&\exp\big(- i\sigma(f,h) \un  + i\j_\pi(h) \big) = e^{-i\sigma(f,h)}
\, W(h) \, .
\end{eqnarray*}
Making repeatedly use of this equation 
and Theorem \ref{RegThm}(v) 
the asserted Weyl relations then follow by an application of
the Trotter product formula, cf.\ Theorem VIII.31 in \cite[p.\ 297]{RS1}, 
\begin{eqnarray*}  W(f+h) & = & \slim_{n \rightarrow \infty}
\Big( W({\textstyle {1 \over n}} f) 
 W({\textstyle {1 \over n}} h)  \Big)^n \\
& = & \slim_{n \rightarrow \infty}
e^{i ((n^2 -n)/2 n^2) \, \sigma(f,h)} \, W(f) W(h) 
= e^{i \sigma(f,h)/2 } W(f) W(h) \, .
\end{eqnarray*} 
Finally, 
\[
W(sf)\al D.=W(sf)\pi\big(\rlf R(\mu,h)\big)\al H._\pi
=\pi\big(\rlf R(\mu+i\sigma(h, sf),h)\big)\al H._\pi \subset \al D.
\]
hence we conclude that $\al D.$ is a core for $\j_\pi(f)$ (cf.  Theorem~VIII.11
in \cite[p 269]{RS1}).

\subsection*{Proof of Corollary~\ref{RegBij}}

Given $\pi\in{\rm Reg}\big(\rsl,\al H.\big),$ then by definition
$\wt\pi$ is regular on $\CCRX.$  To see that the correspondence $\pi\mapsto\wt\pi$
is a bijection, we verify that Eq.~(\ref{Laplace1}) defines its inverse.
This is obvious, because for $\lambda > 0$ one obtains 
by spectral theory 
\[
-i\int_0^\infty e^{-\lambda t}\wt\pi(\delta\s -tf.)\,dt
=-i\int_0^\infty e^{-\lambda t}\exp\big(-it\j_\pi(f)\big)\,dt=\big(i\lambda\un-\j_\pi(f)\big)^{-1}
=\pi(\rlf)\, , 
\]
and similarly for $\lambda < 0$.
It is also clear from the definition that $\pi\mapsto\wt\pi$ respects direct sums,
and since $\big\{\exp\big(i\j_\pi(f)\big)\,\mid\,f\in X \big\}''=
\big\{\big(i\lambda\un-\j_\pi(f)\big)^{-1}\,
\mid\,f\in X,\,\lambda\in\R\backslash 0\}''$
we see that it takes irreducible representations to irreducibles, and conversely.
The corresponding bijection $\omega\mapsto\wt\omega$ for regular states is given by 
\[
\wt\omega(A):=\big(\Omega_\omega,\,\wt\pi_\omega(A)\Omega_\omega\big),\quad A\in\CCRX
\]
from which the claims follow.

\subsection*{Proof of Proposition~\ref{RegAlg}}

(i) If $\pi$ is faithful and factorial, then $P_f=\un$ for all $f\in X\backslash 0$,
or else $\pi(\rlf)=0$ by Theorem~\ref{Ideals1}(iii) for some $f$ which contradicts
with the faithfulness of $\pi$.
But by Theorem~\ref{Ideals1}(i),  $P_f$ is the projection onto
 ${\big(\ker\pi\big(R(1,f)\big)\big)^\perp}$
so $P_f=\un$ implies ${\ker\pi\big(R(1,f)\big)}=\{0\}$ and as this holds
for all $f$ it follows that $\pi$ is regular.
\chop
(ii) The same argument as in the proof of 
Theorem~\ref{Relemen}(iii) will apply here
because if $\pi$ is regular, then the spectrum of each $\j_\pi(f)$ is  
all of $\R$, and this was the only property of 
the Fock representation used in the proof of 
Theorem~\ref{Relemen}(iii).\chop
 (iii) 
 $\lim\limits_{\lambda\to\infty}i\lambda\,\omega\big(\rlf C\big)
 =\big(\Omega_\omega,\,\slim\limits_{\lambda\to\infty}i\lambda\,
\pi_\omega\big(\rlf C\big)
 \Omega_\omega\big)=\big(\Omega_\omega,P_f\pi_\omega(C)\Omega_\omega\big)$
 for $C\in\rsl$, where we made use of
Theorem~\ref{Ideals1}(ii). Since $P_f$ commutes with
$\pi_\omega(\rsl)$,  we see that for $A,B \in \rsl$
 \[ \big(\pi_\omega(B)\Omega_\omega, \pi_\omega(A)\Omega_\omega\big)
 = \omega(B^*A) = 
 \lim_{\lambda\to\infty}i\lambda\,\omega\big(\rlf B^*A\big)
 =\big(\pi_\omega(B)\Omega_\omega,P_f\pi_\omega(A)\Omega_\omega\big) \, .
 \]
 This implies that $P_f=\un$,  which in turn implies via 
Theorem~\ref{Ideals1}(i) and (ii)
 that $\pi_\omega(\rlf)$ is invertible. Since this holds for all $f\in X\backslash 0$
 it means that $\pi_\omega$ is regular. The converse is trivial.
 
\subsection*{Proof of Proposition \ref{Xdecomp}}

(i)
We have $X_R:=\set f\in X,\ker\pi\big(R(1,f)\big)=\{0\}.$. Let
$f,\;g\in X_R$ then by
Theorem~\ref{RegThm}(v) we have that
$\ker\pi\big(R(1,\nu f +g)\big)=\{0\}$ for $\nu \in\R,$ \ie 
$\nu f +g\in X_R$ and thus $X_R$ is a linear space.\chop
Now let $f\in X_S:=X\backslash X_R$ and $g\in X_R$.
If $f+g\in X_R$ then obviously $f\in X_R$ which is a contradiction,
so $f+g\in X_S$.\chop
(ii)
Recall that
$X_T:=\set f\in X,\ker\pi\big(R(1,f)\big)=\{0\}\;\hbox{and}\;\pi\big(R(1,f)\big)^{-1}\in
{\al B.(\al H.)}.\subset X_R$. 
Let $f,\;g\in X_T$ then by  Theorem~\ref{RegThm}(i) there are selfadjoint operators
$\j_\pi(f)$ and $\j_\pi(g)$ such that
 $\pi(\rlf)={\big(i\lambda\un-\j_\pi(f)\big)^{-1}}$ and
 $\pi(R(\mu,g))={\big(i\mu\un-\j_\pi(g)\big)^{-1}}$, and by definition of
$X_T$ both field operators are bounded. Then by Theorem~\ref{RegThm}(v) we have that
${\j_\pi(\nu f +g)}=\nu 
\j_\pi(f)+\j_\pi(g)$ (hence it is bounded) for $\nu \in\R$,
and thus ${i\un-\j_\pi(\nu f +g)}={\big(R(1,\nu f +g)\big)^{-1}}$ is bounded,
\ie $\nu f +g\in X_T$, and so $X_T$ is a linear space.\chop
Let  $f\in X_T$ and $g\in X$ with $\sigma(f,g)\not=0$.
Since $f\in X_T$ we can set $c:={\big\|\pi\big(R(1,f)\big)^{-1}\big\|}$.
Then 
\begin{equation} \label{Normbracket} 
\big\| \big[\pi\big(R(1,f)\big)^{-1} ,\pi\big(R(1,g)\big)^{n} \big]
\big\| \leq 2c \, \big\| \pi\big(R(1,g)\big)^{n} \big\| \, . 
\end{equation} 
On the other hand,
since $\;{\rm ad}\,\pi\big(R(1,f)\big)^{-1}=-{\rm ad}\,
\j_\pi(f)$ is an inner derivation of 
$\al B.(\al H.)$, we obtain from Theorem~\ref{RegThm}(vii) 
\[
\Big[\pi\big(R(1,f)\big)^{-1},\,\pi\big(R(1,g)\big)^n\Big]
=-in\,\sigma(f,g)\,\pi\big(R(1,g)^{n+1}\big)\,.
\]
Take the norm of this and use inequality~(\ref{Normbracket}) to find
for all $n\in\N$ that
$$
\label{Normgn}
n\,|\sigma(f,g)|\,\big\|\pi\big(R(1,g)\big)^{n+1}\big\|
\leq 2c\,\big\|\pi\big(R(1,g)\big)^n\big\|\;.
$$
Now assume that  $\pi\big(R(1,g)\big)\not=0$. Then
since it is normal (cf. Theorem~\ref{Relemen}(i)) we have: 
$$
\big\|\pi\big(R(1,g)\big)^n\big\| = \big\|\pi\big(R(1,g)\big) \big\|^n
\neq 0
\, , 
$$
and consequently 
$$
n\,|\sigma(f,g)| 
\leq 2c\,\big\|\pi\big(R(1,g)\big)\big\|^{-1}, \quad n \in \N \, .
$$
However this is impossible because $\sigma(f,g)\not=0$, 
hence $\pi\big(R(1,g)\big)=0$. \chop
(iii)
If $\pi$ is factorial, then by Theorem~\ref{Ideals1}(iii) for any $f\in X$
the projection onto $\ker\pi\big(R(1,f)\big)$ is $0$ or $\un$.  
Since $f\in X_S$ implies that $\ker\pi\big(R(1,f)\big)\not=0$ it follows that 
the projection onto the kernel is $\un$, \ie $\pi\big(R(1,f)\big)=0$ for
$f\in X_S$. \chop
From part~(ii) we have that $X_T\subset X_R^\perp$, and by the preceding step 
$\pi\big(R(1,X_S)\big)=0$. Thus by $X=X_R\cup X_S$ it follows that
any $\pi\big(R(1,f)\big)$ with $f\in X_T$ commutes with all  $\pi\big(R(1,g)\big),$ $g\in X$,
hence with $\pi\big(\rsl\big)$.
Since $\pi$ is factorial, its centre is trivial, hence
$\pi\big(R(1,f)\big)$ is a multiple of the identity.\chop
 Finally, 
by (ii) we know that ${\sigma(X_T,X_R)}=0$, \ie $X_T\subseteq{X_R\cap X_R^\perp}$.
Conversely, let $h\in{X_R\cap X_R^\perp}$ then $\pi(R(\lambda,h))$ commutes with
$\pi(R(\mu,X_R))$ and as $\pi(R(\lambda,X_S))=0$ it commutes with 
$\pi(\rsl)$ and hence since $\pi$ is factorial $\pi(R(\lambda,h))\in\C\un\backslash 0$,
\ie $h\in X_T$. \chop
(iv) Given the basis $\big\{q_1,\ldots,\,q_n\big\}$ of $X_T$ then by 
Lemma~\ref{SympFacts}(iii)  in the Appendix 
there are conjugates $\big\{p_1,\ldots,\,p_n\big\}\subset X$
which augments it into a symplectic basis of their span $Q$, and by part~(ii)
above, all these $p_i$ are in $X_S$. 
Obviously $Q$ is nondegenerate.
Then by Lemma~\ref{SympFacts}(ii) we have the decomposition
$X=Q\oplus Q^\perp$ into nondegenerate spaces.
Since $X_T\subset Q$  we have that
$X_T\cap \big(Q^\perp\cap X_R\big)=\{0\}$, \ie
$Q^\perp\cap X_R\subset\{0\}\cup\big(X_R\backslash X_T\big)$.
Now we have the linear decomposition $X_R={(Q^\perp\cap X_R)\mathrel{\dot{\mathord{+}}}X_T},$
\ie any $f\in X_R$ has a unique decomposition $f=f_1+f_2$ with $f_1\in Q^\perp\cap X_R$ 
and $f_2\in X_T\,.$ Specifically, we have 
\[
f_1:=f-\sum_{k=1}^n\sigma(p_k,f)\,q_k\in X_R\qquad\hbox{and}\qquad
f_2:=\sum_{k=1}^n\sigma(p_k,f)\,q_k\in X_T\,.
\]
From (ii) we see that $f_1\in {\{q_1,\ldots,\,q_k\}^\perp}$ and by construction
$f_1\in {\{p_1,\ldots,\,p_k\}^\perp}$ and thus $f_1\in Q^\perp\cap X_R$.
We can now show that
$Q^\perp\cap X_R$ is nondegenerate. If it is not, then there is
a nonzero $h\in Q^\perp\cap X_R$ such that ${\sigma\big(h,\,Q^\perp\cap X_R\big)}=0$.
Then ${\sigma\big(h,\, X_R\big)}=0$ via the decomposition above, using
${\sigma\big(h,\, X_T\big)}=0$ by $h\in Q^\perp\subset X_T^\perp$ and
${\sigma\big(h,\, Q^\perp\cap X_R\big)}=0$.
 Now for $k\in X_S=X\backslash X_R$ we get 
$\pi(R(\lambda,k))=0$ by part~(iii) above, and combining this with the previous fact
gives
\[
\big[\pi(R(\lambda,h)),\,\pi(R(\mu,f))\big]=0, \quad f\in X\,.
\]
Since $\pi$ is factorial this means that $\pi(R(\lambda,h))\in\C\un$ which contradicts 
the fact that $h\in X_R\backslash X_T$. Thus $Q^\perp\cap X_R$ is nondegenerate,
and as $Q$ is nondegenerate as well, we have by
Lemma~\ref{SympFacts}(ii) the decomposition
$X=Q\oplus(Q^\perp\cap X_R)\oplus(Q^\perp\cap X_R^\perp)$ into nondegenerate spaces,
since the decomposition of $X_R$ above implies that
\[
\big(Q\oplus(Q^\perp\cap X_R)\big)^\perp=(Q+X_R)^\perp=Q^\perp\cap X_R^\perp.
\]
Since we have a partition $X=X_R\cup X_S$ and $X_R\subset Q\oplus (Q^\perp\cap X_R)$
we conclude that $(Q^\perp\cap X_R^\perp)\subset X_S\cup\{0\}$.

\subsection*{Proof of Proposition \ref{NonregApprox}}
   
(i) Recall that from $\pi_\omega$ we obtain the decomposition in Eq.~(\ref{Qdecomp})
\[
X=Q\oplus(Q^\perp\cap X_R)\oplus(Q^\perp\cap X_R^\perp)
\]
into nondegenerate subspaces. We first show that there is a regular state
$\wt\omega$ which coincides with $\omega$ on the *--algebra
\[
\hbox{*--alg}\big\{\rlf\,\mid\,f\in Y:=Q^\perp\cap X_R,\;\lambda\in\R\backslash 0\big\}\,.
\]
Since $Y$ is nondegenerate we can use Theorem~\ref{RegThm}(viii) to define a
representation $\pi_1:\overline{\Delta(Y,\sigma)}\to\al B.(\al H._\omega)$ by
$\pi_1(\delta_f):=\exp(i\j_{\pi_\omega}(f))$ for $f\in Y$, and hence a 
regular state $\omega_1(A):={\big(\Omega_\omega,\,\pi_1(A)\Omega_\omega\big)}$
for $A\in \overline{\Delta(Y,\sigma)}$.
{}From the decomposition $X=Y\oplus Y^\perp$ we obtain the 
(minimal) tensor product
$\overline{\Delta(X,\sigma)}=\overline{\Delta(Y,\sigma)}\otimes
\overline{\Delta(Y^\perp,\sigma)}$ \cite{Man}. Define the regular
state $\wt\omega:=\omega_1\otimes\omega_2$ on $\overline{\Delta(X,\sigma)}$,
where $\omega_2$ is any regular state on $\overline{\Delta(Y^\perp,\sigma)}$.
By  Corollary~\ref{RegBij} this corresponds to a regular state $\wt\omega$ on
$\rsl$ (via spectral theory) and it is clear that $\wt\omega$ coincides with
$\omega$ on the *--algebra generated by ${R(\lambda,Y)}$.

Next, we want to construct from $\wt\omega$ a sequence of regular states 
$\omega_n:=\wt\omega\circ\gamma_n$ which will
converge to $\omega$ in the w*--topology, where we now define
the automorphisms $\gamma_n$.

Start with a basis
${\big\{q_1,\ldots,\,q_t\big\}}$ of $X_T$ and augment it by
${\big\{p_1,\ldots,\,p_t\big\}}\subset X_S$ into a symplectic
basis of $Q:={\rm Span}{\big\{q_1,p_1;\ldots;\,q_t,p_t\big\}},$
cf. Proposition~\ref{Xdecomp}(iv). Let
${\big\{q_{t+1},p_{t+1};\ldots;\,q_r,p_r\big\}}$ be a symplectic basis
of $Y:=Q^\perp\cap X_R$,
and let ${\big\{q_{r+1},p_{r+1};\ldots;\,q_s,p_s\big\}}$ be a symplectic basis
of $Q^\perp\cap X_R^\perp\subset X_S$.
Thus we get a symplectic basis of $X$:
\[
\big\{q_1,p_1;\ldots;\,q_t,p_t;\,q_{t+1},p_{t+1};\ldots;\,q_r,p_r;\,
q_{r+1},p_{r+1};\ldots;\,q_s,p_s\big\}
\]
which coincides with the decomposition of $X$ above.
We decompose the elements 
\begin{eqnarray*}
f &=& \sum_{j = 1}^s (x_j q_j + y_j p_j) \in X\quad\hbox{according to}\quad
f = f_T + f_{Q \backslash T} + f_R + f_S\quad\hbox{where:}   \\[1mm]
f_T &:=&\sum_{j=1}^tx_jq_j\in X_T,\qquad  
 f_{Q\backslash T}
:=\sum_{j=1}^ty_jp_j\in Q\backslash X_T\subset X_S, \\[1mm]  
f_R &:=& \sum_{j=t+1}^r \big(x_jq_j+y_jp_j\big)\in Q^\perp\cap X_R\,,\qquad
 f_S:=\sum_{j=r+1}^s \big(x_jq_j+y_jp_j\big)\in Q^\perp\cap X_R^\perp\subset X_S\,.
\end{eqnarray*}
We define now for $n\in\N$ the symplectic transformation $T_n^{(i)}$ of 
$f =  \sum_{j = 1}^s (x_j q_j + y_j p_j) \in X$ by
\[ T_n^{(i)}(f) ={x_i\over n}\,q_i +n\,y_i\,p_i 
+\sum_{j=1\atop j\not=i}^s\big(x_jq_j+y_jp_j\big) 
\]
and the automorphisms $\alpha_n^{(i)}\in\aut\rsl$ by 
$\alpha_n^{(i)}\big(\rlf\big):=R\big(\lambda,\,T_n^{(i)}(f)\big)$ for $f\in X$,
cf. Theorem~\ref{Relemen}(v). Next, using Proposition~\ref{ImAut} we define
the automorphisms $\beta_n\in\aut\rsl$ by
$\beta_n\big(\rlf\big):=R\big(\lambda+i\sigma(h_n + k,f),\,f\big)$ for all
$f\in X$, $\lambda\in\R\backslash 0$ where
\[
h_n:=n\sum_{j=r+1}^s\big(n^{j-r}q_j+n^{j+s-2r}p_j\big)\quad
\hbox{and}\quad k:= \sum_{j=1}^tb_jp_j
\]
with $b_j :=(\Omega_\omega, \j_{\pi_\omega}(q_j) \Omega_\omega ) $; note that 
the fields 
$\j_{\pi_\omega}(q_j)$, $j = 1, \dots t$ are multiples of the identity according to 
Proposition \ref{Xdecomp}(iii).
We will make use of the fact that 
\[
\sigma(h_n+ k,\,f)=n\sum_{j=r+1}^s(x_jn^{j+s-2r}-y_jn^{j-r})
+ \sum\limits_{j=1}^tb_jx_j\;
\maprightt \infty,n.\pm\infty
\]
if any of the coefficients $x_j,\,y_j$ are nonzero for $r+1\leq j\leq
s$, \ie if $f_S \neq 0$. For in that 
case $\sigma(h_n+ k,\,f)$ is a polynomial in $n$ with degree
$\geq 2$.  

Now define $\gamma_n:=\alpha_n^{(1)}\cdots\alpha_n^{(t)}\beta_n\in\aut\rsl$, so
we have
\begin{eqnarray*}
\gamma_n\big(\rlf\big)&=& R\Big(\lambda+i\sigma(h_n+k,f),\;\sum_{j=1}^t\big(
{x_j\over n}\,q_j+n\,y_j\,p_j\big) + \sum_{j=t+1}^s\big(x_jq_j+y_jp_j\big)
\Big) \\[1mm]
&=& R\Big(\lambda+i\mu^f+i\lambda_n^f,\; \f 1,n.f_T+n\,f_{Q\backslash T}
+f_R+f_S\Big) \\[1mm]
\hbox{where}\qquad\mu^f&:=&\sigma(k,\,f)=
\big(\Omega_\omega, \j_{\pi_\omega} (f_T) \Omega_\omega \big)\,,
\qquad\lambda_n^f:=\sigma(h_n,\,f)\maprightt\infty,n.\cases{\pm\infty
& if $f_S\not=0\,;$ \cr 0 & if $f_S=0$\,. } 
\end{eqnarray*}
Since $\gamma_n$ maps resolvents to resolvents, $\pi_{\tilde\omega}\circ\gamma_n$
is also regular, where ${\wt\omega}$ is the regular state obtained
above. 
Next, we wish to find the limits $\slim\limits_{n\to\infty}\pi_{\tilde\omega}\circ\gamma_n
\big(\rlf\big)$. \chop
$\bullet$ The case $f_{Q\backslash T}\not=0$. 
We begin by noting that for any real polynomial $\xi_n$ in $n$ of
degree $\geq 2$ or for $\xi_n = 0$, $n \in \N$ one has
\begin{equation}
\label{RnQTto0} 
\pi_{\tilde\omega}\big(R(\lambda +i \xi_n,\, n\,f_{Q\backslash T})\big)
=\int{dP(\mu)\over i\lambda - \xi_n - n \mu}\;\maprightt\infty,n.\;0
\end{equation}
in strong operator topology, where $dP$ is the spectral measure 
of $\j_{\pi_{\tilde\omega}}(f_{Q\backslash T})$. Now from 
Eq.~(\ref{Rsum}) we get that
\begin{eqnarray*}
&&\pi_{\tilde\omega}\Big(R\big(\f\lambda,2.+i\mu^f,\,\f 1,n.f_T+
f_R+f_S\big)\,R\big(\f\lambda,2.+i\lambda_n^f,\,
n\,f_{Q\backslash T}\big)\Big)  \\[1mm]
&&=\pi_{\tilde\omega}\Big(R\big(\lambda+i\mu^f+i\lambda_n^f,\; \f 1,n.f_T+n\,f_{Q\backslash T}
+f_R+f_S\big)\,\Big[R\big(\f\lambda,2.+i\mu^f,\,\f 1,n.f_T+
f_R+f_S\big)+ \\[1mm]
&&\quad+R\big(\f\lambda,2.+i\lambda_n^f,\,
n\,f_{Q\backslash T}\big)+i\sigma(f_T,f_{Q\backslash T})\,
R\big(\f\lambda,2.+i\mu^f,\,\f 1,n.f_T+
f_R+f_S\big)^2\,R\big(\f\lambda,2.+i\lambda_n^f,\,
n\,f_{Q\backslash T}\big)\Big]\Big) \, .
\end{eqnarray*}
If we let $n\to\infty$ in the strong operator topology, and use the fact that
multiplication is jointly continuous on bounded sets in the strong operator topology,
then by Eq.~(\ref{RnQTto0}) for $\xi_n=\lambda^f_n$,
the left hand side of this is zero, and we get:
\[
0=\slim_{n\to\infty}\pi_{\tilde\omega}\Big(R\big(\lambda+i\mu^f+i\lambda_n^f,\; \f 1,n.f_T+n\,f_{Q\backslash T}
+f_R+f_S\Big)\,R\big(\f\lambda,2.+i\mu^f,\,\f 1,n.f_T+
f_R+f_S\big)\Big)\,.
\]
Consider the last factor. If $f_R+f_S=0$ then by spectral theory 
of $\j_{\pi_{\tilde\omega}}(f_T)$ we get that 
\[
\slim\limits_{n\to\infty} \, 
\pi_{\tilde\omega}\big(R\big(\f\lambda,2.+i\mu^f,\,\f 1,n.f_T\big)\big)
={R\big(\f\lambda,2.+i\mu^f,\,0\big)} = {1\over i \lambda/2 - \mu^f}\, {\un} \,,
\]
and hence conclude that
\[
0=\slim_{n\to\infty}\pi_{\tilde\omega}\Big(R\big(\lambda+i\mu^f+i\lambda_n^f,\; \f 1,n.f_T+n\,f_{Q\backslash T}
+f_R+f_S\big)\Big)=\slim\limits_{n\to\infty}\pi_{\tilde\omega}\circ\gamma_n
\big(\rlf\big)\,.
\]
If $f_R+f_S\not=0$ then since the resolvent  of 
$\j_{\pi_{\tilde\omega}}(f_T)$ commutes with that of $\j_{\pi_{\tilde\omega}}(f_R+f_S)$
we can use their joint spectral theory:
\begin{eqnarray*}
&& \pi_{\tilde\omega}\Big(R\big(\f\lambda,2.+i\mu^f,\,\f 1,n.f_T+
f_R+f_S\big)\Big)=\int{dP(\rho, \sigma)\over 
i\f\lambda,2.-\mu^f-\f 1,n.\, \rho - \sigma} \\[1mm]
&&\maprightt\infty,n.\;\int{dP(\rho,\sigma)\over i\f\lambda,2.-\mu^f- \sigma}=
\pi_{\tilde\omega}\big(R\big(\f\lambda,2.+i\mu^f,\,f_R+f_S\big)\big)
\end{eqnarray*}
where the limit is in the strong operator topology.
Thus we obtain
\[
0=\slim_{n\to\infty}\pi_{\tilde\omega}\Big(R\big(\lambda+i\mu^f+i\lambda_n^f,\; \f 1,n.f_T+n\,f_{Q\backslash T}
+f_R+f_S\Big)\Big)\,\pi_{\tilde\omega}\big(R\big(\f\lambda,2.+i\mu^f,\,f_R+f_S\big)\big)\,.
\]
Since $\pi_{\tilde\omega}$ is regular, the last factor has dense range, hence using 
uniform boundedness of the resolvents in the first factor, we find again that
\[
0=\slim_{n\to\infty}\pi_{\tilde\omega}\Big(R\big(\lambda+i\mu^f+i\lambda_n^f,\; \f 1,n.f_T+n\,f_{Q\backslash T}
+f_R+f_S\big)\Big)=\slim_{n\to\infty}\pi_{\tilde\omega}\circ\gamma_n
\big(\rlf\big)\,.
\]
$\bullet$ The case $f_{Q\backslash T}=0$. We have that if $f_T\not=0$ then
\begin{equation}
\label{RnTto0}
\pi_{\tilde\omega}\big(R(\lambda+i\mu^f+i\lambda_n^f,\, \f 1,n.f_T)\big)
=\int{dP(\rho)\over i\lambda-\mu^f-\lambda_n^f-\f 1,n. \rho}
\;\maprightt\infty,n.\;
\cases{\hspace*{2mm} 0 & if $f_S\not=0$ \cr \f 1,
i\lambda-\mu^f.\,\un &if $f_S=0$.\cr} 
\end{equation}
in strong operator topology, where $dP$ is the spectral measure 
of $\j_{\pi_{\tilde\omega}}(f_T)$. When $f_T=0$ then $\mu^f=0$ and  we get:
\[
\pi_{\tilde\omega}\big(R(\lambda+i\mu^f+i\lambda_n^f,\, \f 1,n.f_T)\big)
={1\over i\lambda-\lambda_n^f}\,\un\;\maprightt\infty,n.\;
\cases{\hspace*{2mm} 0 & if $f_S\not=0$ 
\cr \f 1,i\lambda.\,\un &if $f_S=0$\,.\cr}
\]
By Eq.~(\ref{Rsum}) we get
\begin{eqnarray*}
&&\pi_{\tilde\omega}\Big(R\big(\f\lambda,2.+i\mu^f+i\lambda_n^f,\,\f 1,n.f_T\big)
\,R\big(\f\lambda,2.,\,f_R+f_S\big)\Big)  \\[1mm]
&&=\pi_{\tilde\omega}\Big(R\big(\lambda+i\mu^f+i\lambda_n^f,\; \f 1,n.f_T
+f_R+f_S\big)\,\Big(R\big(\f\lambda,2.+i\mu^f+i\lambda_n^f,\,\f 1,n.f_T\big)+
R\big(\f\lambda,2.,\,f_R+f_S\big)\Big)\Big) \, .
\end{eqnarray*}
Now if we let $n\to\infty$ in the strong operator topology, and use Eq.~(\ref{RnTto0})
we find that if $f_S\not=0$ then
\[
0=\slim_{n\to\infty}\pi_{\tilde\omega}\Big(R\big(\lambda+i\mu^f+i\lambda_n^f,\; \f 1,n.f_T
+f_R+f_S\big)\Big)\,\pi_{\tilde\omega}\Big(R\big(\f\lambda,2.,\,f_R+f_S\big)\Big)
\]
hence since $\pi_{\tilde\omega}$ is regular,
\[
0=\slim_{n\to\infty}\pi_{\tilde\omega}\Big(R\big(\lambda+i\mu^f+i\lambda_n^f,\; \f 1,n.f_T
+f_R+f_S\big)\Big)=\slim_{n\to\infty}\pi_{\tilde\omega}\circ\gamma_n
\big(\rlf\big)\,.
\]
If $f_S=0$ then
\[
{ 1\over i\f\lambda,2.-\mu^f}\,\pi_{\tilde\omega}\Big(R\big(\f\lambda,2.,\,f_R\big)\Big)
=\slim_{n\to\infty}\pi_{\tilde\omega}\Big(R\big(\lambda+i\mu^f,\; \f 1,n.f_T
+f_R\big)\,\Big({ 1\over i\f\lambda,2.-\mu^f}\,\un+
R\big(\f\lambda,2.,\,f_R\big)\Big)\Big)\,.
\]
Since the second factor on the right hand side is invertible (as it is a normal
operator with continuous spectrum), we obtain via spectral theory
of $\j_{\pi_{\tilde\omega}}(f_R)$ that the equation rearranges to
\[
\pi_{\tilde\omega}\Big(R\big(\lambda+i\mu^f,\,f_R\big)\Big)
=\slim_{n\to\infty}\pi_{\tilde\omega}\Big(R\big(\lambda+i\mu^f,\; \f 1,n.f_T
+f_R\big)=\slim_{n\to\infty}\pi_{\tilde\omega}\circ\gamma_n\big(\rlf\big)
\]
when $f_R\not=0$, and the same equation follows trivially when $f_R=0$.
Summarizing the reasoning above, we conclude for a general $f$ that
\begin{eqnarray*}
&&\slim_{n\to\infty}\,\pi_{\tilde\omega}\circ\gamma_n\big(\rlf\big)=
\cases{0 & if $f_{Q\backslash T}+f_S\not=0$ \cr
\pi_{\tilde\omega}\Big(R\big(\lambda+i\mu^f,\,f_R\big)\Big)
& if $f_{Q\backslash T}+f_S=0$ . \cr}  
\end{eqnarray*}
Hence, for $m \in \N$, 
\begin{eqnarray*}
&& \slim_{n\to \infty}\,\pi_{\tilde\omega}\circ\gamma_n\big(
R(\lambda_{(1)},\,f_{(1)})\cdots R(\lambda_{(m)},\,f_{(m)})\big) \\[1mm]
&&=\cases{0 & if $f_{(i)\,Q\backslash T}+f_{(i)\,S}\not=0$ for any $i$, \cr
\pi_{\tilde\omega}\big(R(\lambda_{(1)}+i\mu^{f_{(1)}},\,f_{(1)\,R})
\cdots R(\lambda_{(m)}+i\mu^{f_{(m)}},\,f_{(m)\,R})\big)
& if $f_{(i)\,Q\backslash T}+f_{(i)\,S}=0$ for all $i$,\cr}  
\end{eqnarray*}
and so
\begin{eqnarray*}
&&\lim_{n\to\infty}\,{\wt\omega}\circ\gamma_n\big(
R(\lambda_{(1)},\,f_{(1)})\cdots R(\lambda_{(m)},\,f_{(m)})\big) \\[1mm]
&&=\cases{0 & if $f_{(i)\,Q\backslash T}+f_{(i)\,S}\not=0$ for any $i$, \cr
{\wt\omega}\Big(R\big(\lambda_{(1)}+i\mu^{f_{(1)}},\,f_{(1)\,R}\big)
\cdots R\big(\lambda_{(m)}+i\mu^{f_{(m)}},\,f_{(m)\,R}\big)\Big)
& if $f_{(i)\,Q\backslash T}+f_{(i)\,S}=0$ for all $i$,\cr}   \\[1mm]
&&=\cases{0 & if $f_{(i)\,Q\backslash T}+f_{(i)\,S}\not=0$ for any $i$, \cr
{\omega}\Big(R\big(\lambda_{(1)}+i\mu^{f_{(1)}},\,f_{(1)\,R}\big)
\cdots R\big(\lambda_{(m)}+i\mu^{f_{(m)}},\,f_{(m)\,R}\big)\Big)
& if $f_{(i)\,Q\backslash T}+f_{(i)\,S}=0$ for all $i$\cr}  
\end{eqnarray*}
since $\omega$ coincides with ${\wt\omega}$ on $\hbox{*-alg}\big\{R(z,f)\,\mid\,
f\in Q^\perp\cap X_R,\,z\in\C\backslash i\R\big\}$.
By Proposition~\ref{Xdecomp}(iii) we know that $\pi_\omega(\rlf)=0$
if $f_{Q\backslash T}+f_S\not=0$ and if $f_{Q\backslash T}+f_S=0$ then
\[
\pi_\omega(\rlf)=\big[ i\lambda\un-\j_{\pi_\omega}(f_T+f_R)\big]^{-1}
=\big[(i\lambda-\mu^f)\un-\j_{\pi_\omega}(f_R)\big]^{-1}=\pi_\omega\big(
R(\lambda+i\mu^f,\,f_R)\big) 
\]
by $\j_{\pi_\omega}(f_T) = \mu_f \un$.  Thus 
\[
\lim_{n\to\infty}\,{\wt\omega}\circ\gamma_n\big(
R(\lambda_{(1)},\,f_{(1)})\cdots R(\lambda_{(m)},\,f_{(m)})\big)
=\omega\big(
R(\lambda_{(1)},\,f_{(1)})\cdots R(\lambda_{(m)},\,f_{(m)})\big)
\]
\ie $\omega(A)=\lim\limits_{n\to\infty}\,{\wt\omega}\circ\gamma_n(A)$
for all  $A\in\al R._0$. Since $\al R._0$ is norm dense in $\rsl$ 
the sequence of regular states
${\wt\omega}\circ\gamma_n$ converges to the nonregular state $\omega$
in the weak *-topology.\chop
(ii) The vector states of the universal regular representation $\pi_r$ include
the regular states ${\wt\omega}\circ\gamma_n$ constructed above. Thus from part~(i)
we see that the given nonregular pure state $\omega$ is in the w*-closure
of the convex hull of the vector states of $\pi_r$ and hence by Fell's theorem,
(cf.~Theorem~1.2 in \cite{Fell} and \cite[p 106]{Emch}) 
 we find that $\ker\pi_r\subseteq\ker\pi_\omega$.\chop
(iii) {}From part (ii) we see that if $\pi_r(A)=0$ then $\pi_\omega(A)=0$ for
all pure states $\omega$ of $\rsl$. However according to standard
results in the theory of C*--algebras, cf.\
Lemma 2.3.23 in \cite{BR2}, the norm of $\rsl$
is $\|A\|={\sup\big\{\omega(A^*A)^{1/2}\,\mid\,\omega\in\ot S.\;\hbox{pure}\big\}}$
and hence $A=0$. Thus $\pi_r$ is faithful and so 
$\rsl\cong\al R._r(X,\sigma)$. 

\subsection*{Proof of Theorem \ref{UniqueR}}

(i)
On $\hbox{*-alg}\big\{\rlf\,\mid\,f\in S,\;\lambda\in\R\backslash 0\big\}
\subset\rsl$
we have the following four norms: ${\|\cdot\|_X}=$ C*-norm of $\rsl$,
${\|\cdot\|_S}=$ C*-norm of $\al R.(S,\sigma)$, and
${\|\cdot\|_{{\rm reg}S}}=$ regular norm of $\al R.(S,\sigma)$, and 
${\|\cdot\|_{{\rm reg}X}}=$ regular norm of $\rsl$.
We first prove that ${\|A\|_{{\rm reg}S}}={\|A\|_{{\rm reg}X}}$
for all $A$ in this *-algebra. Since $S$ is finite dimensional
and nondegenerate, we have
$X=S\oplus S^\perp$
by Lemma~\ref{SympFacts}(ii).
By definition of the regular seminorm on $\rsl$ we have
\[
\|A\|_{{\rm reg}X}=\|\pi_r(A)\|=
\sup\big\{\|\pi_\omega(A)\|\,\mid \omega\in\ot S._r\big(\rsl\big)\big\}
=\sup\big\{\sqrt{\omega(A^*A)} \,\mid\,\omega\in\ot S._r\big(\rsl\big)\big\}\,.
\]
Thus this will coincide with the regular seminorm of $\al R.(S,\sigma)$ on
$\hbox{*-alg}\big\{\rlf\,\mid\,f\in S,\;\lambda\in\R\backslash 0\big\}$
if we can show that each $\omega\in{\ot S._r\big(\al R.(S,\,\sigma)\big)}$ extends to
a regular state of $\rsl.$ By the bijection of Corollary~\ref{RegBij} this will be the case if 
each regular state of  $\ccr S,\sigma.\subset\CCRX$ extends to a regular
state of $\CCRX.$ Now we know by Manuceau~\cite{Man} that
$\CCRX=\ccr S,\sigma.\otimes\ccr S^\perp,\sigma.,$ so if $\omega_1$ is a regular state
of $\ccr S,\sigma.$ then $\omega_1\otimes\omega_2$ will be a regular extension
of $\omega_1$ to $\CCRX$ if we choose $\omega_2$ to be a regular state of
$\ccr S^\perp,\sigma.,$ which is of course possible. Thus it follows that
${\|A\|_{{\rm reg}S}}={\|A\|_{{\rm reg}X}}$.

All (regular) states of $\rsl$ restrict to  
(regular) states of ${\hbox{*-alg}\big\{\rlf  
\mid f\in S,\;\lambda\in\R\backslash 0\big\},}$
so it follows that for all $A$ in this *-algebra: 
\[
\|A\|_{S}\geq\|A\|_{X}\geq\|A\|_{{\rm reg}X}=\|A\|_{{\rm reg}S} \, .
\]
However, since  $\|A\|_{S}=\|A\|_{{\rm reg}S}$ by 
Proposition~\ref{NonregApprox}(iii), it follows that 
$\|A\|_{S}=\|A\|_{X}$, which establishes the claim.
Hence the containment
${\hbox{*-alg}\big\{\rlf\,\mid\,f\in S,\;\lambda\in\R\backslash 0\big\}}
\subseteq\rsl$ extends to an isomorphism of $\al R.(S,\sigma)$ with a subalgebra
of $\rsl$, and we indicate this as containment.\chop
(ii) 
Since the net $S\to\al R.(S,\sigma)$ has the partial ordering of containment
of the finite dimensional nondegenerate spaces $S \subset X$ (by (i)), and 
$\rsl$ is generated by all $\al R.(S,\sigma)$,
it is clear that it is the inductive limit of the net.\chop
(iii)
{}From the inductive limit in part (ii), it suffices to verify that $\pi_r$
restricts to an isomorphism on each $\al R.(S,\sigma)$, but this holds by
$\|A\|_{S}=\|A\|_{{\rm reg}X}$ mentioned above.

\subsection*{Proof of Theorem \ref{RegFaith}}

Let $\pi_i:\rsl\to\al B.(\al H._i)$, $i=1,\,2$ be regular representations.
Since $\CCRX$ is simple, we have a *-isomorphism
$\alpha:\pi_1\big(\CCRX\big)\to\pi_2\big(\CCRX\big)$ by 
$\alpha\big(\pi_1(A)\big):=\pi_2(A)$ for all $A\in\CCRX$.
We want to extend $\alpha$ to 
$\pi_1\big(\rsl\big) \subset \pi_1\big(\CCRX\big)''$. 
Let $S\subset X$ be a finite dimensional nondegenerate subspace, then 
by the von Neumann uniqueness theorem, both of $\pi_i\restriction\ccr S,\sigma.$ are normal
to the Fock representation of $\ccr S,\sigma.$, hence 
$\pi_1\restriction\ccr S,\sigma.$ is normal to $\pi_2\restriction\ccr S,\sigma.$.
Then by Theorem~2.4.26 in \cite[p 80]{BR1} we conclude that 
$\alpha:\pi_1\big(\ccr S,\sigma.\big)\to\pi_2\big(\ccr S,\sigma.\big)$
is normal and extends to a *-homomorphism
$\alpha:\pi_1\big(\ccr S,\sigma.\big)''\to\pi_2\big(\ccr S,\sigma.\big)''$
by strong operator continuity. Now $\pi_i\big(\ccr S,\sigma.\big)''
\supset\pi_i\big(\al R.(S,\sigma)\big)$, and in fact by the Laplace 
transform~(\ref{Laplace1}), for each $A\in\al R._0(S)$ there is a sequence
$\{A_n\}\subset\ccr S,\sigma.$ such that $\pi_i(A_n)\maprightt \infty,n.\pi_i(A)$
in the strong operator topology. This means that
$\alpha\left(\pi_1\big(\al R._0(S)\big)\right)\subseteq\pi_2\big(\al R._0(S)\big)$
and thus $\alpha$ restricts to a *-homomorphism
$\alpha:\pi_1\big(\al R.(S,\sigma)\big)\to \pi_2\big(\al R.(S,\sigma)\big)$
and so 
\[
\big\|\pi_2(A)\big\|=\big\|\alpha(\pi_1(A))\big\|\leq\big\|\pi_1(A)\big\|
\, ,
\quad A\in\al R.(S,\sigma)\,.
\]
By symmetry of the argument we also get that
$\big\|\pi_1(A)\big\|\leq\big\|\pi_2(A)\big\|$ and hence that
$\big\|\pi_1(A)\big\|=\big\|\pi_2(A)\big\|$ for all $A\in\al R.(S,\sigma)$.
Let the regular representation $\pi_2=\pi_r$ which is faithful by
Theorem~\ref{UniqueR}(iii), then we have obtained that 
$\big\|\pi_1(A)\big\|=\|A\|$ for all $A\in\al R.(S,\sigma)$
and for all finite dimensional nondegenerate subspaces $S\subset X$.
Since by Theorem~\ref{UniqueR}(ii) we know that $\rsl$ is the inductive limit
of all the $\al R.(S,\sigma)$, it follows that $\pi_1$ is faithful
on all of $\rsl$.

\subsection*{Proof of Proposition \ref{Uaut}}

Note first that if $\alpha$ corresponds to a symplectic transformation,
then so does its inverse. Moreover, $\alpha$ and $\alpha^{-1}$
preserve both the set of regular states
${\ot S._r\big(\rsl\big)}$ and the strongly regular states
${\ot S._{sr}\big(\rsl\big)}$, respectively. Let 
${\ot S.}$ be either one of these sets of states 
and put $\pi_{\ot S.}:=\bigoplus\limits_{\omega\in\ot S.}\pi_\omega$. 
Since  both $\alpha$ and $\alpha^{-1}$ preserve $\ot S.$, 
one obtains a bijection of $\ot S.$ 
by $\omega\mapsto\omega\circ\alpha$.
Hence $\pi_{\ot S.}\circ\alpha$
is just $\pi_{\ot S.}$ where its direct summands have been permuted.
 Such a permutation of direct summands can be done by
 conjugation of a unitary, thus
 $\pi_{\ot S.}$ is unitarily equivalent to $\pi_{\ot S.}\circ\alpha$.

\subsection*{Proof of Theorem \ref{TensorAlg}}

By Theorem~\ref{UniqueR}(i) and (ii) the C*--algebras generated in $\al R.(X,\,\sigma)$ by  
${\{\rlf\,\mid\,f\in S,\,\lambda\in\R\backslash 0\}}$ and
${\{\rlf\,\mid\,f\in S^\perp,\,\lambda\in\R\backslash 0\}}$
are $\al R.(S,\sigma)$ and $\al R.(S^\perp,\sigma)$.
We already have that $\CCRX=\ccr S,\sigma.\otimes\ccr S^\perp,\sigma.$
by Manuceau~\cite{Man}.
Consider a representation $\pi=\pi_1\otimes\pi_2$ of $\CCRX$ where 
$\pi_1$ (resp. $\pi_2$) is a regular representation of 
$\ccr S,\sigma.$ (resp. $\ccr S^\perp,\sigma.$). Then $\pi$ is regular,
hence extends to a representation of $\al R.(X,\sigma)$ by
$\pi\big(\al R.(X,\sigma)\big)\subset\pi\big(\CCRX\big)''$ as discussed before,
and likewise for $\pi_1$ and $\pi_2$.
Moreover, spectral theory respects tensor products, so
 if $A\in\al R.(S,\sigma)$ and $B\in\al R.(S^\perp,\sigma)$ then
 $\pi(A)=\pi_1(A)\otimes\un$ and $\pi(B)=\un\otimes\pi_2(B)$ hence
$\pi(AB)=\pi(A)\pi(B)=\pi_1(A)\otimes\pi_2(B)$.
By Theorem~\ref{UniqueR}(iii) we can choose $\pi_1$ and $\pi_2$ to
be faithful, hence 
\[
\|AB\|\geq\|\pi(AB)\|=\|\pi_1(A)\|\cdot\|\pi_2(B)\|=\|A\|\cdot\|B\|
\]
 so we conclude that $AB=0$ implies that at least one of $A$ and $B$ must be zero.
 Since $\al R.(S,\sigma)$ and $\al R.(S^\perp,\sigma)$ are commuting subalgebras
 of $\al R.(X,\sigma)$ we conclude from this via
 Exercise~2 in Takesaki~\cite[p 220]{Tak}, that
 \[
 C^*\big(\al R.(S,\sigma)\cup \al R.(S^\perp,\sigma)\big)
 \cong\al R.(S,\sigma)\otimes \al R.(S^\perp,\sigma)\,.
 \]
To see that the containment
$\rsl\supset C^*\big(\al R.(S,\sigma)\cup \al R.(S^\perp,\sigma)\big)$
is in general proper, we present a simple example.
Let ${\rm dim}(X)=4$, and choose a symplectic basis 
${\big\{q_1,p_1;\,q_2,p_2\big\}}$ and let 
$S:={\rm Span}\{q_1,p_1\big\}$, hence $S^\perp={\rm Span}\{q_2,p_2\big\}$
and $X=S\oplus S^\perp$. Choose a fixed regular state $\omega$ on $\rsl$
and define the automorphisms $\beta_n\in\aut\rsl$ by
$\beta_n\big(\rlf\big):={R\big(\lambda+i\sigma(h_n,f),\,f\big)}$
where $h_n:=n(p_1-p_2)-n^2(q_1-q_2)$, making use of Proposition~\ref{ImAut}.
So if $f=x_1q_1+x_2q_2+y_1p_1+y_2p_2$ then 
$\sigma(h_n,f)=n(x_1-x_2)+n^2(y_1-y_2)$. Thus
\begin{eqnarray*}
\slim_{n\to\infty}\,\pi_\omega\circ\beta_n\big(\rlf\big)
&=&\slim_{n\to\infty}\,\pi_\omega\big(R\big(\lambda+i\sigma(h_n,f),\,f\big)\big) \\[1mm]
&=& \slim_{n\to\infty}\,\int{dP(t)\over i\lambda-\sigma(h_n,f)-t}
=\cases{\pi_\omega\big(\rlf\big) & if $x_1=x_2$ and $y_1=y_2$\cr
0 & if $x_1\not=x_2$ or  $y_1\not=y_2$ \cr}
\end{eqnarray*}
where $dP$ is the spectral measure of $\j_{\pi_\omega}(f)$.
Now proceeding as at the end of the proof of Proposition~\ref{NonregApprox}(i),
we conclude that the w*-limit $\wt\omega:=\lim\limits_{n\to\infty}\omega\circ\beta_n$ defines a state on
$\rsl$ such that
\[
\wt\omega(\rlf)=\cases{\omega\big(\rlf\big) & if $x_1=x_2$  and $y_1=y_2$ \cr
0 & if $x_1\not=x_2$ or  $y_1\not=y_2$\,.\cr} 
\]
Thus $\wt\omega(\rlf)=0$ if $f\in S\backslash 0$ or  $f\in S^\perp\backslash 0$, and
by Theorem~\ref{Ideals1}(iv) we have $\rlf\in\ker\pi_{\tilde\omega}$ 
for such $f$.
Since $\ker\pi_{\tilde\omega}$
is a closed two--sided ideal, we get that
$C^*\big(\al R.(S,\sigma)\cup \al R.(S^\perp,\sigma)\big)\subset\ker\pi_{\tilde\omega}$.
However $\wt\omega\big(R\big(\lambda,\,q_1+q_2+p_1+p_2)\big)=\omega
\big(R\big(\lambda,\,q_1+q_2+p_1+p_2)\big)\not=0$,
hence ${R\big(\lambda,\,q_1+q_2+p_1+p_2)}$ is not an element of 
$C^*\big(\al R.(S,\sigma)\cup 
\al R.(S^\perp,\sigma)\big)$ and the containment is proper.

\subsection*{Proof of Theorem \ref{Nonsep}}

(i) 
First, consider the case when $\sigma(f,h)\not=0$.
{}From Proposition~\ref{RDirac}(i) applied to $C=\{f\}$ we see that there is 
a state $\omega$ such that ${\omega\big(R(1,f)\big)}=-i$ (note that the proof
of  Proposition~\ref{RDirac} is logically independent from this Theorem).
Moreover, by Proposition~\ref{RDirac}(ii) we then have that 
${\omega\big(R(1,h)\big)}=0$ and thus by 
Theorem~\ref{Ideals1}(iv) we have
 $\pi_\omega\big(\big[\rsl R(1,h)\big]\big)=0$ but $\pi_\omega\big(R(1,f)\big)\not=0$,
 and hence $R(1,f)\not\in \big[\rsl R(1,h)\big]$.\chop
Next, we prove the claim for the case $\sigma(f,h)=0$.
Augment $f,\,h$ to a symplectic basis $\{f,p_f;\,h,p_h\}$
using Lemma~\ref{SympFacts}(iii)\,.
Let $S:={\rm Span}\{f,p_f; h,p_h\}\subset X$ then by Lemma~\ref{SympFacts}(ii)
 we have that $X=S\oplus S^\perp$ and likewise we get that
 $S=S_1\oplus S_2$ where $S_1:={\rm Span}\{f,p_f\}$ and
 $S_2:={\rm Span}\{h,p_h\}$. Then
 \[
 \al R.(X,\sigma)\supset\al R.(S_1,\sigma)\otimes\al R.(S_2,\sigma)\otimes\al R.(S^\perp,\sigma)
 \]
 by Theorem \ref{TensorAlg}\,.
 But then we can choose a product  state $\omega=\omega_1\otimes\omega_2\otimes\omega_3$
 of  $\al R.(S_1,\sigma)\otimes\al R.(S_2,\sigma)\otimes\al R.(S^\perp,\sigma)$,
 such that
 $\omega_1$ is a Fock state of $\al R.(S_1,\sigma)$, $\omega_2(R(1,h))=0$ 
 (which is possible by
 Theorem~\ref{Ideals1}(iv)), and $\omega_3$ is regular on $\al R.(S^\perp,\sigma)$.
 Extend $\omega$ by the Hahn--Banach theorem to a state on $\rsl$ (still denoted $\omega$).
 Then $\omega\big(R(1,h)\big)=0$ and $\omega\big(R(1,f)\big)\not=0$
 so the statement now follows as in the preceding step. \chop 
(ii)
Consider first the case $\sigma(f,h)\not=0$. Choosing a state $\omega$
as in the preceding step  we obtain
\[
\big\|R(1,f)-R(1,h)\big\|\geq\big|\omega\big(R(1,f)-R(1,h)\big)\big|
=|-i|=1 \;.
\]
{}For the case $\sigma(f,h)=0$ we augment $(f,h)$ to a symplectic basis
$(f,p_f; h,p_h)$ and put  $S:={\rm Span}\{f,p_f; h,p_h\}\subset X$.
By Theorem~\ref{UniqueR}(i)
we have the containment ${\al R.(S,\,\sigma)}\subset\al R.(X,\,\sigma)$. 
Let $\pi$ be the Schr\"odinger representation of ${\al R.(S,\,\sigma)}$. 
Since $\pi$ is regular, by Theorem~\ref{RegFaith}
it is faithful, hence applying the joint spectral theory to the two
resolvents $\pi(R(1,f))$, $\pi(R(1,h))$ we obtain
\[
\big\|R(1,f)-R(1,h)\big\|=\big\|\pi\big(R(1,f)-R(1,h)\big)\big\|=
\sup_{\rho,\,\sigma\in\R}\,\Big|{1\over i-\rho}-{1\over i-\sigma}\Big| = 1 \, .
\]
(iii)
Finally, note that given any $f,\,h$ as above, we can define
$f_\xi:= \xi f+(1- \xi )h$, $\xi \in[0,1]$ to obtain an uncountable family such that
if $\xi \not= \zeta$ then $f_\xi \not\in\R f_\zeta$.
Since then the $R(1,f_\xi)$ are all far apart by 
$\big\|R(1,f_\xi)-R(1,f_\zeta)\big\|\geq 1$ 
it follows that $\rsl$ is nonseparable.

\subsection*{Proof of Theorem \ref{CompactId}}

(i)  Since all irreducible regular representations are unitarily equivalent,
we may assume that $\pi_0$ is the Schr\"odinger representation on 
$L^2(\R^n)$ w.r.t. 
the given symplectic
basis. Taking into account the commutation relations (\ref{Rccr})
of the resolvents, we obtain
\[
\pi_0\Big(\big(R(\lambda_1,p_1)R(\mu_1,q_1)\big)\cdots
\big(R(\lambda_n,p_n)R(\mu_n,q_n)\big)\Big) =
\prod_{j=1}^n \, (i \lambda_j - Q_j )^{-1} \cdot
\prod_{k=1}^n \, (i \mu_k - P_k )^{-1} \, ,
\]
where $Q_j := \phi_{\pi_0}(p_j)$, 
$P_k := \phi_{\pi_0}(q_k)$ are the familiar position
and momentum operators on $L^2(\R^n)$. 
Now for any pair $A, B$
of continuous, bounded and square integrable functions on $\R^n$ 
the operator $A(Q_1, \dots Q_n) \cdot B(P_1, \dots P_n)$ is 
in the Hilbert--Schmidt class cf.~\cite{RS3} Theorem~XI.20.
Thus the above product of
resolvents is in the Hilbert--Schmidt class and hence 
a compact operator. 

\noindent (ii) It is well--known that if a C*--algebra acts irreducibly 
on a Hilbert space, and
contains any nonzero element of the compact operators, then it contains 
all compact operators (cf. \cite{Dix} Theorem~4.1.10 or \cite{Mur} Theorem~2.4.9).
Thus by (i), $\pi_0\big(\rsl\big)$ contains all of $\al K.(\al H._0),$ and so, since 
$\pi_0$ is faithful, $\rsl$ contains an ideal $\al K.$ isomorphic to
$\al K.(\al H._0)$. Uniqueness follows from the fact that up to unitary
equivalence (which preserves the compacts), $\pi_0$ is unique.\chop
(iii)
Since $\al K.$ is a proper closed two--sided ideal of $\rsl,$ each 
representation $\pi:\rsl\to\al B.(\al H._\pi)$ has a unique decomposition
$\pi=\pi_1\oplus\pi_2$ where $\pi_1$ is nondegenerate when restricted to
$\al K.$ and $\pi_2(\al K.)=0$. Now if $\pi$ is regular,
then from the products of resolvents in (i), we obtain a sequence $\{I_n\}$ in
$\al K.$ such that $\pi(I_n)\to\un$ in the strong operator topology, 
cf.\ Theorem \ref{RegThm}.     Hence
$\pi$ is nondegenerate on $\al K.$. Conversely, let $\pi$ be nondegenerate on
$\al K.$. If $\pi$ is not regular, \ie there is an $f\in X$ such that
$\ker\pi(\rlf)\not=0$, then we know that $\ker\pi(\rlf)$ reduces $\pi(\rsl)$,
\ie we can decompose $\pi=\pi_0\oplus\pi_R$ where $\pi_0(\rlf)=0$ and
$\ker\pi_R(\rlf) = 0$. Since $\rlf$ will occur in some products
of resolvents in $\al K.$ as in (i),  
$\pi_0$ when restricted to $\al K.$ has 
nonzero kernel. However $\al K.$ is simple, so $\pi_0(\al K.)=0$,
and this contradicts the assumption that $\pi$ is nondegenerate on
$\al K.$. Thus $\pi$ is regular.\chop
(iv) If $n=1$ then by (i) we see that if $\sigma(f,g)\not=0$ then  
$R(\lambda,f)R(\mu,g)\in\al K..$ Thus in the factor algebra
$\rsl\big/\al K.$, all products of noncommuting pairs in the generating set of resolvents
will be put to zero, and it is clear that only commuting products survive the
factoring. If $n>1$ then the products $R(\lambda,f)R(\mu,g)$ are generally not
in $\al K.$. Note that if $X=S\oplus S^\perp$ then $\al R.(S,\sigma)$ imbeds as
$\al R.(S,\sigma)\otimes\un$ of the subalgebra
$\al R.(S,\sigma)\otimes\al R.(S^\perp,\sigma)$ and so nonzero commutators of
elements of $\al R.(S,\sigma)$
are of the form $B\otimes\un$, which cannot be compact in $\pi_0$ for $B\not=0$.
\chop
(v) Let  $\pi$ be a regular representation and let
$\{I_n\}$ be a sequence in $\al K.$ as in (iii).
Thus, if $A\al K.=0$ then 
$0=\slim\limits_{n \to\infty}\pi(A I_n)=\pi(A)$. 
However $\pi$ is faithful by Theorem~\ref{RegFaith},
hence $A=0$. 
\chop
(vi) Since $\al K.$ is simple, it is clear that it is a minimal nonzero
ideal. Let $\al J.\subset \rsl$ be a nonzero closed two--sided ideal. 
Then so is 
$\al J.\cap \al K. = \al J.\cdot \al K.\not=\{0 \}$,
where the latter inequality follows from (v). So, since 
$\al J.\cap \al K.$  is a nonzero
ideal in $\al K.$ and $\al K.$ is simple, we get that 
$\al J.\cap \al K.=\al K.$ which is obviously in $\al J.$.

\subsection*{Proof of Proposition \ref{Dynamics}}

 Let  $U_0(t) = e^{itH_0}$, $t \in \R$ where $H_0 = P^2$; then since
 $H_0$ is quadratic in $P$, $\Ad U_0(t)$ induces a symplectic transformation
 on the resolvent algebra.  Thus  
 \[
 U_0(t) \, \pi_0(\rsl) \, U_0(t)^{-1} \subset \pi_0(\rsl) \, , \quad
 t \in \R \, .
 \]
 {}To prove that this inclusion still holds if 
 $U_0(t)$ is replaced by $U(t):= e^{itH}$ where $H = P^2 + V(Q)$, 
 we consider the cocycle $\Gamma_V(t) := U(t) \, U_0(t)^{-1}$, $t \in \R$. 
 It will suffice to show that the $\Gamma_V(t) - \un$ 
 are compact for all $V \in C_0(\R)$ since then
 $\Gamma_V(t) \in \pi_0(\rsl)$, $t \in \R$ and hence 
 \[
 U(t) \pi_0(\rsl) U(t)^{-1} = 
 \Gamma_V(t) \, U_0(t) \pi_0(\rsl) U_0(t)^{-1} \,
    \Gamma_V(t)^{-1} \subset \ \pi_0(\rsl) \, ,
 \] 
 by $\Gamma_V(t)^{-1} = \Gamma_V(t)^* \in \pi_0(\rsl)$. 
 
 We start with the Dyson perturbation  series of  $\Gamma_V(t)$ given by 
 \[
 \Gamma_V(t) = \sum_{n=0}^\infty \, i^n \int_0^t \! dt_n \int_0^{t_n}
 \! dt_{n-1} \cdots  \int_0^{t_2} \! dt_1 \, V_{t_1} \,
 V_{t_2} \cdots V_{t_n} \, ,
 \]
 where  $V_t :=  U_0(t) V(Q)  U_0(t)^{-1}$,
 cf.~\cite{BR1} Theorem~3.1.33.
 The integrals
 are defined in the strong operator topology. It is an immediate
 and well--known consequence of this formula that the cocycles
 $\Gamma_V(t)$, $t \in \R$ are continuous in $V$. More precisely,
putting $\| V \| = \| V(Q) \|$ 
  
 \[ \| \Gamma_{V_1}(t) - \Gamma_{V_2}(t) \| \leq
 \| V_1 - V_2 \| \ {e^{|t| (\| V_1 \| + \| V_2 \|)} - 1 \over 
 {\scriptstyle \| V_1 \| + \| V_2 \|} } \, . 
 \]
 So it suffices to prove compactness of $\Gamma_{V}(t) - \un$ 
 for a subspace of functions $V$
 which are dense in $C_0(\R)$, and we will take the space 
 $\big\{V\in\al S.(\R)\,\mid\,\int \! dx \, V(x) = 0\big\}$.
 {}For functions  in this space, the
 operators $\int_0^t ds \, V_s$ are  Hilbert--Schmidt. To 
 see this, consider
  their integral kernels in 
 momentum space, which are given by 
 \[ 
 \R \ni (u,v) \, \mapsto \, \Big(\int_0^t \! ds \, V_s \Big)(u,v) = 
 {i \over \sqrt{2\pi}} \, 
 {1 - e^{\, it (u^2 - v^2)} \over {(u^2 - v^2)}} \, 
 \wt{V}(u - v) \, ,
 \] 
 where $\wt{V}$ denotes the Fourier transform of $V$.
 Then  the square of the 
 Hilbert--Schmidt norm of the operator is:
 \[ 
 \Big\| \int_0^t ds \, V_s  \Big\|_2^2  = 
 \int \! du \! \int \! dv \, \Big| 
 \Big( \int_0^t \! ds \, V_s \Big) (u,v) \Big|^2 =
 |t| \int \! dw \, {| \wt{V}(w) |^2 \over {2 |w|}} 
 < \infty \, .
 \]
 The latter bound follows from the fact that $ \wt{V}$ is
 a test function which vanishes at the origin. Thus
 the strong operator continuous functions    
 \[
 \R^{n-1} \ni (t_2, \dots t_n) \, \mapsto \, 
 \int_0^{t_2} \! dt_1 \, V_{t_1} \,
 V_{t_2} \cdots V_{t_n} \,,\; n-1 \in \N
 \]
 have values in the Hilbert--Schmidt class and their Hilbert--Schmidt
 norms are bounded by 
 \[
 \Big\| \int_0^{t_2} \! dt_1 \, V_{t_1} \,
 V_{t_2} \cdots V_{t_n} \Big\|_2^2 \leq |t_2|  
 \int \! dw \, {| \wt{V}(w) |^2 \over {2 |w|}} \, \| V \|^{2n-2} \, . 
 \]
 In particular, these norms are uniformly bounded on compact subsets of 
 $\R^{n-1}$. But the integral of any  strong operator 
 continuous 
 Hilbert--Schmidt valued function with uniformly bounded 
 Hilbert--Schmidt norm is again in the Hilbert--Schmidt class. So 
 we conclude that each term in the above Dyson expansion
 is a Hilbert--Schmidt operator, except for the first term $\un$ 
 corresponding to $n=0$. Since the Dyson series converges absolutely in 
 norm, this shows that $\Gamma_V(t) - \un$ is a compact
 operator for the restricted class of potentials $V$. 
 The statement for arbitrary $V \in C_0(\R)$  then follows 
 from the continuity of $\Gamma_V(t)$ in $V$. 

We mention as an aside that the operators 
$\big( U(t) W U(t)^{-1} - U_0(t) W U_0(t)^{-1} \big)$,
$t \in \R$ are compact
for any bounded operator $W$ as a consequence of the
preceding result. Thus if a norm--closed and irreducible 
subalgebra ${\cal W}$ of the algebra of all bounded operators 
is to be stable under the action of the given family of 
dynamics it must contain the compact operators. 
The Weyl algebra, being simple, does not have this feature
and therefore does not admit interesting dynamics.

\subsection*{Proof of Proposition \ref{Affil}}

By the very definition of the resolvent algebra we have  
$(i \lambda \un - P)^{-1} \in \pi_0(\rsl)$, $\lambda \in \R \backslash
\{0\}$. As $ \pi_0(\rsl)$ is a C*--algebra, any continuous, 
asymptotically vanishing function of $P$ is therefore an element 
of $\pi_0(\rsl)$ as well, cf.\ Proposition
\ref{Czero}. In particular, the resolvent of
the free Hamiltonian $H_0 = P^2$ is contained in $\pi_0(\rsl)$.
Now for $H = P^2 + V(Q)$ we have
\[ 
(i \lambda \un - H)^{-1} = (i \lambda \un - H_0)^{-1}
+ (i \lambda \un - H)^{-1} \, V(Q) \, (i \lambda \un - H_0)^{-1} \, ,
\quad \lambda \in \R \backslash
\{0\} \, .
\] 
It follows from standard arguments, cf.\ \cite{AgCo, RS4},
that for the given family of potentials $V$ the operators 
$V(Q) \, (i \lambda \un - H_0)^{-1}$ are compact. Hence the 
resolvent of $H$ is contained in  $\pi_0(\rsl)$.

\subsection*{Proof of Proposition \ref{Notaffil}}

Let $U(n) := e^{\, i n H}$, $n \in \N$ and let $\Omega$
be any normalized vector in the underlying Hilbert
space. We define a corresponding sequence of states 
$\omega_n$ on $\rsl$, putting
\[
\omega_n (R) := \big( \Omega, U(n) \pi_0(R)  U(n)^{-1} \Omega \big) \, , \quad
R \in \rsl \, .
\] 
As the Hamiltonian $H = P^2 - Q^2$ is quadratic, one obtains 
by a routine calculation, $f_\pm \in \R$,  
\[ U(n) \big( i \lambda \un - f_+ \, (P + Q) - 
f_- \, (P - Q) \big)^{-1} U(n)^{-1} = 
\big( i \lambda \un - e^{2n} f_+ \, (P + Q) - e^{-2n} 
f_- \, (P - Q) \big)^{-1} \, .
\]
Hence, by the same reasoning as in the proof of 
Proposition \ref{NonregApprox}, one finds that  
\begin{eqnarray*}
&&\slim_{n\to\infty}\, U(n) \big( i \lambda \un - f_+ \, (P + Q) - 
f_- \, (P - Q) \big)^{-1} U(n)^{-1} = 
\cases{0 & if $f_+ \not=0$ \cr
{1 \over i \lambda} \, \un
& if $f_+ =0$ \cr} 
\end{eqnarray*}
It follows that the sequence of states $\omega_n$, $n \in \N$ converges 
pointwise on $\rsl$ and that its limit $\omega_\infty$ induces a
one--dimensional representation of $\rsl$. \chop
Assume now that there is some pseudo--resolvent
$R_\lambda \in \rsl$ such that $\pi_0(R_\lambda) = 
(i \lambda \un - H)^{-1}$, $\lambda \in \R \backslash 0$. 
By the preceding result and the resolvent equation for 
$R_\lambda$ we then have 
$\omega_\infty (R_\lambda) = (i\lambda - \nu)^{-1}$ for some 
$\nu \in \R \cup \{ \infty \}$. On the other hand it follows from
the definition of the states $\omega_n$ that 
$\omega_n(R_\lambda) = (\Omega, (i \lambda \un - H)^{-1} \, \Omega)$,
$n \in \N$. Hence we conclude that
\[
(\Omega, (i \lambda \un - H)^{-1} \, \Omega) =  (i\lambda - \nu)^{-1} \,
, \quad \lambda \in \R \backslash 0 \, .
\]
But this is impossible since $H$ has continuous spectrum as $-Q^2$ is a repulsive potential. 

\subsection*{Proof of Proposition \ref{InteractDynamics}}

(i) In this somewhat lengthy proof we
will state intermediate
 results in italics if they are of interest in their own right.
 We begin by gathering  notation and elementary facts. For
 $\Lambda \subset \Z\ni l$, let
\begin{eqnarray*} 
&& H_\Lambda^{(0)} := \sum_{l \in \Lambda} (P^2_l + Q^2_l) \, , \qquad
 U_\Lambda^{(0)} (t) := e^{ \, it H_\Lambda^{(0)}} \, , \ t \in \R \, , \\[1mm]
&& K_l^{(0)} := 
 {1 \over 2 } 
 \big( (P_l - P_{l + 1} )^2 + (Q_l - Q_{l + 1} )^2 \big) \, , \qquad
 V_l^{(0)} (t) := e^{ \, it K_l^{(0)}} \, , \ t \in \R \, .
\end{eqnarray*} 
Since $U_\Lambda^{(0)} (t)\in \pi_0( {\cal R}(X_{\Lambda}, \sigma))^{\prime \prime}$ 
and $\Ad U_\Lambda^{(0)} (t)$ induce symplectic transformations 
 on $ {\cal R}(X_{\Lambda}, \sigma)$, 
 we have for any $\Lambda_0 \subset \Lambda$ that  
 \[
 U_\Lambda^{(0)} (t) \, R_0  \, U_\Lambda^{(0)} (t)^{-1}
 = U_{\Lambda_0}^{(0)} (t) \, R_0  \, U_{\Lambda_0}^{(0)} (t)^{-1} 
 \in \pi_0({\cal R}(X_{\Lambda_0}, \sigma)) \, ,
 \qquad R_0 \in \pi_0({\cal R}(X_{\Lambda_0}, \sigma)) \, .
 \]
Furthermore, by
 \[
 2 (P_l^2 + Q_l^2 + P_{l+1}^2 +  Q_{l+1}^2) =
 ( (P_l - P_{l + 1} )^2 + (Q_l - Q_{l + 1} )^2 ) +
 ( (P_l + P_{l + 1} )^2 + (Q_l + Q_{l + 1} )^2 ) 
 \] 
 we have for $l, l+1 \in \Lambda$ that 
 \begin{equation}
 \label{ULambdaV}
U_\Lambda^{(0)} (t) \, R_l  \, U_\Lambda^{(0)} (t)^{-1}
 = V_l^{(0)} (t) \, R_l  \, V_l^{(0)} (t)^{-1} 
 \in \pi_0({\cal R}(Y_l, \sigma)) \, ,
 \quad R_l \in \pi_0({\cal R}(Y_l, \sigma)) \, ,
\end{equation}
 where 
 $Y_l := {{\rm Span}\{p_l - p_{l + 1},
 q_l -  q_{l + 1} \}}$.
  
 After these preparations we start our proof of (i).
 We intend to show that the cocycles 
 $\Gamma_\Lambda(t) := U_\Lambda(t) \, U_\Lambda^{(0)} (t)^{-1}$, $t \in \R$
 are in $\pi_0({\cal R}(X_{\Lambda}, \sigma))$. 
 As in the proof of
 Proposition~\ref{dynamics} consider
 the Dyson expansion of $\Gamma_\Lambda(t)$
 which now takes the form 
 \begin{eqnarray*}
 &&\Gamma_\Lambda (t) = \sum_{n=0}^\infty \, i^n \int_0^t \! dt_n \int_0^{t_n}
 \! dt_{n-1} \cdots  \int_0^{t_2} \! dt_1 \, V_{\Lambda, t_1} \,
 V_{\Lambda, t_2} \cdots V_{\Lambda, t_n} \,, \\[1mm] 
 V_{\Lambda, t} &:=& \sum_{l, l+1 \in \Lambda} V_{l, t} \, , \quad  
 V_{l, t}  := U_\Lambda^{(0)} (t) \, V(Q_l - Q_{l+1}) \,
 U_\Lambda^{(0)}(t)^{-1} =  V_l^{(0)} (t) \, V(Q_l - Q_{l+1}) \, V_l^{(0)}
 (t)^{-1} \, ,
 \end{eqnarray*}
 where the latter equality follows from equation~(\ref{ULambdaV})   via 
 $ V(Q_l - Q_{l+1}) \in \pi_0({\cal R}(Y_l, \sigma))$.
 We will show that each summand in 
 the Dyson expansion is in  
 $\pi_0({\cal R}(X_{\Lambda}, \sigma))$. 
 Consider the first non--trivial term. 
 {}For its building blocks we have: \\[1mm]

\noindent
 {\it The functions $\R \ni t \mapsto \int_0^t  ds \, V_{l,s}$ are
 H\"older continuous in the norm topology and their values are  
  in $\pi_0({\cal R}(Y_l, \sigma))$, $l \in \Z$. } \\[1mm]
 
 \noindent
 The H\"older continuity is obvious by
 $\| \int_0^{t_1}  ds \, V_{l,s} -  \int_0^{t_2}  ds \, V_{l,s} \|
 \leq |t_1 - t_2| \| V \|$.
 Since $V_{l,s} \in \pi_0({\cal R}(Y_l, \sigma))$ and 
 the integral is defined in the strong operator topology,
 it is also clear that 
 $\int_0^t  ds \, V_{l,s} \in \pi_0({\cal R}(Y_l, \sigma))^{\prime
 \prime}$. But $\pi_0({\cal R}(Y_l, \sigma))^{\prime \prime}$ 
 is a factor, so for the  second part of
 the statement it suffices to prove that 
 $\int_0^t  ds \, V_{l,s} \upharpoonright {\cal H}_0(Y_l) 
 \in \pi_0({\cal R}(Y_l, \sigma))  \upharpoonright {\cal H}_0(Y_l) $,
 where ${\cal H}_0(Y_l) = \overline{\pi_0({\cal R}(Y_l, \sigma)) \Omega_0}$.
 This will be done by proving that  
 $\int_0^t  ds \, V_{l,s} \upharpoonright {\cal H}_0(Y_l)$ 
 is  compact, 
 using the fact that
 $\pi_0 \upharpoonright 
 {\cal R}(Y_l, \sigma)$ on ${\cal H}_0(Y_l)$ is equivalent
 to the Schr\"odinger representation of ${\cal R}(Y_l, \sigma)$.
  
Define the canonical operators
 $Q = {1 \over \sqrt{2}} (Q_l - Q_{l+1}) \upharpoonright  {\cal
 H}_0(Y_l) $ and
 $P = {1 \over \sqrt{2}} (P_l - P_{l+1}) \upharpoonright  {\cal
 H}_0(Y_l)$,
 then 
 \[ 
 C := \int_0^t  ds \, V_{l,s} \upharpoonright {\cal H}_0(Y_l)
 = \int_0^t  ds \, e^{\, is (P^2 + Q^2) } \, V (\sqrt{2} \,  Q)
 \, e^{\,- is (P^2 + Q^2) } \, .
 \] 
Let $\Phi_n \in {\cal H}_0(Y_l)$ be the orthonormal
basis of eigenvectors of 
 $P^2 + Q^2$ corresponding to the eigenvalues $2n +1 $,  $n=0,\,1,\,2,\ldots$, 
then for $n\not= m$ we have
 \[ C_{m \,  n} := 
 \Big( \Phi_m , \int_0^t  ds \, e^{\, is (P^2 + Q^2) } \, V (\sqrt{2} \,  Q)
 \, e^{\,- is (P^2 + Q^2) } \Phi_n \Big) 
 = { e^{ 2it (m - n)} - 1 \over 2 i  (m - n)} \, 
 \big( \Phi_m,  V (\sqrt{2} \,  Q) \Phi_n  \big) \, ,
 \]
 where for $m = n$ the fraction has to be replaced by $t$. 
 We need to estimate the  
 matrix elements of the potential on the rhs, and will first consider
 potentials $V \in C_0(\R)$ with compact support.
 Then
 \[
 |\big( \Phi_m,  V (\sqrt{2} \,  Q) \Phi_n  \big)|
 \leq \|V(\sqrt{2} Q) \, e^{Q^2} \| \, \| e^{-{1 \over 2} Q^2} \Phi_m
 \| \, \| e^{-{1 \over 2} Q^2} \Phi_n \|\,.
 \]
{}From the standard representation of $\Phi_n$ in 
 configuration space by Hermite functions,
 \[
 x \mapsto \Phi_n(x) = (-1)^n \, (2^n n!)^{-1/2} \pi^{-1/4} \, 
 e^{x^2 /2} { d^n \over dx^n} e^{-x^2} \, ,
 \]
one gets by proceeding to Fourier transforms and
making use of Parseval's Theorem 
  \[ \| e^{-{1 \over 2} Q^2} \Phi_n \|^2 = {1 \over \sqrt{2} }
  \, { (2 n)! \over 2^{2n} (n!)^2 } \, .
 \]
Then the estimate (cf.~\cite{Rob})
\[
\sqrt{2\pi}\,n^{n+1/2}e^{-n+{(12n+1)}^{-1}}<
n!< \sqrt{2\pi}\,n^{n+1/2}e^{-n+{(12n)}^{-1}}\quad\hbox{for}\quad n\geq 1
\]
implies that    
\[
\| e^{-{1 \over 2} Q^2} \Phi_n \|^2\leq n^{-1/2}\quad\hbox{for}\quad n\geq 1 \, ,
\]
and for $n=0$ one has
${\| e^{-{1 \over 2} Q^2} \Phi_0 \|^2}\leq 1\,.$ 
 Thus one obtains
\[
 |  C_{m \, n} | \leq  
 \cases{
K / (m^{1/4} n^{1/4} \, |m-n| )   & \mbox{if $m \neq n$ and $m,\, n > 0$} \cr
  K / n^{5/4} & \mbox{if $m=0, n > 0$} \cr
K / m^{5/4} & \mbox{if $n=0, m > 0$} \cr
  |t| K / n^{1/2} & \mbox{if $m = n > 0$}  \cr
 |t| K & \mbox{if $m=n=0$} \, . \cr}
\]  
where $K:=\big\|V(\sqrt{2}\,Q)e^{Q^2}\big\|$.  

These bounds will enable us to show that $C$ is a compact operator. 
Note first that the operators $C_N:=C\cdot P_N,$ $N \in\N$ where $P_N$
is the projection onto ${\rm Span}\{\Phi_1,\ldots,\Phi_N\}$, are compact
as $C$ is bounded and $P_N$ is finite rank.
Hence to prove that $C$ is a compact operator, it suffices to show
that $\| C - C_N \| \rightarrow 0$ as $N \rightarrow \infty$. 
This can be accomplished by Schur's test according to which 
the norm of an operator $D$ satisfies the bound  
$\|D\| \leq \sqrt{ab}$ if there exist $a,b \in \R_+$ such that
$\sup_n \, \sum_m |D_{m \, n}| < a$ and $\sup_m \, \sum_n |D_{m \, n}|
< b$ where $D_{m \, n}$ denotes its matrix element w.r.t. a given orthonormal basis of
the Hilbert space on which it acts. 
Using the preceding bounds on the matrix elements
$ C_{m \, n}$ one can show by a routine computation that 
$\sup_{n \geq N} \, \sum_m |C_{m \, n}| \rightarrow 0$ as $N \rightarrow \infty$
and similarly $\sum_{n \geq N} |C_{m \, n}| \rightarrow 0$
uniformly in $m$ as $N \rightarrow \infty$.
It follows  that $\| C - C_N \| \rightarrow 0$ for $N \rightarrow \infty$. 
Thus 
 $C$ is compact for the restricted
 class of potentials $V$ and this result extends to arbitrary 
 potentials by norm continuity of $C$ in $V \in C_0(\R)$. 
 
 Since $Y_l \subset X_\Lambda$ for $l, l+1 \in \Lambda$ we 
 thus have shown that 
 $\int_0^t \! ds \, V_{\Lambda, s} = \int_0^t \! ds \, 
 \sum\limits_{l, l+1 \in \Lambda} V_{l, s} \in
 \pi_0({\cal R}(X_\Lambda , \sigma))$.
 {}For the proof that the remaining terms in the Dyson expansion are
 also contained in $\pi_0({\cal R}(X_\Lambda, \sigma))$ we make use of the
 following  fact. \\[1mm]
 {\it For a C*--algebra $ {\cal C} $ on a Hilbert space ${\cal H}$, 
 let $\R \ni s \mapsto A_s \in  {\cal C}$ be H\"older continuous 
 on compact subsets of $\R$ in the norm topology, 
 let $\R \ni s \mapsto B_s \in {\cal B}({\cal H})$ 
 be continuous in the strong operator topology, 
 and let $\int_0^t \! ds \,  B_s \in  {\cal C}$,
 $t \in \R$. Then $\R \ni t \mapsto C_t := \int_0^t \! ds \,  B_s A_s$ is
 H\"older continuous in the norm topology on compact subsets of
 $\R$ and $ C_t \in  {\cal C} $, $t \in \R$.} \\[1mm]
 \noindent We prove this. 
 By the assumptions, the operators 
 \[
 C_t^{(n)}  := \sum_{m=0}^{n-1} \int_{tm/n}^{t(m+1)/n} \! ds  
 \,  B_s \, A_{t(m+1/2) / n} 
 \]
 are in $ {\cal C} $ for any $n \in \N$. But
 \[
 \| C_t  - C_t^{(n)} \|  \leq {|t| \over n} \sup_{0 \leq s \leq t}  \| B_s  \|
 \sum_{m=0}^{n-1} 
 \sup_{mt/n\leq s\leq t(m+1)/n} \|A_s   - A_{t(m+1/2) / n} \| 
 \leq c \, {|t|^{1+h} \over n^h} \, ,  
 \]
 where in the second inequality 
 we used the assumption that $ \| A_{s_1} - A_{s_2} \| \leq c^\prime \, |s_1 -
 s_2|^h$ on compact subsets of $\R\,,$ and $c$ is a constant depending on
 $B_s\,,$ $c'$ and $h$. 
 Thus $C_t$ can be approximated in norm by elements of $ {\cal C} $ 
 and hence is in $\al C.$ as well.  
 {}For the H\"older--continuity, note that
 \[ \| C_{t_1} - C_{t_2} \| \leq | t_1 - t_2 | \sup_{t_1 \leq s \leq
 t_2} \| B_s A_s \| \leq c \, | t_1 - t_2 |  
 \] 
 on compact subsets of $\R$. 
 
 We now use the 
  preceding two results to prove by induction that all
 terms in the Dyson expansion are in $\pi_0({\cal R}(X_\Lambda
 , \sigma))$. For $n = 1$ we have already shown that 
 $t \mapsto S_{\, t}^{(1)} := \int_0^t \! ds \, V_{\Lambda, s} $ is H\"older continuous
 in norm (on compact subsets of $\R$)
 with values in $\pi_0({\cal R}(X_\Lambda, \sigma))$.
Assume that
 \[ 
 t \mapsto S_{\, t}^{(n)} := \int_0^t \! dt_n \int_0^{t_n}
 \! dt_{n-1} \cdots  \int_0^{t_2} \! dt_1 \, V_{\Lambda, t_1} \,
 V_{\Lambda, t_2} \cdots V_{\Lambda, t_n}
 \]
 has these properties too. For the inductive
 step, note that the H\"older continuity of $ S_{\, t}^{(n+1)}$
 follows from the estimate 
 \begin{eqnarray*} 
 \|  S_{\, t}^{(n+1)} -  S_{\, t^\prime}^{(n+1)} \|
 & \leq & \Big| 
 \int_{t^\prime}^t \! dt_{n+1} \, \Big| \int_0^{t_{n+1}}
 \! dt_{n}  \cdots  \Big| \int_0^{t_2} \! dt_1 \Big| \cdots \Big|\Big| \,  
 \sup_{t_1, \dots , t_n} \|  V_{\Lambda, t_1} \,
 V_{\Lambda, t_2} \cdots V_{\Lambda, t_{n+1}} \|  \\ 
 & \leq &  {\big|\epsilon(t) \, |t|^{n+1} - 
\epsilon(t^\prime) \, |t'|^{n+1} \big| \over (n + 1) \, !} \
 | \Lambda |^{n+1} \|V\|^{n+1} \, ,
 \end{eqnarray*}
 where $|\Lambda|$ is the number of points in $\Lambda$
and $\epsilon$ the sign--function.
 But $ S_{\, t}^{(n+1)} = \int_0^t dt_{n+1}  
 V_{\Lambda, t_{n+1}}  S_{\,t_{n+1} }^{(n)}$,
 so it follows from the induction basis and 
 hypothesis by an application of the preceding general
 result that $ S_{\, t}^{(n+1)} \in  \pi_0({\cal R}(X_\Lambda,
 \sigma))$, completing the induction. 
 Since the Dyson series converges absolutely in
 norm we conclude that also $\Gamma_\Lambda(t) \in \pi_0({\cal R}(X_\Lambda,
 \sigma))$, $t \in \R$. Since the adjoint action
 of $ U_\Lambda^{(0)} (t) $ leaves  $\pi_0({\cal R}(X_\Lambda,
 \sigma))$ invariant it is then clear that 
 \[
 U_\Lambda(t) \pi_0({\cal R}(X_\Lambda,
 \sigma)) U_\Lambda(t)^{-1} = 
 \Gamma_\Lambda(t) U_\Lambda^{(0)}(t) \pi_0({\cal R}(X_\Lambda,
 \sigma))  U_\Lambda^{(0)}(t)^{-1} \Gamma_\Lambda(t)^{-1}
 \subset \pi_0({\cal R}(X_\Lambda, \sigma)) \, .
 \]
\chop
(ii) We use a standard argument from the theory of spin systems \cite{BR1}
for this part.
By the net structure of $\rsl$ 
it suffices to prove the claim for the sets 
$\Lambda_0 := \{ l \in \N \mid |l| \leq n_0 \} \subset \Z$
and henceforth we fix  such a $\Lambda_0$, hence an $n_0,$
and let $R_0 \in \pi_0({\cal R}(X_{\Lambda_0}, \sigma))$.

To prove the claimed convergence, we start with
the Dyson  perturbation series for the adjoint action
of the cocycle  $\Gamma_\Lambda(t):=U_\Lambda(t) \, U_\Lambda^{(0)} (t)^{-1}$: 
\begin{equation} \label{PerturbSeries}
\Gamma_\Lambda(t) R_0 \Gamma_\Lambda(t)^{-1} =
R_0 + \sum_{n=1}^\infty i^n  \!
\int_0^t \! dt_1  \! \int_0^{t_1}  \! dt_2 \cdots  \! \int_0^{t_{n-1}}
\! \! dt_n \,
[V_{\Lambda, t_n}, [V_{\Lambda, t_{n-1}},  \cdots ,[V_{\Lambda, t_1},
R_0 ] \cdots ]\, ] \, .
\end{equation} 
Since 
$V_{l, \,t} \in  \pi_0({\cal R}(S_l , \sigma))$, $t \in \R$
it commutes with
$\pi_0({\cal R}(X_{\Lambda_n}, \sigma))$
if $l > n_0 + n$ or $l+1 < -n_0 - n$
where $\Lambda_n :=  \{ l \in \N \ | \ |l| \leq n_0 + n \}$. So
 by $R_0 \in \pi_0({\cal R}(X_{\Lambda_0}, \sigma))$ we have
\[
[V_{\Lambda, t_n}, [V_{\Lambda, t_{n-1}},  \cdots ,[V_{\Lambda, t_1},
R_0 ] \cdots ]\, ]=
[V_{\Lambda\cap\Lambda_n, t_n}, [V_{\Lambda\cap\Lambda_{n-1}, t_{n-1}},  \cdots ,
[V_{\Lambda\cap\Lambda_1, t_1},
R_0 ] \cdots ]\, ] 
\]
and hence 
\begin{eqnarray*}
&& \big\|[V_{\Lambda, t_n}, [V_{\Lambda, t_{n-1}},  \cdots ,[V_{\Lambda, t_1},
R_0 ] \cdots ]\, ]\big\|  \\[1mm]
&& 
\leq 2^n\big\|\sum_{l\in\Lambda\cap\Lambda_n}V_{l, 0}\big\|\cdots
\big\|\sum_{l^\prime \in\Lambda\cap\Lambda_1}V_{l^\prime, 0}\big\|\,\|R_0\|  
\leq 4^n(n_0 + 1) \cdots (n_0+ n) \| V \|^n \| R_0 \|\,.
\end{eqnarray*}
Moreover, if
$\Lambda, \Lambda^\prime$ are regions which both contain
$\Lambda_N$
for some $N \in \N$, then the first $N$ terms in the
respective Dyson series coincide and so
\[
\| 
\Gamma_\Lambda(t)R_0 \Gamma_\Lambda(t)^{-1} -
\Gamma_{\Lambda^\prime}(t) R_0 \Gamma_{\Lambda^\prime}(t)^{-1} \|
\leq \sum_{n=N+1}^\infty {|t|^n \over n \, !} 4^n (n_0 + 1) \cdots (n_0
+ n) \| V \|^n \| R_0 \| \, .
\]
The upper bound exists if $|t| < { 1 \over 4||V||\ }$ and it tends to $0$ 
as $N \rightarrow \infty$. Thus for 
 $R_0 \in \pi_0({\cal R}(X_{\Lambda_0}, \sigma))$
and sufficiently small $|t|$ the nets 
$\{ \Gamma_{\Lambda} (t) R_0 \Gamma_{\Lambda} (t)^{-1} \}_{\Lambda
\subset  \Z}$ converge (uniformly) as $\Lambda \nearrow \Z$. 
Since $U_{\Lambda}^{(0)} (t) R_0 U_{\Lambda}^{(0)} (t)^{-1} =
U_{\Lambda_0}^{(0)} (t) R_0 U_{\Lambda_0}^{(0)} (t)^{-1}
\in \pi_0({\cal R}(X_{\Lambda_0}))$  by the remarks at the
beginning of this proof and hence 
\[
U_{\Lambda} (t) R_0 U_{\Lambda} (t)^{-1} 
=  \Gamma_{\Lambda} (t) \, U^{(0)}_{\Lambda_0} (t)  R_0
 U^{(0)}_{\Lambda_0} (t)^{-1}  \, \Gamma_{\Lambda} (t)^{-1}
\]
we conclude that the nets
$\{ U_{\Lambda} (t) R_0 U_{\Lambda} (t)^{-1} 
 \}_{\Lambda
\subset  \Z}$ also converge in norm as $\Lambda \nearrow \Z$. However
$\bigcup\limits_{\Lambda_0 \subset \Z} \, \pi_0({\cal R}(X_{\Lambda_0}, \sigma))$  
is norm dense in $\pi_0(\rsl)$ and the adjoint action of unitary operators
is norm continuous, so the existence of the norm limits
\begin{equation} \label{FullDynamics}
\beta_t(R) := \mbox{n--} \! \! \lim_{\Lambda \nearrow \Z} 
 U_{\Lambda} (t) \pi_0(R) U_{\Lambda} (t)^{-1} \, , \quad
R \in \rsl
\end{equation}
follows if  $|t| < { 1 \over 4||V||\ }$. 
Moreover, $\beta_t(\pi_0({\cal R}(X, \sigma)) 
\subset  \pi_0({\cal R}(X, \sigma))$, $|t| < { 1 \over 4||V||\ }$   
by part (i) above.
By the group property 
$ U_{\Lambda} (s +t) =  U_{\Lambda} (s)  U_{\Lambda} (t)$
one finds by repeated application of 
the preceding two results that these statements hold for
arbitrary $t \in \R$  which proves this part.
\chop 
(iii)
Recalling that the representation $\pi_0$ of $\rsl$
is faithful, one can 
define an automorphic action of the group $\R$ on  $\rsl$
induced by the set of Hamiltonians
$\{H_\Lambda\}_{\Lambda \in \Z}$,
putting for $t \in \R$
\[ \alpha_t(R) := \pi_0^{-1} \big( \beta_t(\pi_0(R))\big) \, ,
\quad R \in \rsl \, .
\]
This completes the proof of the proposition.

\subsection*{Proof of Proposition  \ref{ContinuousAction}}

Let $\Lambda_0 = \{ l \in \N \ | \ |l| \leq n_0 \}\subset\Z$
 and let 
$R_0 \in  \pi_0({\cal R}(X_{\Lambda_0}, \sigma))$. 
In the preceding proof of Proposition \ref{InteractDynamics}
we established the existence
of the norm limits
\[
\gamma_t(R_0) := \mbox{n--}\! \! \lim_{\Lambda \nearrow \Z}
\Gamma_\Lambda(t) R_0 \Gamma_\Lambda(t)^{-1} \, , \quad 
R_0 \in  \pi_0({\cal R}(X_{\Lambda_0}, \sigma)) \, .
\]
{}From the  expansion (\ref{PerturbSeries}) and
the remarks subsequent to it  we also obtain for 
$|t|, |t^\prime| < { 1 \over 4 \| V \| }$ the uniform bound
for $\Lambda \subset \Z$ 
\[ 
\| \Gamma_\Lambda(t) R_0 \Gamma_\Lambda(t)^{-1}
- \Gamma_\Lambda(t^\prime) R_0 \Gamma_\Lambda(t^\prime)^{-1} \|
\leq 
\sum_{n=1}^\infty 
{\big|\epsilon(t) \, |t|^n - \epsilon(t^\prime) \, |t'|^n\big| 
\over n \, !} 4^n (n_0 + 1) \cdots (n_0 + n) \| V \|^n \| R_0 \| \, .
\]
Combining these results it follows that $\gamma_t$ 
acts norm continuously on the elements of 
$ \pi_0({\cal R}(X_{\Lambda_0}, \sigma))$ 
for $|t| <  { 1 \over 4 \| V \| }$,  
and this continuity 
property extends to all of $\pi_0({\cal R}(X, \sigma))$
since $\bigcup\limits_{\Lambda_0 \subset \Z} \, \pi_0({\cal R}(X_{\Lambda_0}, \sigma))$  
is norm dense in $\pi_0(\rsl)$. Next, putting
\[ 
\beta_{\, t}^{(0)} (R) := \mbox{n--}\! \! \lim_{\Lambda \nearrow \Z}
\, 
U_\Lambda^{(0)}(t) \, R \, U_\Lambda^{(0)}(t)^{-1} \, , \quad R \in 
 \pi_0({\cal R}(X, \sigma)) 
\] 
we infer from the remarks at the beginning of the 
proof of Proposition \ref{InteractDynamics} that 
\[
\beta_{\, t}^{(0)} \upharpoonright
\pi_0({\cal R}(X_{\Lambda_0}, \sigma)) = 
({\Ad}U_{\Lambda_0}^{(0)}(t)) \upharpoonright
\pi_0({\cal R}(X_{\Lambda_0}, \sigma)) \subset 
\pi_0({\cal R}(X_{\Lambda_0}, \sigma)) \, .
\] 
It follows that $\beta_{\, t}^{(0)}$, $t \in  \R$  acts 
norm continuously 
on the (up to  multiplicity) compact
operators in $\pi_0({\cal K}_{\Lambda_0}) \subset 
\pi_0({\cal R}(X_{\Lambda_0}, \sigma)) $, $\Lambda_0 \subset \Z$
and hence also on their inductive limit $\pi_0({\cal K})$. \chop  
Now for the full time evolution (\ref{FullDynamics}) we have
$\beta_t = \gamma_t \circ \beta_{\, t}^{(0)}$, $t \in \R$.
So summarizing the   preceding results 
we conclude that $\beta_t$ acts pointwise norm continuously
on  $\pi_0({\cal K})$ for small $t$. 
In view of $\|\beta_t(\pi_0(K)) - \beta_{t^\prime}(\pi_0(K))\|
= \|\beta_{t - t^\prime} (\pi_0(K)) - \pi_0(K) \|$ this statement 
extends to arbitrary $t \in \R$ and it is then clear that
$\alpha_t$ acts pointwise norm continuously 
on  ${\cal K}$, $t \in \R$, proving the statement. 

\subsection*{Proof of Lemma \ref{RegularState}}

In addition to the operators $H_n$, $\widetilde{H}_n$, $n \in \N$ 
introduced in the main
text we will consider here also the
operators  $H_{n \backslash m}  := H_{\Lambda_{n \backslash m}}$
corresponding to the sets
$\Lambda_{n \backslash m} := \{ l \in \Z \, | \, m < |l| \leq n \}$
and their  renormalized versions
$\widetilde{ H}_{n \backslash m} =
 H_{n \backslash m} - E_{n \backslash m} \un$, where
$ E_{n \backslash m}$ is the smallest eigenvalue of 
$H_{\Lambda_{n \backslash m}}$. 
Since $\widetilde{H}_m$ and  $\widetilde{H}_{n \backslash m}$ commute and
the potential $V$ is bounded, the domains of these operators
are related by ${\cal D}(\widetilde{H}_n) = {\cal D}(\widetilde{H}_m)\cap {\cal
  D}(\widetilde{H}_{n \backslash m})$, and on the latter domain
we have the operator equality 
\[
\widetilde{H}_n = \widetilde{H}_m + 
\widetilde{H}_{n \backslash m} + V(Q_{-m -1} -Q_{-m}) + 
V(Q_m - Q_{m +1}) + \big(E_m + E_{n \backslash m} - E_n \big) \un \, .
\] 
Let $\Omega$ be a normalised joint eigenvector for $ \widetilde{H}_m$ and
$\widetilde{ H}_{n \backslash m}$ for the eigenvalue zero, 
then ${(\Omega,\widetilde{H}_n\Omega)}\geq 0$ implies via the last equation
that:
\[
\big(E_m + E_{n \backslash m} - E_n \big) \geq - 2 \, \| V \| \, .
\] 
Consequently $ \widetilde{H}_n \geq \widetilde{H}_m - 4 \, \| V \| \un $
and hence for their resolvents we have:
\[
\big((\mu +4\|v\|)\un+\wt{H}_{n}\big)^{-1}\leq
\big(\mu \un+\wt{H}_{m}\big)^{-1} \leq \mu^{-1} \un \quad
\mbox{for all} \,  \mu>0;\;
 m < n, \ m,n \in \N \, .
\]
Recalling that $\wt{H}_{m}$ is affiliated with 
${\cal K}_{\Lambda_m} \subset {\cal R}({\Lambda_m}, \sigma)$, let 
$\widetilde{R}_m(\mu) \in {\cal K}_{\Lambda_m}$, $\mu > 0$ be the 
corresponding pseudo resolvent,  
\ie $\pi_0(\widetilde{R}_m(\mu)) = \big(\mu \un+\wt{H}_{m}\big)^{-1}$, 
$m \in \N$. 
Since $\omega_\infty$ is a w*--limit point of $\{\omega_n\}_{n\in\N}$
there is a subsequence  $\{\omega_{n_k} \}_{k \in \N}$ such that
\[ 
\omega_\infty(\widetilde{R}_m(\mu)) = \lim_{k \rightarrow \infty}
 \big( \Omega_{n_k}, \pi_0(\widetilde{R}_m(\mu))
\Omega_{n_k} \big)
= \lim_{k \rightarrow \infty}
 \big( \Omega_{n_k},  \big(\mu \un+\wt{H}_{m}\big)^{-1} 
\Omega_{n_k} \big)
\] 
for fixed $m$ and $\mu$.
The preceding operator inequality and the fact that 
$\wt{H}_{n_k} \Omega_{n_k} = 0$, $k \in \N$ imply 
\[ (\mu + 4 \| V \|)^{-1} =
\lim_{k \rightarrow \infty}
 \big( \Omega_{n_k},  \big((\mu + 4 \| V \|) \un+\wt{H}_{n_k}\big)^{-1} 
\Omega_{n_k} \big)
\leq 
\lim_{k \rightarrow \infty}
 \big( \Omega_{n_k},  \big(\mu \un+\wt{H}_{m}\big)^{-1} 
\Omega_{n_k} \big) \leq \mu^{-1} \, .
\] 
Hence $ (\mu + 4 \| V \|)^{-1} \leq \omega_\infty(\widetilde{R}_m(\mu)) \leq
\mu^{-1}$, $\mu > 0$ and so $\lim\limits_{\mu \to\infty}
\, \omega_\infty(\mu \, \widetilde{R}_m(\mu)) = 1$. 
As $\| \mu \, \widetilde{R}_m(\mu) \|
\leq 1 $ and $\widetilde{R}_m(\mu) \in {\cal K}_{\Lambda_m}$ it follows that   
$\big\|\omega_\infty\restriction\al K._{\Lambda_m}\big\|=1$,  $m \in \N$.
Thus
the GNS representation of $\omega_\infty$ is nondegenerate on $\al K._{\Lambda_m}$
and hence by Theorem~\ref{CompactId}\,(iii), $\omega_\infty$ is regular on 
$\al R.(X_{\Lambda_m},\sigma)$. This holds for all 
$m \in \N$ hence $\omega_\infty$ is regular.

\subsection*{Proof of Proposition \ref{RDirac}}

(i) First, observe that if $\omega\in\ot S._D$ then $(i\lambda\rlf-\un)
\in {\cal N}_\omega\subset\ker\omega$ for $f\in C$ so putting $\lambda=1$ 
we get that $\omega(R(1,f))=-i$. Next, assume that 
$\omega(R(1,f))=-i$ for all $f\in C$. Then, using
$\lambda\rlf=R(1,\f 1,\lambda.f)$ we get that 
${(i\lambda\rlf-\un)}\in\ker\omega$ for $\lambda \in \R$.
Now by Eqs.~(\ref{Rinvol}) and (\ref{Resolv})
\[
\big(i\lambda\rlf-\un\big)^*\big(i\lambda\rlf-\un\big)
= \f i\lambda,2.\,\big(R(-\lambda,f)-\rlf\big)+\hbox{$\un$}
     \]
and hence 
$\omega\left(\big(i\lambda\rlf-\un\big)^*\big(i\lambda\rlf-\un\big)\right)
=0$ for $f\in C$, \ie $\omega\in\ot S._D$. Note that as  
$\big(i\lambda\rlf-\un\big)$ is a normal operator one also has 
 $\omega\left(\big(i\lambda\rlf-\un\big) \big(i\lambda\rlf-\un\big)^* \right)
=0$. 

(ii)
It suffices to prove that if $\sigma(g,C)\not=0$ for $g\in X,$ 
then $\omega(R(\mu,g))=0$ for all $\mu \in\R\backslash 0$,
since then $\pi_\omega(R(\lambda,g))=0$ by  Theorem~\ref{Ideals1}(iv).
Let $f\in C$ and $\omega\in\ot S._D$ then by 
$(i\lambda\rlf-\un) \in {\cal N}_\omega\cap{\cal N}_\omega^*$
and Eq.\ (\ref{Rccr}) we have
\[
0 = \omega\left(\big[(i \lambda \rlf - 1),\,R(\mu,g)\big]\right) 
= - \omega\big(\lambda \, \sigma(f,g)\,\rlf\,R(\mu,g)^2\rlf\big) 
= {\textstyle {1 \over \lambda}} \, \sigma(f,g)\,
\omega\big(R(\mu,g)^2\big) 
\]
so $\omega\big(R(\mu,g)^2\big)=0$. 
Now by Eq.\ (\ref{Resolv}) and the continuity properties of the 
resolvents we have
\[ i {\textstyle {d \over d \mu }}  R(\mu, g) = 
i \, \mbox{n--}\!\!\lim_{\lambda \rightarrow \mu} \, (\mu - \lambda)^{-1} \,
(R(\mu, g) - R(\lambda, g)) =
\mbox{n--}\!\!\lim_{\lambda \rightarrow \mu}
 \, R(\mu, g)  R(\lambda, g) = R(\mu, g)^2 \, .
\]
Hence $i \, {\textstyle {d \over d \mu }} \omega( R(\mu, g))
= \omega\big(R(\mu,g)^2\big) = 0$ 
and consequently $\mu \mapsto \omega( R(\mu, g)) = \mbox{const}$. But
$|\omega(R(\mu, g))| \leq 
\| R(\mu, g) \| \leq {1 \over |\mu|}$ and thus 
$\omega( R(\mu, g)) = 0$.

(iii) {}From part (ii) note that if $\sigma(C,C)\not=0,$ then for 
$f,\,g\in C$ with $\sigma(f,g)\not=0$ we have that ${\omega\big(i\lambda\rlf-\un\big)}
=0$ for all $\lambda$ implies that $\omega(R(\mu,g))=0,$ which contradicts 
with the requirement that  ${\omega\big(i\lambda R(\mu,g)-\un\big)}=0$.
Thus $\sigma(C,C)\not=0$ implies that $\ot S._D=\emptyset$.

{}For the converse let  $\sigma(C,C)=0$, then
$\al R.(C):={\rm C}^*\big\{\rlf\,\mid\,\lambda\in\R\backslash 0,\;
f\in C\} + \C \un$ is a unital commutative C*--algebra.
It is easily checked that the linear map 
from all polynomials in the resolvents in $\al R.(C)$
to the complex numbers given by  
\[
\omega(R(\lambda_1, f_1) \cdots R(\lambda_n, f_n))
:= \prod_{k=1}^n (1/i \lambda_k) \, , \quad 
\lambda_1 \dots \lambda_n \in \R\backslash 0 \, , \
f_1 \dots f_n \in C  
\]
is a *--homomorphism. Hence it extends to 
a character on $ \al R.(C)$, \ie a pure state,
and then, by the Hahn--Banach theorem, to a state
on $\rsl$. By its very construction,
$\omega\in\ot S._D$. 

\subsection*{Proof of Proposition \ref{Diracfacts}}

The proofs of the facts listed here already appeared elsewhere~\cite{GrLl},
but we recall them here for completeness.

(i) Recall that $\al N.:=[\rsl\al C.]
=\cap\; \{{\cal N}_{\omega}\mid\omega\in{\ot S.}_D \}$
so since ${\cal N}_{\omega}\subset{\rm Ker}\,\omega$ for all $\omega$
it is clear that  $\al D.=\al N.\cap\al N.^*\subset
\bigcap\, \{ {\rm Ker}\,\omega\mid \omega\in 
  {\got S}_{{D}} \}$. Since $\al D.$ is the intersection of a closed left
  ideal with its adjoint, we see that $\al D.$ is a C*--algebra.
  Next, let $\al A.$ be any C*--algebra in 
  $\bigcap\, \{ {\rm Ker}\,\omega\mid \omega\in 
  {\got S}_{{D}} \}$. Since $\al A.$ is a C*--algebra,
  $A\in\al A.$ implies that $AA^*\in\al A.\ni A^*A$ hence
  $A\in {\cal N}_\omega\cap {\cal N}_\omega^*$ for all $\omega\in 
  {\got S}_{{D}}$. Thus $\al A.\subseteq
  \bigcap\{ {\cal N}_\omega\cap {\cal N}_\omega^* \mid \omega\in 
  {\got S}_{{D}} \}=\al D.$ and so $\al D.$ is the maximal
  C*--algebra annihilated by all Dirac states.\chop
 To see that $\al D.$ is hereditary, use
Theorem~3.2.1 in \cite{Mur} and the fact that 
$\al N.=[ {\rsl}\al C.]$ is a closed left ideal of $\rsl$.

(ii) Since $\al D.$ is a two--sided ideal for 
the relative multiplier algebra $M_{\rsl}(\al D.)$ 
of  $\al D.$ in $\rsl$ it is
obvious that $M_{\rsl}(\al D.)\subset\WD$. Conversely, consider
$B\in\WD$, then for any $D\in\al D.$, we have
$BD=DB + D'\in\al N.$ with $D'$ some element of $\al D.$,
where we used $\rsl \al D.=\rsl (\al N.\cap\al N.^*)\subset\al N.$.
Similarly we see that $DB\in\al N.^*$. But then
$\al N.\ni BD=DB+D'\in\al N.^*$, so $BD\in\al N.\cap\al N.^*=\al D.$.
Likewise $DB\in\al D.$ and so $B\in M_{\rsl}(\al D.)$.

(iii) Since $\al C.\subset\al D.$ we see from the definition of $\WD$
that $F\in\WD$ implies that $[F,\,\al C.]\subset\al D.$.
Conversely, let $[F,\,\al C.]\subset\al D.$
for some $F\in\rsl$.
Now $F[\rsl \al C.]\subset[\rsl \al C.]$ and $F[\al C. \rsl]=
[F\al C. \rsl]\subset
[(\al C. \rsl +\al D.)\rsl]\subset[\al C. \rsl]$ because
$\al C.F+\al D.\subset\al C.F+ \al N.^* \subset[\al C. \rsl]$.
Thus $F\al D.=F\big([\rsl \al C.]\cap[\al C. \rsl]\big)
\subset[\rsl \al C.]\cap[\al C. \rsl]=\al D.$.
Similarly $\al D.F\subset\al D.$, and thus by (ii) we see
$F\in\WD$.

(iv) $\al D.\subset\WD$ so by (i) it is the unique maximal C*--algebra
annihilated by all the states $\omega\in\ot S._D(\WD)=
\ot S._D\restriction\WD$ (since $\al C.\subset\WD$).
Thus $\al D.=[\al OC.]\cap[\al CO.]$. But $\al C.\subset\al D.$,
so by (ii) $[\al OC.]\subset\al D.\subset[\al OC.]$ and so
$\al D.=[\al OC.]=[\al CO.]$.

\section{Appendix: Symplectic Spaces}

We collect here some basic facts for symplectic spaces 
which are required for the
preceding proofs. In this section $X$ will be a real linear space with a 
nondegenerate symplectic form $\sigma:X\times X\to\R$, and for any
subspace $S\subset X$ its symplectic complement
will be denoted by $S^\perp:={\set f\in X,\sigma(f,S)=0.}$.
By $X=S_1\oplus S_2\oplus\cdots\oplus S_n$ we will mean that all
$S_i$ are nondegenerate and $S_i\subset S_j^\perp$ if $i\not=j,$
and each $f\in X$ has a unique decomposition
$f=f_1+f_2+\cdots+f_n$ such that $f_i\in S_i$ for all $i$.
\begin{lem}
\label{SympFacts}
\begin{itemize}
\item[(i)] If $X$ is countably dimensional, then it has a symplectic basis, \ie 
a basis ${\big\{q_1,\,p_1;\,q_2,\,p_2;\ldots\big\}}$ such that
$\sigma(p_i,q_j)=\delta_{ij}$ and $0=\sigma(q_i,q_j)=\sigma(p_i,p_j)$ for all $i,\,j$.
\item[(ii)] {}For any symplectic space $X$ we have that if $S$ is a nondegenerate
finite--dimensional subspace, then $X=S\oplus S^\perp$
\item[(iii)] {}For any symplectic space $X$ and a finite linearly independent subset
${\big\{q_1,\,q_2,\ldots,\,q_k\big\}}\subset X$ such that $\sigma(q_i,q_j)=0$ 
for all $i,\,j,$ there is a set
${\big\{p_1,\,p_2,\ldots,\,p_k\big\}}\subset X$ such that 
$B:={\big\{q_1,\,p_1;\,q_2,\,p_2;\ldots;\,q_k,\,p_k\big\}}$
is a symplectic basis for ${\rm Span}(B)$.
\end{itemize}
\end{lem}
\begin{beweis}
(i)  Let $(e_n)_{n \in \N}$ be a linear basis of $X$. 
We construct the basis elements $p_n, q_n$ inductively as follows. 
If $p_1,\ldots, p_k$ and $q_1, \ldots, q_k$ are already chosen, pick a minimal
$m$ with $e_m \not\in {\rm Span}\{p_1,\ldots, p_k, q_1,\ldots, q_k\}$ and put 
$$ p_{k+1} := e_m - \sum_{i=1}^k\big(\sigma(e_m, q_i)p_i +\sigma(p_i, e_m) q_i\big) $$
to ensure that this element is $\sigma$-orthogonal to all previous ones. 
Then pick $l$ minimal, such that $\sigma(p_{k+1}, e_l) \not=0$, put
$$ \tilde q_{k+1} := e_l - \sum_{i=1}^k\big( \sigma(e_l, q_i)p_i +\sigma(p_i, e_l) q_i\big) $$
and pick $q_{k+1} \in \R \tilde q_{k+1}$ with
$\sigma(p_{k+1}, q_{k+1}) = 1$.  
This process can be repeated \textit{ad infinitum} 
and produces the required 
basis of $X$ because for each $k$, the span of 
$\big\{p_1,\ldots, p_k,\, q_1,\ldots, q_k\big\}$ contains at least 
$\{e_1,\ldots, e_k\}$. \chop
(ii)
Since $S$ is finite dimensional and nondegenerate, we can choose 
by (i) a symplectic basis
${\big\{q_1,\,p_1;\,q_2,\,p_2;\ldots;\,q_k,\,p_k\big\}}$ for it.
Given any $v\in X$ then
\[
v_S:=\sum_{i=1}^k\big(\sigma(v,q_i)\,p_i +\sigma(p_i,v)\,q_i\big)
       \in S
\]
and $v-v_S\in S^\perp$, \ie $\sigma(v-v_S,S)=0$.
Thus $X={\rm Span}\{S\cup S^\perp\}$, and as $\sigma$
is nondegenerate $S\cap S^\perp=\{0\}$. Moreover, if
$0=v+w$ where $v\in S$ and $w\in S^\perp$, then 
$v=-w\in S\cap S^\perp=\{0\}$, and hence any decomposition
of an $x\in X$ as $x=x_1+x_2$ where $x_1\in S,$ $x_2\in S^\perp$
is unique. Thus $X=S\oplus S^\perp$. \chop
(iii) 
We first find via the method of part (i),
symplectic pairs ${\big\{\wt{q}_1,r_1;\ldots;\wt{q}_k,\,r_k\big\}}\subset X$
such that the nondegenerate subspaces
$S_j:= {{\rm Span}\big\{\wt{q}_1,r_1;\ldots;\wt{q}_j,\,r_j\big\}}\supset\{q_1,\ldots, q_j\}$
but $q_{j+1}\not\in S_j$.
We construct the basis elements $\wt{q}_i,\, r_i$ inductively as follows.
If $r_1,\ldots, r_j$ and $\wt{q}_1, \ldots, \wt{q}_j$ are already chosen,
put
$$ \wt{q}_{j+1} := q_{j+1} - \sum_{i=1}^k\big(\sigma(q_{j+1}, \wt{q}_i)r_i +
\sigma(r_i, q_{j+1}) \wt{q}_i\big) $$
to ensure that $\wt{q}_{j+1}\in S_j^\perp$. By (ii), $X=S_j\oplus S_j^\perp$
hence $S_j^\perp$ is nondegenerate, so there is an element $r_{j+1}\in S_j^\perp$
such that $\sigma(r_{j+1},\wt{q}_{j+1})=1$.
It follows that $q_{j+2}\not\in S_{j+1}$ and that
$\{q_1,\ldots, q_{j+1}\}\subset S_{j+1}$.
This process can be repeated to produce the required symplectic bases.
Next, we want to show that in $S_k$ we can choose 
${\big\{p_1,\,p_2,\ldots,\,p_k\big\}}$ such that 
${\big\{q_1,\,p_1;\,q_2,\,p_2;\ldots;\,q_k,\,p_k\big\}}$
is a symplectic basis for $S_k$.
Now $\{q_1,\ldots, q_k\}\subset\{q_1,\ldots, q_k\}^\perp$
where henceforth the symplectic complements are all taken in $S_k$.
We claim that the containment ${\{q_2,\ldots, q_k\}^\perp}\supset
{\{q_1,\,q_2,\ldots, q_k\}^\perp}$ is proper. The map
$\varphi:S_k\to S_k^*$ by $\varphi_x(y):=\sigma(x,y)$ is a linear isomorphism
by nondegeneracy of $\sigma$. Then for any set $R\subset S_k$ we have 
$\varphi\big(R^\perp\big)=R^0$
\ie the annihilator of $R$ in $S_k^*$, hence 
${\rm dim}(R^\perp)={\rm dim}(R^0)=2k-{\rm dim}\big({\rm Span}(R)\big)$.
Thus ${\rm dim}{\{q_1,\ldots, q_j\}^\perp}=2k-j$ from which the claim follows.
Thus there is an $r\in{\{q_2,\ldots, q_k\}^\perp}\backslash
{\{q_1,\,q_2,\ldots, q_k\}^\perp}$ such that $\sigma(r,q_1)\not=0$.
In particular, let $p_1$ be that multiple of $r$ such that
$\sigma(p_1,q_1)=1$. Let $T_1:={\rm Span}\{q_1,p_1\}$
then ${\{q_2,\ldots, q_k\}}\subset T_1^\perp,$ and by (ii)
we have $S_k=T_1\oplus T_1^\perp$ where $T_1^\perp$ is nondegenerate.
 Thus we can now
repeat this procedure in $T_1^\perp$ starting from $q_2$ to obtain
$p_2$. This procedure will exhaust $S_k$ to produce the 
desired symplectic 
basis ${\big\{q_1,\,p_1;\,q_2,\,p_2;\ldots;\,q_k,\,p_k\big\}}.$ 
\end{beweis}

\section*{Acknowledgements.}

We wish to thank Professor H. Araki for his careful reading of the manuscript,
and suggestions which improved the arguments.
DB wishes to thank the Department of Mathematics of the 
University of New South Wales and HG wishes to thank the 
Institute for Theoretical Physics of the University of G\"ottingen
for hospitality and financial support which facilitated this research.
The work was also supported in part by the UNSW Faculty Research Grant Program.

\bigskip

\providecommand{\bysame}{\leavevmode\hbox to3em{\hrulefill}\thinspace}

\end{document}